\documentclass[11pt, a4paper, abstract=yes]{scrartcl}

\usepackage[utf8]{inputenc}
\usepackage[T1]{fontenc}
\usepackage{graphicx}
\usepackage{longtable}
\usepackage{wrapfig}
\usepackage{rotating}
\usepackage[normalem]{ulem}
\usepackage{amsmath}
\usepackage{amssymb}
\usepackage{hyperref}
\usepackage{amsthm,amssymb,amsmath}
\usepackage{hyperref}
\usepackage{graphicx}
\usepackage[T1]{fontenc}\usepackage{lmodern}
\theoremstyle{definition}
\newtheorem{definition}{Definition}
\newtheorem{mexample}[definition]{Example}
\theoremstyle{remark}
\newtheorem{remark}[definition]{Remark}
\theoremstyle{plain}
\newtheorem{lemma}[definition]{Lemma}

\newtheorem{proposition}[definition]{Proposition}
\newtheorem{theorem}[definition]{Theorem}
\newtheorem{corollary}[definition]{Corollary}

\newtheorem{conjecture}[definition]{Conjecture}

\newcommand{\vek}[1]{\boldsymbol{#1}}

\author{Jan Snellman\thanks{Department of Mathematics, Linköping University, 581 83 Linköping, Sweden. \texttt{jan.snellman@liu.se}. ORCID: \href{https://orcid.org/0009-0002-6676-5068}{0009-0002-6676-5068}}}
\usepackage[T1]{fontenc}
\usepackage{orcidlink}
\usepackage{booktabs}
\usepackage{etoolbox}
\AtEndEnvironment{mexample}{\hfill$\Diamond$}
\AtEndEnvironment{remark}{\hfill$\triangleright$}
\AtEndEnvironment{definition}{\hfill$\triangleleft$}
\AtEndEnvironment{observation}{\hfill$\triangledown$}
\theoremstyle{remark}
\newtheorem{observation}[definition]{Observation}
\theoremstyle{plain}
\hypersetup{hidelinks}
\newcommand{\Koz}{\operatorname{Koz}}
\newcommand{\shift}{\operatorname{shift}}
\newcommand{\conv}{\operatorname{conv}}
\newcommand{\indeg}{\operatorname{indeg}}
\newcommand{\tp}{\operatorname{top}}
\newcommand{\rasimp}{\angle}
\newcommand{\RA}{\operatorname{RA}}
\newcommand{\kozvec}{\vek\kappa}

\newcommand{\FF}{\mathrm{FF}}
\newcommand{\NV}{\operatorname{NV}}
\newcommand{\vol}{\operatorname{vol}}
\newcommand{\relint}{\operatorname{relint}}

\newcommand{\asc}{\operatorname{asc}}

\newcommand{\fl}{\operatorname{fl}}
\usepackage{tikz}
\usepackage{marginnote}
\definecolor{leanfull}{HTML}{2E7D32}
\definecolor{leancond}{HTML}{1565C0}
\definecolor{leannone}{HTML}{8A8A8A}
\newcommand{\leandiscof}[1]{\tikz[baseline=-0.6ex]\fill[#1] (0,0) circle (0.6ex);}
\newcommand{\leanfull}{\leandiscof{leanfull}}
\newcommand{\leancond}{\tikz[baseline=-0.6ex]{\fill[leancond] (0,0) circle (0.38ex);\draw[leancond,line width=0.5pt] (0,0) circle (0.72ex);}}
\newcommand{\leanpartial}{\tikz[baseline=-0.6ex]{\fill[leanfull] (0,-0.6ex) arc (-90:90:0.6ex) -- cycle;\draw[leanfull,line width=0.5pt] (0,0) circle (0.6ex);}}
\newcommand{\leannone}{\tikz[baseline=-0.6ex]\draw[leannone,line width=0.6pt] (0,0) circle (0.55ex);}
\newcommand{\leanmark}[1]{\leavevmode{\reversemarginpar\marginnote{\raggedleft#1}}}
\KOMAoptions{captions=nooneline}
\author{Jan Snellman \orcidlink{0009-0002-6676-5068}\thanks{Department of Mathematics, Linköping University, 581 83 Linköping, Sweden. \texttt{jan.snellman@liu.se}.}}
\date{November 19, 2026}
\title{Hilbert functions over the exterior algebra and \(f\)-vectors in higher rank\\\medskip
\large III: Lattice points of the Kozlov polytope}
\hypersetup{
 pdfauthor={},
 pdftitle={Hilbert functions over the exterior algebra and \(f\)-vectors in higher rank},
 pdfkeywords={},
 pdfsubject={},
 pdfcreator={},
 pdflang={English}}
\let\autocite\cite

\begin{document}

\maketitle
\begin{abstract}
Let \(E\) be the exterior algebra on \(n\) generators and \(F=\bigoplus_{i=1}^rEg_i\) a graded free
\(E\)-module with \(\deg g_i=d_i\). Part I of this series identified the convex hull of the Hilbert
functions of the quotients \(F/M\), \(M\) a proper monomial submodule, as the Kozlov polytope
\(\Koz(n;\vek d)\), a Minkowski sum of shifted Kozlov simplices, and described its vertices at
rank two, conjecturing its facets there; Part II counted the Hilbert functions. Here we study the
lattice points of the polytope.

We show that
\(\Koz(n)\) is a lecture hall simplex and that every point of \(\Koz(n;\vek d)\) has a unique
chained decomposition into points of \(\Koz(n)\). Hence the lattice points of \(\Koz(n;\vek d)\) are
exactly the sums of shifted lattice points of \(\Koz(n)\), the polytope has the integer
decomposition property, every lattice point is a Hilbert function when \(n\le3\), and lattice-point
counts and Ehrhart polynomials are products of transfer matrices of the kind used in Part II. We
determine the codegree of every Kozlov polytope and show that none is Gorenstein for \(n\ge3\).
Finally we compare discrete and continuous invariants: the number of Hilbert functions against the
volume and the lattice-point count, where we prove Kupreyeva's observed limit for the fraction of
lattice points that are \(f\)-vectors and extend it to every rank. We also compare the centroid and
second moments of the Hilbert functions with those of the polytope.
\end{abstract}


\medskip

\setcounter{tocdepth}{2}
\tableofcontents
\section{Introduction}
\label{sec-intro}
\subsection{Where Parts I and II leave things}
\label{sec:org8b31deb}

Let \(k\) be a field, \(E\) the exterior algebra on \(e_1,\dots,e_n\), and \(F\) the graded free
\(E\)-module
\begin{equation}\label{eq:F}
  F=\bigoplus_{i=1}^rEg_i,\qquad \deg g_i=d_i,\qquad 0=d_1\le d_2\le\dots\le d_r .
\end{equation}
A monomial
submodule \(M\subseteq F\) is a direct sum \(I_1g_1\oplus\dots\oplus I_rg_r\) of monomial ideals,
and the Hilbert function of \(F/M\) is a sum of shifted \(f\)-vectors of simplicial complexes. If
\(S_n\) denotes the set of \(f\)-vectors \((1,n,f_2,\dots,f_n)\) of the complexes on \([n]\) that
use every vertex, the achievable Hilbert functions of quotients by proper monomial submodules form
the sumset
\begin{equation}\label{eq:ff}
  \FF(n;\vek d)\;=\;S_n+\shift^{d_2}(S_n)+\dots+\shift^{d_r}(S_n)\;\subseteq\;\mathbb Z^{N+1},
  \qquad N=n+d_r ,
\end{equation}
whose elements we call \emph{ff}-vectors; here \(\shift^d\) prepends \(d\) zeros \eqref{eq:shift}. In Part II this set is written \(\FF_{\mathrm{pr}}(n;\vek d)\),
and \(\FF(n;\vek d)\) denotes there the larger sumset of \(f\)-vectors of all complexes on \([n]\),
the void complex included. In this paper \(\FF\) always means the proper set. Part I
\autocite{SnellmanArticle1} proved that its convex hull is the \textbf{Kozlov polytope}
\begin{equation}\label{eq:koz}
  \Koz(n;\vek d)\;=\;\sum_{i=1}^r\shift^{d_i}\bigl(\Koz(n)\bigr),
\end{equation}
a Minkowski sum of shifted copies of the Kozlov simplex \(\Koz(n)=\conv(S_n)\)
\autocite[Theorem 39]{SnellmanArticle1}, and determined its vertices and, conjecturally, its facets
at rank two \autocite[Section 9]{SnellmanArticle1}. Part II \autocite{SnellmanArticle2} counted
\(\FF(n;\vek d)\) by a product of transfer matrices, one per gap \(d_i-d_{i-1}\). Part I left the
vertices at every rank to this paper \autocite[Section 9.3]{SnellmanArticle1}; they, and the facets at
every rank, are instead determined in a forthcoming manuscript \autocite{SnellmanArticle4}, which also
proves Part I's two facet conjectures \autocite[Conjectures 52 and 55]{SnellmanArticle1}.

\begin{figure}[htbp]
\centering
\includegraphics[width=\linewidth]{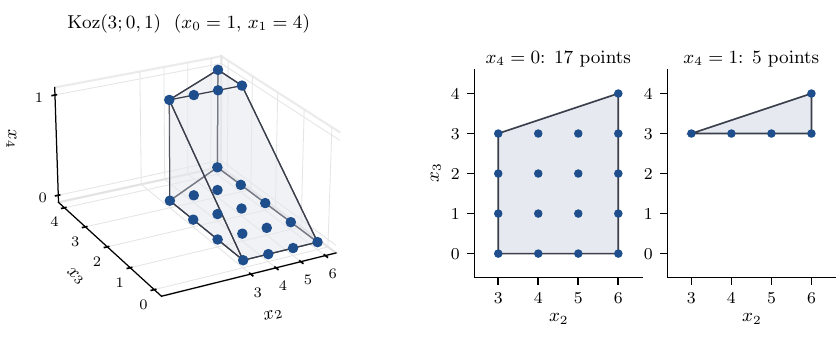}
\caption{$\Koz(3;0,1)=\Koz(3)+\shift(\Koz(3))$, a $3$-polytope with $7$ vertices in the
coordinates $(x_2,x_3,x_4)$ ($x_0=1$, $x_1=4$), with its slices $x_4=0$ and $x_4=1$. All $22$
lattice points are ff-vectors.}
\label{fig:koz3-d01}
\end{figure}

For example, for \(n=3\) and \(\vek d=(0,1)\) the ff-vectors are the \(22\) vectors
\((1,4,f_2,f_3,f_4)\) obtained by adding the \(f\)-vector of a graph on \([3]\), or of the full
triangle, to a shifted one; their convex hull \(\Koz(3;0,1)\) is a \(3\)-polytope with \(7\) vertices and
\(6\) facets (Figure \ref{fig:koz3-d01}), and every one of its lattice points is an ff-vector.

The question that remains, and the subject of this paper, is that of the lattice points.
A hull is coarse. Already for \(n=4\) the simplex \(\Koz(4)\) contains twenty lattice points of
which only sixteen are \(f\)-vectors, and the proportion falls quickly with \(n\). Which lattice
points are Hilbert functions, and how are they distributed?
\subsection{What this paper does}
\label{sec:orgb441878}

The answers, in the order of the paper. Section \ref{sec-lattice} is about the lattice points
themselves; Sections \ref{sec-volume} to \ref{sec-moments} compare discrete and continuous invariants,
one comparison to a section.

\begin{itemize}
\item \textbf{Lattice points} (Section \ref{sec-lattice}). \(\Koz(n)\) is, up to a unimodular map, the
\(\vek s\)-lecture hall simplex with \(\vek s=(\binom n2,\dots,\binom nn)\)
(Proposition \ref{prop-lecture-hall}). Every point of \(\Koz(n;\vek d)\) has a unique \emph{chained}
decomposition into points of \(\Koz(n)\), integral when the point is (Theorem \ref{thm-chained}).
Consequences: lattice points add, the polytope has the integer decomposition property, every
lattice point is an ff-vector when \(n\le3\), and an effective membership test.
\item \textbf{Volume and lattice-point count} (Section \ref{sec-volume}). The normalised volume of
\(\Koz(n;\vek d)\), from a tiling by products of simplices (Theorem \ref{thm-staircase}), and the
number of Hilbert functions against the number of lattice points.
\item \textbf{Kupreyeva's limit} (Section \ref{sec-kupreyeva}). As the number of vertices grows, the fraction of
lattice points that are \(f\)-vectors tends to a limit that Kupreyeva observed; we prove it,
explain it as a volume ratio (Theorem \ref{thm-kupreyeva-limit}), and extend it to every rank
(Theorem \ref{thm-kupreyeva-rank-r}).
\item \textbf{Ehrhart theory} (Section \ref{sec-ehrhart}). Ehrhart polynomials by transfer matrices; their
linear and quadratic coefficients at every rank, read off ranks one and two; the codegree of
every Kozlov polytope, and the fact that none with \(n\ge3\) is Gorenstein
(Theorem \ref{thm-codegree}, Corollary \ref{cor-gorenstein}).
\item \textbf{Centroids} (Section \ref{sec-centroid}). The centroid of the Hilbert functions against that of
the polytope: exact symmetries at rank one, and at equal shifts the convergence of the former to
the latter after scaling (Theorem \ref{thm-centroid-dilation}), with the whole uniform measure
(Remark \ref{rem-weak-convergence}).
\item \textbf{Second moments} (Section \ref{sec-moments}). The Hilbert functions are more spread than the
polytope for small \(n\), and concentrate for large \(n\) (Observation \ref{obs-concentration}).
\end{itemize}

Which lattice points are Hilbert functions is answered here through the canonical decomposition
(Corollary \ref{cor-membership}). There is also an intrinsic answer, a growth condition on the
sequence itself, valid for all graded submodules of \(F\) and not only for proper monomial ones:
the module version of the Kruskal--Katona theorem due to Amata and Crupi
\autocite[Theorem 4.2]{AmataCrupi2020KK}. A new proof of it, with a complete proof of the step in which
the free generators of largest degree are filled first, is the subject of a separate paper
\autocite{SnellmanArticle5}.
\subsection{What is proved, and what is only observed}
\label{sec:org3e2c560}

The distinction is kept explicit throughout, as in Parts I and II. \textbf{Proved} means proved in this
paper, or proved modulo a precise citation. \textbf{Verified} means checked by exact computation, along
a route independent of the formula being tested, over the range stated; such statements appear as
Observations or Conjectures and are not claimed to be proved. Section \ref{sec-methods} says what
was computed, how, and where the code is.

\begin{center}\footnotesize
\begin{tabular}{@{}ll@{}}
\toprule
result & status\\
\midrule
Prop.~\ref{prop-lecture-hall}, \(\Koz(n)\) as a lecture hall simplex & proved\\
Thm.~\ref{thm-hstar}, its \(h^*\)-vector & Savage--Schuster; proof included\\
Cor.~\ref{cor-zero-shift}, lattice points of dilates & proved (\(\FF\) half: Part II)\\
Cor.~\ref{cor-deficiency}, the non-ff lattice points of \(r\Koz(n)\), \(n=4,5\) & proved (finite computation)\\
Thm.~\ref{thm-chained}, the canonical chained decomposition, and its corollaries & proved; verified\\
Thm.~\ref{thm-staircase}, the staircase tiling and the volume & proved; verified\\
Thm.~\ref{thm-beta-lattice}, Thm.~\ref{thm-alexander-core}, centroid symmetries & proved\\
Thm.~\ref{thm-n3-equal-shifts}, \(n=3\), equal shifts & proved\\
Thm.~\ref{thm-kupreyeva-limit}, Kupreyeva's limit & proved (two proofs); verified\\
Thm.~\ref{thm-kupreyeva-rank-r}, its rank-\(r\) form & proved (uses Part II); verified\\
Thm.~\ref{thm-centroid-dilation}, centroid and covariance at equal shifts & proved\\
Prop.~\ref{prop-dilates}, the lattice-point centroid of \(t\Koz(n)\), \(n\le5\) & proved (finite computation)\\
Thm.~\ref{thm-codegree}, Cor.~\ref{cor-gorenstein}, codegree; never Gorenstein & proved; verified\\
Prop.~\ref{prop-quadratic-coeff}, the \(t^2\)-coefficient at rank \(r\) & proved; verified\\
Obs.~\ref{obs-hstar-unimodal}, \ref{obs-positivity}, \ref{obs-fill}, \ref{obs-beta-composite}, \ref{obs-pushforward-fails}, \ref{obs-concentration}, \ref{obs-moments-rank-r} & verified only\\
Conj.~\ref{conj-nv-monotone} & conjectured\\
\bottomrule
\end{tabular}
\end{center}
\section{Setup}
\label{sec-setup}
\subsection{Notation}
\label{sec:org1257f4a}

\begin{center}\footnotesize
\begin{tabular}{@{}lll@{}}
\toprule
symbol & meaning & where\\
\midrule
\(E\), \(F\), \(\vek d\), \(N\) & exterior algebra, free module, shift vector, \(N=n+d_r\) & \eqref{eq:F}, \eqref{eq:ff}\\
\(\delta_i\), \(\shift^d\) & gap \(d_i-d_{i-1}\); prepending \(d\) zeros & \eqref{eq:gaps}, \eqref{eq:shift}\\
\(S_n\), \(\mathcal L_n\) & \(f\)-vectors of proper complexes on \([n]\); lattice points of \(\Koz(n)\) & \S\ref{sec-setup}\\
\(\FF(n;\vek d)\) & the ff-vectors (Part II's \(\FF_{\mathrm{pr}}\)) & \eqref{eq:ff}\\
\(\Koz(n;\vek d)\), \(\RA(\vek c;\vek d)\), \(\rasimp(\vek c)\) & Kozlov polytope; right-angle polytope and simplex & \eqref{eq:koz}, Def.~\ref{def-ra}\\
\(\vek c\), \(\tilde F_j\), \(W_i\) & leg vector \((c_0,\dots,c_n)\), \(c_0=1\); vertex; window & Def.~\ref{def-ra}\\
\(\kozvec_n\), \(\vek 1\) & Part I's Kozlov vector \((\binom n1,\dots,\binom nn)\); all-ones vector & \eqref{eq:koz-as-ra}, \S\ref{sec-setup}\\
\(\mu_k\) & LYM coordinate \(x_k/c_k\), in particular \(x_k/\binom nk\) & \eqref{eq:mu}\\
\(\fl\), \(\tp\), \(\indeg\) & full up to, top degree, initial degree \(\fl+1\) & \eqref{eq:fl}--\eqref{eq:indeg}\\
\(\NV\), \(\vol\) & normalised and Euclidean volume & \eqref{eq:NV}\\
\(L_P\), \(h^*\) & Ehrhart polynomial (not to be confused with \(\mathcal L_n\)), \(h^*\)-vector & \eqref{eq:ehrhart}\\
\(L^\circ_P\) & lattice points of \(t\relint P\); reciprocity, codegree & \eqref{eq:reciprocity}\\
\(\operatorname{fill}(n;\vek d)\) & fill rate \(\lvert\FF\rvert/\#\text{lattice points}\) & \eqref{eq:fill}\\
\(a^{(k)}\) & \(k\)-th pseudo-power; Kruskal--Katona & \eqref{eq:kk}, \eqref{eq:kk-cond}\\
\(\nu_k\) & \(x_k/\dim F_k\), the rank-\(r\) analogue of \(\mu_k\) (not barycentric) & \eqref{eq:nu}\\
\(c(X)\), \(\Lambda_t\) & uniform centroid of \(X\); lattice points of \(t\Koz(n)\) & \eqref{eq:centroid}, \eqref{eq:Lambda}\\
\bottomrule
\end{tabular}
\end{center}

Vectors are set in bold, as in Part II \autocite{SnellmanArticle2} (Part I marks them with an arrow):
\(\vek d\), \(\vek p\), \(\vek c\), and \(\vek 1=(1,\dots,1)\), the all-ones vector, whose length is
always fixed by the context. One exception is deliberate: a point of the polytope under study, of
\(\Koz(n;\vek d)\), \(\RA(\vek c;\vek d)\) or a dilate, is written \(x\), not in bold, with
coordinates \(x_k\), while the parts of its decompositions are bold, \(\vek p_i\).

Throughout, \(n\) and the \textbf{shift vector} \(\vek d=(d_1,\dots,d_r)\), the degrees of the generators of
the free module \(F\) of \eqref{eq:F}, are fixed when nothing else is said. The shift vector is always sorted,
\(0=d_1\le d_2\le\dots\le d_r\); its entries are the \textbf{shifts}, \(r\) is its length, and its \textbf{gaps}
are
\begin{equation}\label{eq:gaps}
  \delta_i=d_i-d_{i-1},\qquad i=2,\dots,r .
\end{equation}

Two lattice polytopes are \textbf{unimodularly equivalent} if an affine bijection between their affine hulls
maps the one onto the other and the lattice points of the one affine hull onto those of the other.
Such a map preserves the normalised volume, the lattice points of every dilate, and so the Ehrhart
polynomial.

A few letters are reused locally with a different meaning (for instance \(u\) and \(T\) in
Appendix \ref{sec-app-chained}); each such use is introduced where it occurs and does not leave its
section.
\subsection{Hilbert functions and ff-vectors}
\label{sec:org9dbff0e}

We use the notation of Parts I and II. \(E=\bigwedge(k^n)\) is graded with
\(\dim_kE_j=\binom nj\). A simplicial complex \(\Delta\) on \([n]\) (an order ideal of the Boolean
lattice) has \(f\)-vector \((f_0,\dots,f_n)\), with \(f_k\) the number of faces of cardinality
\(k\). The \emph{void} complex, with no faces at all, has \(f_0=0\); every other complex is \emph{non-void}
and has \(f_0=1\). The exterior face ideal
\(I_\Delta=(e_\sigma:\sigma\notin\Delta)\) satisfies \(\dim_k(E/I_\Delta)_k=f_k(\Delta)\). The
complex is \emph{proper} if it uses all \(n\) vertices, that is \(f_1=n\)
\autocite[Definition 10]{SnellmanArticle1}. We write \(S_n\) for the set of \(f\)-vectors of proper
complexes, so \(\lvert S_n\rvert=1,2,5,16,70,457,4908\) for \(n=1,\dots,7\), and \(\mathcal L_n\)
for the set of lattice points of \(\Koz(n)\). For an integer \(d\ge0\) the \textbf{shift} prepends \(d\)
zeros,
\begin{equation}\label{eq:shift}
  \shift^d(x_0,\dots,x_m)=(\underbrace{0,\dots,0}_{d},x_0,\dots,x_m),
\end{equation}
and for \(F\) as in \eqref{eq:F}, \(\FF(n;\vek d)\) is the sumset \eqref{eq:ff} of the introduction
\autocite[Proposition 35]{SnellmanArticle1}. We call its elements ff-vectors, or achievable Hilbert
functions, or, where no confusion is possible, Hilbert functions; a lattice point of
\(\Koz(n;\vek d)\) that is not an ff-vector is \emph{non-achievable}.

The set \(S_n\) is described by the Kruskal--Katona theorem \autocite{Katona1968,Kruskal1963}. For
integers \(k\ge1\) and \(a\ge1\), the \emph{\(k\)-th Macaulay representation} of \(a\) is the unique
expression
\begin{equation}\label{eq:macaulay}
  a=\binom{a_k}k+\binom{a_{k-1}}{k-1}+\dots+\binom{a_s}s,\qquad a_k>a_{k-1}>\dots>a_s\ge s\ge1,
\end{equation}
and the \emph{\(k\)-th pseudo-power} of \(a\) is
\begin{equation}\label{eq:kk}
  a^{(k)}=\binom{a_k}{k+1}+\binom{a_{k-1}}{k}+\dots+\binom{a_s}{s+1},\qquad 0^{(k)}=0 .
\end{equation}
A vector \((1,n,f_2,\dots,f_n)\) of non-negative integers lies in \(S_n\) if and only if
\begin{equation}\label{eq:kk-cond}
  f_{k+1}\le f_k^{(k)},\qquad 1\le k\le n-1 .
\end{equation}
That is, \(f_k\) faces of cardinality \(k\) have at most \(f_k^{(k)}\) faces of cardinality \(k+1\)
above them, and this bound is attained.
\subsection{Right-angle simplices and the two coordinate systems}
\label{sec:org860a5b9}

\begin{definition}
Let \(\vek c=(c_0,c_1,\dots,c_n)\) be a vector of positive reals with \(c_0=1\), the \textbf{leg vector}. The
\textbf{right-angle simplex} \(\rasimp(\vek c)\subseteq\mathbb R^{n+1}\) is the convex hull of the points
\begin{equation}\label{eq:Ftilde}
  \tilde F_j=(c_0,c_1,\dots,c_j,0,\dots,0),\qquad j=1,\dots,n .
\end{equation}
For a shift vector \(\vek d\) and
\(N=n+d_r\) the \textbf{right-angle polytope} is
\begin{equation}\label{eq:ra}
  \RA(\vek c;\vek d)=\sum_{i=1}^r\shift^{d_i}\bigl(\rasimp(\vek c)\bigr)\subseteq\mathbb R^{N+1} .
\end{equation}
The \textbf{window} of summand \(i\) is
\begin{equation}\label{eq:window}
  W_i=[d_i,d_i+n] ,
\end{equation}
and every point of \(\shift^{d_i}(\rasimp(\vek c))\) is supported there.
\label{def-ra}
\end{definition}

For \(\vek c=(\binom n0,\binom n1,\dots,\binom nn)\) these are \(\Koz(n)\) and \(\Koz(n;\vek d)\), and
\(\tilde F_j\) is the \(f\)-vector of the \((j-1)\)-skeleton of the full simplex on \([n]\).

This is \autocite[Definition 36]{SnellmanArticle1}, with two changes. First, the leg vector may be
real, which costs nothing. Second, the objects live in the space of the Kozlov polytope, the
\emph{degree-\(0\)-inclusive} space \(\mathbb R^{N+1}\) with coordinates \(x_0,\dots,x_N\). Part I works in
the \emph{degree-\(0\)-dropped} space \(\mathbb R^N\), with coordinates \(x_1,\dots,x_N\): the same names
with \(x_0\) forgotten. Let \(\pi\) be the projection forgetting \(x_0\). Part I's right-angle simplex
with leg vector \((c_1,\dots,c_n)\) is \(\pi(\rasimp(\vek c))\), and Part I's \(\tilde F_j\) is the
image of this one under \(\pi\). For the Kozlov polytope the two pictures are related by
\autocite[Section 8]{SnellmanArticle1}
\begin{equation}\label{eq:koz-as-ra}
  \pi\bigl(\Koz(n;\vek d)\bigr)=\RA^{\mathrm I}\bigl(\kozvec_n;\vek d\bigr)+\sum_{i\,:\,d_i\ge1}\vek e_{d_i},
  \qquad x_0=\#\lbrace i:d_i=0\rbrace\ \text{on }\Koz(n;\vek d),
\end{equation}
where \(\RA^{\mathrm I}(\kozvec_n;\vek d)\subseteq\mathbb R^N\) is Part I's right-angle polytope, whose
leg vector is the \textbf{Kozlov vector} \(\kozvec_n=(\binom n1,\dots,\binom nn)\), and \(\vek e_k\) is the
\(k\)-th standard basis vector of \(\mathbb R^N\), in the coordinates \(x_1,\dots,x_N\). Each summand
with \(d_i\ge1\) contributes its constant coordinate \(x_0=1\) in global coordinate \(d_i\), which
\(\pi\) keeps; the summands with \(d_i=0\) contribute it to \(x_0\), which \(\pi\) forgets. The vertex
and facet descriptions of Part I and of \autocite{SnellmanArticle4} are stated in Part I's coordinates.
In this paper \(\rasimp\), \(\RA\) and \(\tilde F_j\) always mean the degree-\(0\)-inclusive objects of
Definition \ref{def-ra}.

The dimension is \autocite[Proposition 45]{SnellmanArticle1}
\begin{equation}\label{eq:dim}
  \dim\Koz(n;\vek d)=\Bigl\lvert\bigcup_{i=1}^r[d_i+2,d_i+n]\Bigr\rvert ,
\end{equation}
which is \(N-1\) when every gap is at most \(n-1\). When
every gap is at least \(n-1\) we say that \(\vek d\) is in the \textbf{product regime}: consecutive windows
then overlap only in coordinates where the later summand is constant, and \(\Koz(n;\vek d)\) is
unimodularly equivalent to \(\Koz(n)^r\) (Proposition \ref{prop-product-regime}).

The rank-one facet description, for a general leg vector, is the following. For
\(\Koz(n)\) it is \autocite[Lemma 13 and Corollary 16]{SnellmanArticle1}.

\begin{proposition}
\leanmark{\leanfull}Let \(\lambda>0\) and write
\begin{equation}\label{eq:mu}
  \mu_k=x_k/c_k\quad(0\le k\le n),\qquad\mu_{n+1}=0 .
\end{equation}
Then
\(\lambda\rasimp(\vek c)=\lbrace x\in\mathbb R^{n+1}:\lambda=\mu_0=\mu_1\ge\mu_2\ge\dots\ge\mu_n\ge0\rbrace\),
and the barycentric coordinates of \(x\) with respect to \(\lambda\tilde F_1,\dots,\lambda\tilde F_n\)
are \((\mu_k-\mu_{k+1})/\lambda\), \(k=1,\dots,n\).
\label{prop-rank-one-H}
\end{proposition}

\begin{proof}
The vertex \(\lambda\tilde F_k\) has \(\mu=\lambda(1,\dots,1,0,\dots,0)\) with \(k+1\) ones, which
satisfies the chain; the chain set is convex, so it contains \(\lambda\rasimp(\vek c)\). Conversely
\(\mu=\sum_{k=1}^n(\mu_k-\mu_{k+1})(1^{k+1},0^{n-k})\) is a combination of the vertices'
\(\mu\)-images with non-negative coefficients summing to \(\mu_1=\lambda\).
\end{proof}

In the Kozlov case the \(\mu_k=x_k/\binom nk\) of \eqref{eq:mu} are the \textbf{LYM coordinates} of Part I, and
\(\Koz(n)=\lbrace1=\mu_0=\mu_1\ge\mu_2\ge\dots\ge\mu_n\ge0\rbrace\) is an order simplex. For a point
\(\vek p\in\rasimp(\vek c)\) we write
\begin{align}
  \fl(\vek p)&=\max\lbrace k:\mu_k=1\rbrace,\label{eq:fl}\\
  \tp(\vek p)&=\max\lbrace k:\mu_k>0\rbrace,\label{eq:tp}\\
  \indeg(\vek p)&=\fl(\vek p)+1,\label{eq:indeg}
\end{align}
(\(\fl\) for ``full up to''). For an \(f\)-vector \(\vek p\in S_n\), \(\tp(\vek p)\) is its top degree,
the largest \(k\) with \(p_k>0\), and \(\indeg(\vek p)\), the least \(k\) with \(p_k<\binom nk\), is
the initial degree of the corresponding ideal; these are Part II's two statistics
\autocite[Section 5]{SnellmanArticle2}. By Proposition \ref{prop-rank-one-H}, \(\vek p\) is a convex
combination of the vertices \(\tilde F_k\) with \(\fl(\vek p)\le k\le\tp(\vek p)\), in which
\(\tilde F_{\fl(\vek p)}\) and \(\tilde F_{\tp(\vek p)}\) have positive weight; we call
\([\fl(\vek p),\tp(\vek p)]\) the \textbf{barycentric support} of \(\vek p\).

\begin{figure}[htbp]
\centering
\includegraphics[width=\linewidth]{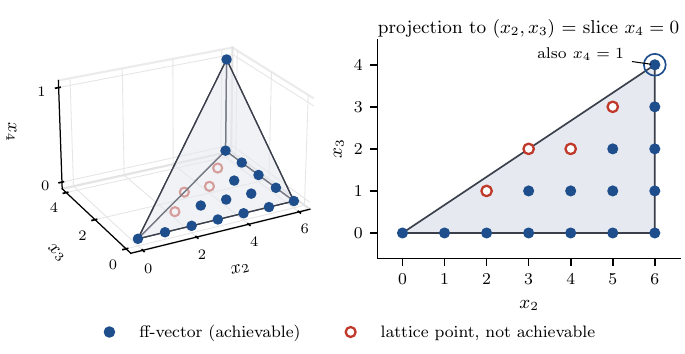}
\caption{The simplex $\Koz(4)$ (coordinates $x_2,x_3,x_4$; $x_0=1$, $x_1=4$) and its projection
to $(x_2,x_3)$. Filled dots: the $16$ $f$-vectors in $S_4$; hollow dots: the $4$ lattice points
$(2,1,0)$, $(3,2,0)$, $(4,2,0)$, $(5,3,0)$ of $\Koz(4)$ which are not $f$-vectors, all of them
violating the Kruskal--Katona bound $x_3\le x_2^{(2)}$. The only point with $x_4=1$ is the apex
$(6,4,1)$.}
\label{fig:koz4}
\end{figure}

\begin{mexample}
\(\Koz(3)\), of dimension \(2\), has vertices \((1,3,0,0)\), \((1,3,3,0)\), \((1,3,3,1)\); in \((x_2,x_3)\) it is the
triangle \(0\le3x_3\le x_2\le3\), whose five lattice points \((0,0),(1,0),(2,0),(3,0),(3,1)\) are
the \(f\)-vectors of the empty graph, one, two and three edges, and the full triangle. So
\(\mathcal L_3=S_3\).

\(\Koz(4)\), of dimension \(3\), is the tetrahedron \(1\ge x_2/6\ge x_3/4\ge x_4\ge0\) in \((x_2,x_3,x_4)\), with
\(20\) lattice points. Sixteen are \(f\)-vectors. The other four,
\[
  (x_2,x_3,x_4)=(2,1,0),\quad(3,2,0),\quad(4,2,0),\quad(5,3,0),
\]
violate the Kruskal--Katona bound \(x_3\le x_2^{(2)}\): one triangle needs three edges, two
triangles need five, three need six, and \(2^{(2)}=0\), \(3^{(2)}=1\), \(4^{(2)}=1\),
\(5^{(2)}=2\). Figure \ref{fig:koz4} shows them.
\label{ex-koz3-koz4}
\end{mexample}

\section{Lattice points: the canonical chained decomposition}
\label{sec-lattice}
From here on the leg vector \(\vek c\) of Definition \ref{def-ra} is integral: \(c_0=1\), and
\(c_1,\dots,c_n\) are positive integers. For \(\vek p\in\rasimp(\vek c)\) we use the LYM coordinates
\(\mu_k=p_k/c_k\) of \eqref{eq:mu} and the statistics \(\fl(\vek p)\) and \(\tp(\vek p)\) of Section \ref{sec-setup}.
\subsection{The Kozlov simplex is a lecture hall simplex}
\label{sec:org7d2fcde}

For positive integers \(\vek s=(s_1,\dots,s_D)\) the \textbf{\(\vek s\)-lecture hall simplex} is
\autocite{HibiOlsenTsuchiya2018,SavageSchuster2012}
\begin{equation}\label{eq:lhs}
  P^{(\vek s)}=\lbrace\lambda\in\mathbb R^D:0\le\lambda_1/s_1\le\dots\le\lambda_D/s_D\le1\rbrace .
\end{equation}

\begin{proposition}
\leanmark{\leanfull}Let \(D=n-1\) and \(\vek s=(c_n,c_{n-1},\dots,c_2)\), so that \(s_i=c_{n+1-i}\). The affine
hull of \(\rasimp(\vek c)\) is \(A=\lbrace p\in\mathbb R^{n+1}:p_0=1,\ p_1=c_1\rbrace\), and the maps
\begin{equation}\label{eq:lh-maps}
\begin{aligned}
  \varphi&\colon A\to\mathbb R^D, & (p_0,p_1,p_2,\dots,p_n)&\mapsto(p_n,p_{n-1},\dots,p_2),\\
  \varphi^{-1}&\colon\mathbb R^D\to A, & (\lambda_1,\dots,\lambda_D)&\mapsto(1,c_1,\lambda_D,\lambda_{D-1},\dots,\lambda_1),
\end{aligned}
\end{equation}
that is \(\lambda_i=p_{n+1-i}\), are mutually inverse affine bijections with
\(\varphi(A\cap\mathbb Z^{n+1})=\mathbb Z^D\) and \(\varphi(\rasimp(\vek c))=P^{(\vek s)}\). So
\(\rasimp(\vek c)\) and \(P^{(\vek s)}\) are unimodularly equivalent. Moreover, with
\(\bar{\vek s}=(s_D,\dots,s_1)\) the reversed sequence, the map
\begin{equation}\label{eq:lh-reversal}
  \operatorname{rev}\colon\mathbb R^D\to\mathbb R^D,\qquad \operatorname{rev}(\lambda)_i=s_{D+1-i}-\lambda_{D+1-i}\quad(1\le i\le D),
\end{equation}
is a unimodular affine bijection with \(\operatorname{rev}(P^{(\vek s)})=P^{(\bar{\vek s})}\); this is the
reversal equivalence of \autocite{HibiOlsenTsuchiya2018}. In particular \(\Koz(n)\) is unimodularly
equivalent to the \(\vek s\)-lecture hall simplex with
\(\vek s=(\binom nn,\dots,\binom n2)=(\binom n0,\dots,\binom n{n-2})\), and to the one with the
reversed sequence \((\binom n2,\dots,\binom nn)\). Conversely, for every \(\vek s\in\mathbb Z_{>0}^D\),
\(P^{(\vek s)}=\varphi(\rasimp(\vek c))\) for \(n=D+1\) and \(\vek c=(1,1,s_D,\dots,s_1)\).
\label{prop-lecture-hall}
\end{proposition}

\begin{proof}
The two maps in \eqref{eq:lh-maps} are inverse to each other by inspection, and both are affine.
The vertices \(\tilde F_j\) \eqref{eq:Ftilde} all have \(p_0=1\) and \(p_1=c_1\), so
\(\rasimp(\vek c)\subseteq A\); both have dimension \(D\), the vectors
\(\tilde F_j-\tilde F_1=(0,0,c_2,\dots,c_j,0,\dots,0)\) (\(2\le j\le n\)) being linearly independent, so
\(A\) is the affine hull. Since \(c_1\) is an integer, \(A\cap\mathbb Z^{n+1}\) is the set of integer
points with \(p_0=1\) and \(p_1=c_1\), which \(\varphi\) maps onto \(\mathbb Z^D\), and \(\varphi^{-1}\)
maps \(\mathbb Z^D\) back into it.

On \(A\), \(\mu_0=\mu_1=1\), so by Proposition \ref{prop-rank-one-H} (at dilation factor \(1\))
\(\rasimp(\vek c)=\lbrace p\in A:1\ge p_2/c_2\ge\dots\ge p_n/c_n\ge0\rbrace\). In the coordinates
\(\lambda_i=p_{n+1-i}\), with \(s_i=c_{n+1-i}\), the chain read from the right is
\(0\le\lambda_1/s_1\le\dots\le\lambda_D/s_D\le1\), which is \eqref{eq:lhs}. So
\(\varphi(\rasimp(\vek c))=P^{(\vek s)}\).

The map \(\operatorname{rev}\) is a signed permutation of the coordinates followed by an integral translation,
so it maps \(\mathbb Z^D\) onto \(\mathbb Z^D\); it is an involution once \(\vek s\) and
\(\bar{\vek s}\) are exchanged. With \(\bar s_i=s_{D+1-i}\),
\(\operatorname{rev}(\lambda)_i/\bar s_i=1-\lambda_{D+1-i}/s_{D+1-i}\), so it carries the chain
\(0\le\lambda_1/s_1\le\dots\le\lambda_D/s_D\le1\) to
\(0\le\operatorname{rev}(\lambda)_1/\bar s_1\le\dots\le\operatorname{rev}(\lambda)_D/\bar s_D\le1\), and \(P^{(\vek s)}\) onto
\(P^{(\bar{\vek s})}\).

For the Kozlov polytope, \(c_k=\binom nk\) and \(\binom nk=\binom n{n-k}\). For the converse,
\(\vek c=(1,1,s_D,\dots,s_1)\) has \((c_n,\dots,c_2)=(s_1,\dots,s_D)=\vek s\).
\end{proof}

So the Ehrhart theory of \(\Koz(n)\) is that of a lecture hall simplex, and much of
Section \ref{sec-ehrhart} at rank one is a citation. The same holds for every integral leg vector. In
particular Kupreyeva's Ehrhart computations for the simplices of \(f\)-vectors with faces of at most
three or four elements \autocite{Kupreyeva2019}, which are the simplices \(\rasimp(\vek c)\) with
\(c_k=\binom pk\) of Section \ref{sec-kupreyeva}, are computations for lecture hall simplices.
\subsection{Chained tuples}
\label{sec:org070e7e1}

\begin{definition}
A tuple \((\vek p_1,\dots,\vek p_r)\in\bigl(\rasimp(\vek c)\bigr)^r\) is \textbf{chained} for
\(\vek d\) if
\begin{equation}\label{eq:chained}
  \tp(\vek p_{i-1})+d_{i-1}\;\le\;\fl(\vek p_i)+d_i\qquad(i=2,\dots,r),
\end{equation}
equivalently \(\tp(\vek p_{i-1})<\indeg(\vek p_i)+\delta_i\), with the gap \(\delta_i\) of
\eqref{eq:gaps}. For
\(f\)-vectors this is the chaining condition of \autocite[Definition 19]{SnellmanArticle2}.
\label{def-chained-points}
\end{definition}

For a summand \(i\) we distinguish its \textbf{own degree} \(m\), the index of a coordinate of
\(\vek p_i\), from the \textbf{global degree} \(k=m+d_i\), the coordinate of \(\RA(\vek c;\vek d)\) in which
it lands; the leg \(c_m\) is the \textbf{capacity} of own degree \(m\).
For a tuple \((\vek p_1,\dots,\vek p_r)\) and a global degree \(k\), the summands \emph{present} in
degree \(k\) are those \(i\) with \(k-d_i\in[0,n]\), and summand \(i\) contributes its own coordinate
\((\vek p_i)_{k-d_i}\) there. We call summand \(i\) \textbf{full} in degree \(k\) if that coordinate is the
whole capacity \(c_{k-d_i}\) (so \(\mu_{k-d_i}=1\)), \textbf{zero} if it is \(0\), and \textbf{partial} if it lies
strictly between.

In words, the definition says: in global degree \(\tp(\vek p_{i-1})+d_{i-1}\), the last degree in
which summand \(i-1\) is non-zero, summand \(i\) must still be full. So at every global degree \(k\)
the summands present there are, in the order of \(i\): zero, then at most one partial, then full.

\begin{theorem}
\leanmark{\leanfull}(\emph{Canonical chained decomposition.}) Let \(\vek c\) be an integral leg vector and \(\vek d\) a
shift vector, and let \(x\in\RA(\vek c;\vek d)\).
\begin{enumerate}
\item (\emph{Existence.}) If \((c_2,\dots,c_n)\) is log-concave, that is \(c_k^2\ge c_{k-1}c_{k+1}\) for
\(3\le k\le n-1\), then \(x=\sum_i\shift^{d_i}(\vek p_i)\) for a chained tuple of points
\(\vek p_i\in\rasimp(\vek c)\).
\item (\emph{Uniqueness.}) There is at most one such chained tuple, and it is given by the
\textbf{water-filling map}: with \(c_m=0\) for \(m\notin[0,n]\),

\begin{equation}\label{eq:water}
\begin{aligned}
  Q_i(k)&=\max\Bigl(0,\;x_k-\sum_{l>i}c_{k-d_l}\Bigr),\qquad Q_0\equiv0,\\
  (\vek p_i)_m&=Q_i(m+d_i)-Q_{i-1}(m+d_i).
\end{aligned}
\end{equation}
\item (\emph{Integrality.}) If \(x\) is integral, so is every \(\vek p_i\) of that tuple.
\end{enumerate}
The leg vector of the Kozlov polytope satisfies the hypothesis of (1), binomial coefficients being
log-concave.
\label{thm-chained}
\end{theorem}

In words: pour the amount \(x_k\) into the summands present in global degree \(k\), filling them
from the last one down, each to its capacity. Then \(Q_i(k)\) is what is left for summands
\(1,\dots,i\) once the summands \(l>i\) are full, and summand \(i\) receives the part of it that
does not go on to summands \(1,\dots,i-1\).

The proof is in Appendix \ref{sec-app-chained}, followed by a comparison with the decomposition
theorem of Brändén and Solus for lecture hall polytopes, to which Theorem \ref{thm-chained} reduces at
\(\vek d=\vek 0\) (Remark \ref{rem-chained-credit}).
\subsection{Consequences}
\label{sec:org8a44cf5}

\begin{corollary}
\leanmark{\leanfull}(\emph{Lattice points add.}) If \((c_2,\dots,c_n)\) is log-concave, then
\begin{equation}\label{eq:lattice-add}
  \RA(\vek c;\vek d)\cap\mathbb Z^{N+1}=\sum_{i=1}^r\shift^{d_i}\bigl(\rasimp(\vek c)\cap\mathbb Z^{n+1}\bigr).
\end{equation}
In particular
\begin{equation}\label{eq:ff-in-lattice}
  \FF(n;\vek d)=\sum_{i=1}^r\shift^{d_i}(S_n)\;\subseteq\;\sum_{i=1}^r\shift^{d_i}(\mathcal L_n)=\Koz(n;\vek d)\cap\mathbb Z^{N+1},
\end{equation}
and the lattice points of \(\Koz(n;\vek d)\) are in bijection with the chained tuples of
elements of \(\mathcal L_n\).
\label{cor-additivity}
\end{corollary}

\begin{proof}
In \eqref{eq:lattice-add}, the right-hand side lies in the left-hand side because a sum of lattice
points is a lattice point and \(\RA(\vek c;\vek d)\) is the Minkowski sum of the
\(\shift^{d_i}\rasimp(\vek c)\). That every lattice point of \(\RA(\vek c;\vek d)\) is such a sum is
parts (1) and (3) of Theorem \ref{thm-chained}. In \eqref{eq:ff-in-lattice}, the first equality is the
definition \eqref{eq:ff}, the inclusion holds because \(S_n\subseteq\mathcal L_n\), and the last
equality is \eqref{eq:lattice-add} for the leg vector of the Kozlov polytope. The bijection is
parts (1)--(3) of Theorem \ref{thm-chained}.
\end{proof}

\begin{corollary}
\leanmark{\leanfull}For \(n\le3\), every lattice point of \(\Koz(n;\vek d)\) is an ff-vector, for every \(r\) and every
\(\vek d\).
\label{cor-small-n}
\end{corollary}

\begin{proof}
\(\mathcal L_2=S_2=\lbrace(1,2,0),(1,2,1)\rbrace\) and \(\mathcal L_3=S_3\) (Example
\ref{ex-koz3-koz4}), so the two sumsets of Corollary \ref{cor-additivity} coincide. For \(n=1\) the
polytope is a point.
\end{proof}

At \(n=4\) this fails already at rank one, by the four points of Example \ref{ex-koz3-koz4}.

\begin{figure}[htbp]
\centering
\includegraphics[width=\linewidth]{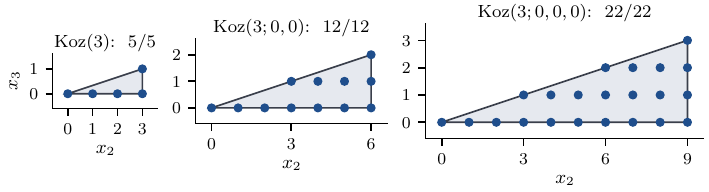}
\caption{The Kozlov polytopes $\Koz(3)$, $\Koz(3;0,0)=2\Koz(3)$ and
$\Koz(3;0,0,0)=3\Koz(3)$, drawn in the free coordinates $(x_2,x_3)$ (here $x_0=r$, $x_1=3r$).
Dots are the ff-vectors, i.e.\ the Hilbert functions of quotients of $E^r$ by proper monomial
submodules. In all three cases every lattice point is achievable: $5$, $12$ and $22$ points
respectively, as Corollary~\ref{cor-small-n} requires.}
\label{fig:koz3-dilates}
\end{figure}

\begin{corollary}
\leanmark{\leancond}(\emph{Which lattice points are ff-vectors.}) A lattice point \(x\) of \(\Koz(n;\vek d)\) lies in
\(\FF(n;\vek d)\) if and only if every part of its canonical decomposition lies in \(S_n\), that is,
satisfies the Kruskal--Katona inequalities.
\label{cor-membership}
\end{corollary}

\begin{proof}
``If'' is trivial. Conversely, by Part II's chaining reformulation for proper submodules
\autocite[Theorems 23 and 41, with Lemma 40]{SnellmanArticle2}, every element of \(\FF(n;\vek d)\) is the
sum of a chained tuple of elements of \(S_n\); that tuple is a chained tuple of points of
\(\Koz(n)\), so by uniqueness it is the canonical one.
\end{proof}

So membership is effective: compute \(r\) water-filling differences and check \(r\) rank-one
Kruskal--Katona conditions. The new information, compared with Part II, is that the canonical
tuple exists for \emph{every} lattice point, achievable or not, so that the non-achievable lattice
points are exactly those whose canonical tuple contains a lattice point of \(\Koz(n)\) that is not
an \(f\)-vector. A second characterisation, directly in terms of the sequence, is Amata and Crupi's
module version of the Kruskal--Katona theorem \autocite[Theorem 4.2]{AmataCrupi2020KK}, which holds
for all graded submodules of \(F\); a new proof of it is given in \autocite{SnellmanArticle5}.

\begin{mexample}
\(\Koz(3;0,1)\) has dimension \(3\), \(7\) vertices \((1,4,x_2,x_3,x_4)\) with \((x_2,x_3,x_4)\) one of
\[\small
\begin{array}{lllllll}
(3,0,0), & (3,3,0), & (3,3,1), & (6,0,0), & (6,4,0), & (6,4,1), & (6,3,1),
\end{array}
\]
and \(22\) lattice points, all ff-vectors. In the same coordinates
\[
  \Koz(3;0,1)=\lbrace3\le x_2\le6,\ 3x_3\le x_2+6,\ 3x_4\le x_3,\ 0\le x_4\le1\rbrace ,
\]
with \(6\) facets: the vertices of this system are exactly the \(7\) points above, and none of its
inequalities is redundant. Some canonical decompositions:
\[
\begin{array}{l}
(1,4,3,0,0)=(1,3,0,0)+\shift(1,3,0,0),\\ (1,4,4,3,1)=(1,3,1,0)+\shift(1,3,3,1),\\
(1,4,6,4,0)=(1,3,3,1)+\shift(1,3,3,0),\\ (1,4,6,4,1)=(1,3,3,1)+\shift(1,3,3,1).
\end{array}
\]
The first part's edge count is \(x_2-3\). It can carry the triangle only if it is all of
\((1,3,3,1)\), which forces \(x_2=6\); the chaining condition \(\tp(\vek p_1)\le\fl(\vek p_2)+1\) then
requires \(\fl(\vek p_2)\ge2\), that is, the second part contains all three edges, so \(x_3\ge4\).
These are the points \((1,4,6,4,0)\) and \((1,4,6,4,1)\) above.
\label{ex-koz3-d01}
\end{mexample}

\begin{figure}[htbp]
\centering
\includegraphics[width=0.582\linewidth]{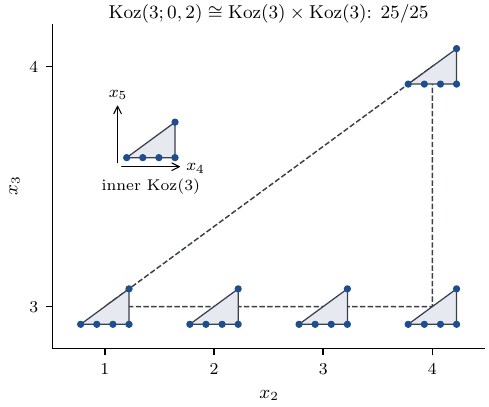}
\caption{In the product regime the Kozlov polytope splits: $\Koz(3;0,2)\cong\Koz(3)\times\Koz(3)$,
with $(x_2,x_3)\in(1,3)+\Koz(3)$ (outer dashed triangle) and $(x_4,x_5)\in\Koz(3)$ (small
triangles). Each of the $5$ outer lattice points carries a full copy of the $5$ inner ones, and
all $25=5\cdot 5$ lattice points are ff-vectors.}
\label{fig:koz3-d02}
\end{figure}

\begin{mexample}
\(\Koz(3;0,2)\) is in the product regime (Section \ref{sec-setup}): its windows overlap only where
the later summand is constant (\(x_0=1\), \(x_1=n\)), so the polytope is unimodularly equivalent to
\(\Koz(3)\times\Koz(3)\) (Proposition \ref{prop-product-regime}), of dimension \(4\), with \(9\)
vertices, \(6\) facets, and \(25\) lattice points, all ff-vectors (Figure \ref{fig:koz3-d02}).
\label{ex-koz3-d02}
\end{mexample}

\begin{figure}[htbp]
\centering
\includegraphics[width=\linewidth]{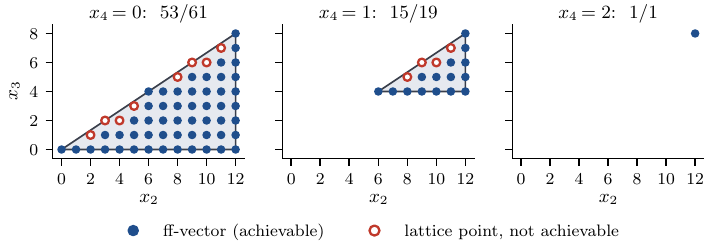}
\caption{The rank-two Kozlov polytope $\Koz(4;0,0)=2\Koz(4)$ ($x_0=2$, $x_1=8$), sliced at
$x_4=0,1,2$. Of its $81$ lattice points, $69$ are ff-vectors (filled) and $12$ are not (hollow).
Each non-achievable point is the topmost lattice point of its column ($x_2$, $x_4$ fixed), just
below the facet $\mu_2=\mu_3$.}
\label{fig:koz4-d00}
\end{figure}

\begin{mexample}
\(\Koz(4;0,1)\) has dimension \(4\), \(11\) vertices \((1,5,x_2,x_3,x_4,x_5)\) with
\((x_2,x_3,x_4,x_5)\) one of
\[\small
\begin{array}{llllll}
(4,0,0,0), & (4,6,0,0), & (4,6,4,0), & (4,6,4,1), & (10,0,0,0), & (10,6,4,0),\\
(10,6,4,1), & (10,10,0,0), & (10,10,4,1), & (10,10,5,0), & (10,10,5,1), &
\end{array}
\]
\(214\) lattice points and \(162\) ff-vectors. In the same coordinates
\[
\begin{array}{l}
\Koz(4;0,1)=\lbrace4\le x_2\le10,\ 3x_3\le2x_2+10,\ 4x_4\le x_3+10,\\
\phantom{\Koz(4;0,1)=\lbrace}3x_4\le2x_3,\ 4x_5\le x_4,\ 0\le x_5\le1\rbrace ,
\end{array}
\]
with \(8\) facets: the vertices of this system are exactly the \(11\) points above, and none of
its inequalities is redundant. The \(52\) non-achievable lattice points split, by which
canonical part fails Kruskal--Katona, into \(28\) whose second part fails and \(24\) whose first
part fails; in none do both parts fail. Two examples, in coordinates \((x_0,\dots,x_5)\):
\begin{gather*}
  (1,5,4,2,1,0)=(1,4,0,0,0)+\shift(1,4,2,1,0),\\
  (1,5,6,7,0,0)=(1,4,2,1,0)+\shift(1,4,6,0,0).
\end{gather*}
The second can be checked by hand: in any representation as \(\vek s_1+\shift(\vek s_2)\) with
\(\vek s_i\in S_4\), degree \(2\) forces \(f_2(\vek s_1)=6-4=2\); two edges carry no triangle, so
\(f_3(\vek s_1)=0\), and degree \(3\) forces \(f_2(\vek s_2)=7>\binom42\), which is impossible.

Two more rank-two polytopes, which reappear in Table \ref{tab:fill}. \(\Koz(4;0,0)=2\Koz(4)\) has \(81\) lattice points, of which
\(69\) are ff-vectors; the \(12\) others are drawn in Figure \ref{fig:koz4-d00}. For
\(\Koz(4;0,2)\), \(246\) of the \(386\) lattice points are ff-vectors, and among the \(140\) others
there are, for the first time in these examples, \(16\) whose two canonical parts both fail.
\label{ex-koz4-d01}
\end{mexample}

\begin{corollary}
\leanmark{\leanfull}(\emph{Integer decomposition property.}) If \((c_2,\dots,c_n)\) is log-concave, then
\(\RA(\vek c;\vek d)\) has the integer decomposition property (IDP) for every \(\vek d\): every lattice
point of \(t\RA(\vek c;\vek d)\) is a sum of \(t\) lattice points of \(\RA(\vek c;\vek d)\). In
particular every Kozlov polytope is IDP.
\label{cor-idp}
\end{corollary}

\begin{proof}
The dilate \(t\RA(\vek c;\vek d)\) is \(\RA(\vek c;\vek d^{\times t})\), where \(\vek d^{\times t}\) repeats
each \(d_i\) \(t\) times:
\[
  \vek d^{\times t}=(\underbrace{d_1,\dots,d_1}_{t},\underbrace{d_2,\dots,d_2}_{t},\dots,\underbrace{d_r,\dots,d_r}_{t}).
\]
By Corollary \ref{cor-additivity} a lattice point of it is a sum of \(rt\) shifted
lattice points of \(\rasimp(\vek c)\); regroup them as \(t\) sums each with one point from every
\(i\), and each such sum is a lattice point of \(\RA(\vek c;\vek d)\).
\end{proof}

At \(\vek d=\vek 0\) this is the integer decomposition property of lecture hall simplices
\autocite{BrandenSolus2020}, which Hibi, Olsen and Tsuchiya had proved for monotone \(\vek s\) and
conjectured in general \autocite{HibiOlsenTsuchiya2018}; the Kozlov sequence
\((\binom n2,\dots,\binom nn)\) is monotone only for \(n\le5\).
\section{Discrete volume versus volume}
\label{sec-volume}
Three sizes of the same polytope compete: \(\lvert\FF(n;\vek d)\rvert\), the number of achievable
Hilbert functions, counted in Part II; the number of lattice points; and the volume. By
Corollary \ref{cor-additivity} the first two are cardinalities of sumsets of the nested sets
\(S_n\subseteq\mathcal L_n\).
\subsection{Normalisation}
\label{sec:org4a7fd9d}

Let \(U=\bigcup_i[d_i+2,d_i+n]\), so that \(D=\lvert U\rvert=\dim\Koz(n;\vek d)\) by \eqref{eq:dim}. Every coordinate outside
\(U\) is constant on \(\Koz(n;\vek d)\), so its affine hull is parallel to the coordinate subspace
\(\mathbb R^U\), and the coordinate projection to \(\mathbb R^U\) is an isometry on it that maps the
lattice points of the affine hull onto \(\mathbb Z^U\), a unimodular equivalence onto its image. We
therefore use two volumes: the \textbf{Euclidean} \(D\)-dimensional volume \(\vol\), which is the one in
Kupreyeva's thesis \autocite{Kupreyeva2019}, and the \textbf{normalised} volume
\begin{equation}\label{eq:NV}
  \NV=D!\,\vol ,
\end{equation}
for which a unimodular simplex has \(\NV=1\) and which equals the sum of the \(h^*\)-vector
(Section \ref{sec-ehrhart}).

\begin{proposition}
\leanmark{\leanpartial}\(\NV(\rasimp(\vek c))=c_2c_3\cdots c_n\). In particular \(\NV(\Koz(n))=\prod_{k=2}^n\binom nk\) and
\(\vol(\Koz(n))=\frac1{n!}\prod_{k=1}^n\binom nk\).
\label{prop-simplex-volume}
\end{proposition}

\begin{proof}
In LYM coordinates \(\rasimp(\vek c)\) is the image of the order simplex
\(\lbrace1\ge\mu_2\ge\dots\ge\mu_n\ge0\rbrace\) under \(\operatorname{diag}(c_2,\dots,c_n)\), and the order
simplex is unimodular: its edge vectors \(e_2,e_2+e_3,\dots\) from \(0\) form a basis of
\(\mathbb Z^{n-1}\). Divide by \((n-1)!\) and note \(\binom n1=n\).
\end{proof}

The Euclidean form is Kupreyeva's \autocite[Section 6.1]{Kupreyeva2019}, proved there by the same
triangular determinant.
\subsection{The staircase tiling}
\label{sec:org940a807}

Assume every gap \(\delta_i\) \eqref{eq:gaps} is at most \(n-1\), so that \(U=[2,N]\) and \(D=N-1\). A
\textbf{staircase index} is \(\vek t=(t_0,\dots,t_r)\) with
\begin{equation}\label{eq:staircase-index}
  1=t_0\le t_1\le\dots\le t_r=N,\qquad t_{i-1}-d_i\ge1,\qquad t_i-d_i\le n ,
\end{equation}
where \(N=n+d_r\) as before. Its \textbf{cell} is
\begin{equation}\label{eq:cell}
  C_{\vek t}=\sum_{i=1}^r\shift^{d_i}\conv\lbrace\tilde F_{\alpha_i},\dots,\tilde F_{\beta_i}\rbrace,
  \qquad \alpha_i=t_{i-1}-d_i,\quad \beta_i=t_i-d_i .
\end{equation}
The face
\(\conv\lbrace\tilde F_\alpha,\dots,\tilde F_\beta\rbrace\) varies exactly in the own coordinates
\((\alpha,\beta]\), so summand \(i\) occupies the global block \((t_{i-1},t_i]\); these blocks partition
\([2,N]\). An \textbf{orthoscheme} is a simplex \(\conv\lbrace v_0,\dots,v_m\rbrace\) whose path edges
\(v_1-v_0,\ v_2-v_1,\ \dots,\ v_m-v_{m-1}\) are pairwise orthogonal; we call these edges its
\textbf{legs}. The
right-angle simplex \(\rasimp(\vek c)\) is one, since \(\tilde F_{k+1}-\tilde F_k=c_{k+1}\vek e_{k+1}\),
and so is each of its faces \(\conv\lbrace\tilde F_\alpha,\dots,\tilde F_\beta\rbrace\). So
\(C_{\vek t}\) is a lattice translate of a product of orthoschemes on complementary coordinate blocks.

\begin{theorem}
\leanmark{\leannone}If \((c_2,\dots,c_n)\) is log-concave and every gap is at most \(n-1\), the cells \(C_{\vek t}\)
tile \(\RA(\vek c;\vek d)\): they cover it and have pairwise disjoint interiors. Hence
\[
  \NV\bigl(\RA(\vek c;\vek d)\bigr)=\sum_{\vek t}\binom{D}{t_1-t_0,\ t_2-t_1,\ \dots,\ t_r-t_{r-1}}
  \prod_{i=1}^r\ \prod_{k=t_{i-1}+1}^{t_i}c_{k-d_i}.
\]
\label{thm-staircase}
\end{theorem}

\begin{proof}
\emph{Cover.} Take the canonical decomposition of \(x\) (Theorem \ref{thm-chained}) and put
\(\ell_i=\fl(\vek p_i)+d_i\), \(\tau_i=\tp(\vek p_i)+d_i\); chaining says
\(\tau_i\le\ell_{i+1}\). Choose \(t_i\in[\tau_i,\ell_{i+1}]\cap[d_{i+1}+1,d_i+n]\),
non-empty because \(\tau_i\le d_i+n\), \(\ell_{i+1}\ge d_{i+1}+1\), \(\tau_i\le\ell_{i+1}\)
and \(d_{i+1}+1\le d_i+n\). Then \(t_{i-1}\le\ell_i\le\tau_i\le t_i\), and the barycentric
support \([\fl(\vek p_i),\tp(\vek p_i)]\) of \(\vek p_i\) lies in \([\alpha_i,\beta_i]\), so
\(x\in C_{\vek t}\).

\emph{Disjointness.} Let \(\vek t\ne\vek t'\), let \(i\) be the first index with \(t_i\ne t'_i\), say
\(t_i<t'_i\), and put \(k=t_i+1\) and \(T=\sum_{l>i}c_{k-d_l}\). In \(C_{\vek t}\) the summands
\(l\le i\) vanish in degree \(k\) and the others are at most full, so \(x_k\le T\). In \(C_{\vek t'}\)
the summands \(l<i\) vanish, the summands \(l>i\) are full (their \(\alpha_l\ge k-d_l\)), and summand
\(i\) is non-negative, so \(x_k\ge T\). The hyperplane \(x_k=T\) separates the cells.

\emph{Volume.} An orthoscheme with integer legs has \(\NV\) equal to the product of its legs
(Proposition \ref{prop-simplex-volume}), and for a product of lattice polytopes on complementary
coordinate blocks \(D!\,\vol\) is the multinomial coefficient times the product of the factors'
\(\NV\).
\end{proof}

If some gap is at least \(n-1\), \(\RA(\vek c;\vek d)\) is a lattice product of the sub-sums over
the maximal runs with gaps at most \(n-2\), and \(\NV\) is the multinomial in the run dimensions
times the product of the runs' \(\NV\). (A gap of exactly \(n-1\) is covered by both descriptions,
which agree: the later summand is constant on the coordinates the two windows share. When every gap
is at least \(n-1\) this is the product regime, Proposition \ref{prop-product-regime} below.) For non-log-concave legs the formula fails, consistently
with the failure of covering: for \(\vek c=(1,2,1,3,1,3)\), \(\vek d=(0,1)\) it gives \(540\) against
the true \(1500\).

\begin{corollary}
\leanmark{\leannone}For \(0\le d\le n-1\), \(m=n-d\) and \(D=n-1+d\),
\[
  \NV\bigl(\Koz(n;0,d)\bigr)=\NV(\Koz(n))\cdot\sum_{j=1}^{m}\binom{n-1+d}{d+j-1}\prod_{k=j+1}^{j+d}\binom nk .
\]
\label{cor-rank2-volume}
\end{corollary}

\begin{proof}
The staircase indices are \((1,d+j,n+d)\), \(1\le j\le m\), and the cells are
\begin{equation}\label{eq:cells-rank2}
  C_j=C_{(1,d+j,n+d)}=\conv\lbrace\tilde F_1,\dots,\tilde F_{d+j}\rbrace+\shift^d\conv\lbrace\tilde F_j,\dots,\tilde F_n\rbrace,
\end{equation}
and
\[
  \NV(C_j)=\binom D{d+j-1}\prod_{k=2}^{d+j}\binom nk\prod_{k=j+1}^n\binom nk .
\]
\end{proof}

At \(d=0\) the sum is \(2^{n-1}\), so \(\NV(2\Koz(n))=2^{n-1}\NV(\Koz(n))\), as it must be; at \(d=n-1\)
there is one cell and \(\NV=\binom{2n-2}{n-1}\NV(\Koz(n))^2\), the product.

\begin{figure}[htbp]
\centering
\includegraphics[width=\linewidth]{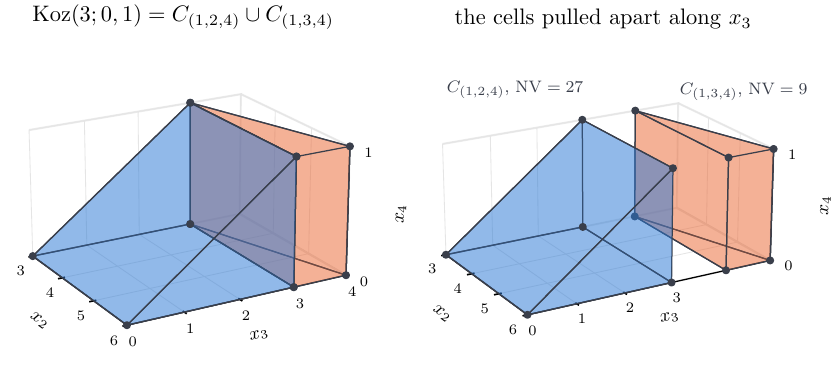}
\caption{The staircase tiling of $\Koz(3;0,1)$ (Example~\ref{ex-volume}), in the coordinates
$(x_2,x_3,x_4)$, with $x_0=1$, $x_1=4$. The cell
$C_{(1,2,4)}=\conv\lbrace\tilde F_1,\tilde F_2\rbrace+\shift\Koz(3)$ is the prism
$[3,6]\times\lbrace0\le3x_4\le x_3\le3\rbrace$, with $\NV=27$; the cell
$C_{(1,3,4)}=\Koz(3)+\shift\conv\lbrace\tilde F_2,\tilde F_3\rbrace$ is the prism
$\lbrace0\le3(x_3-3)\le x_2-3\le3\rbrace\times[0,1]$, with $\NV=9$. They meet in the plane
$x_3=3$, the separating hyperplane $x_k=T$ of the proof of Theorem~\ref{thm-staircase}, with $k=3$.}
\label{fig:koz3-d01-tiling}
\end{figure}

\begin{mexample}
\(\Koz(3;0,1)\) has dimension \(D=3\) and two cells (Figure \ref{fig:koz3-d01-tiling}), \(j=1,2\) in the notation of
Corollary \ref{cor-rank2-volume}: \(j=1\), staircase index \((1,2,4)\), blocks
\(\lbrace2\rbrace\mid\lbrace3,4\rbrace\), \(\NV=\binom31\cdot3\cdot(3\cdot1)=27\); and \(j=2\), staircase
index \((1,3,4)\), blocks \(\lbrace2,3\rbrace\mid\lbrace4\rbrace\), \(\NV=\binom32\cdot(3\cdot1)\cdot1=9\).
So \(\NV=36\) and \(\vol=36/3!=6\). \(\Koz(4;0,1)\) has dimension \(D=4\) and three cells: \(j=1\) gives
\(\binom41\cdot6\cdot(6\cdot4\cdot1)=576\), \(j=2\) gives \(\binom42\cdot(6\cdot4)\cdot(4\cdot1)=576\), \(j=3\)
gives \(\binom43\cdot(6\cdot4\cdot1)\cdot1=96\); \(\NV=1248\), \(\vol=1248/4!=52\).
\label{ex-volume}
\end{mexample}

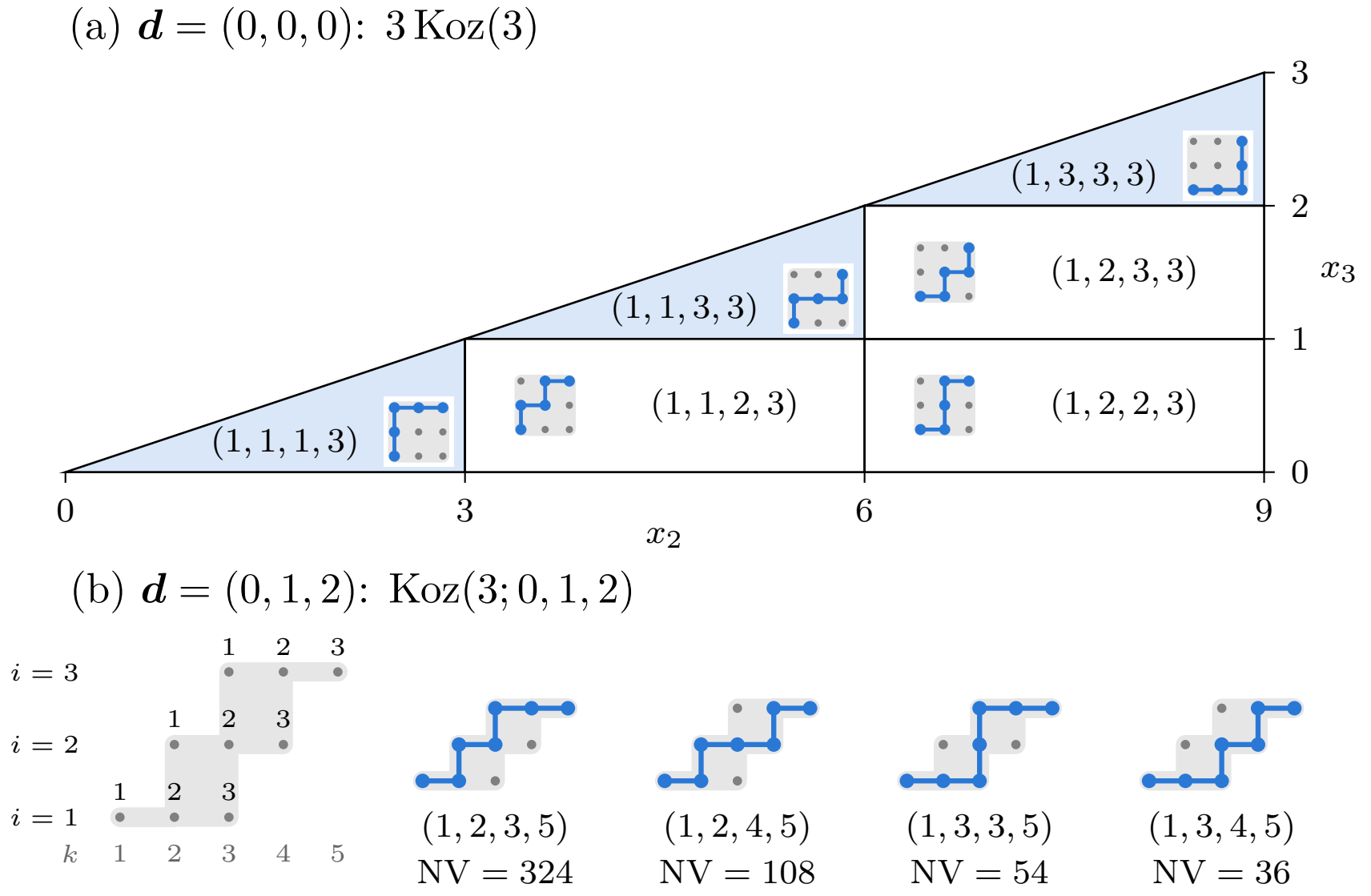
\begin{figure}[htbp]
\centering
\definecolor{stairblue}{HTML}{2A78D6}
\newcommand{\stairglyph}[5]{%
  \filldraw[band] (1,0) -- (#1,0) -- (#1,1) -- ({#2+#1},1) -- ({#2+#1},2) -- ({#3+#1},2)
    -- ({#3+1},2) -- ({#3+1},1) -- ({#2+1},1) -- ({#2+1},0) -- cycle;
  \foreach \x in {1,...,#1} {\node[dot] at (\x,0) {}; \node[dot] at ({\x+#2},1) {};
    \node[dot] at ({\x+#3},2) {};}
  \draw[spath] (1,0) -- (#4,0) -- (#4,1) -- (#5,1) -- (#5,2) -- ({#3+#1},2);
  \pgfmathtruncatemacro{\stairN}{#3+#1}%
  \foreach \x in {1,...,#4} {\node[vdot] at (\x,0) {};}
  \foreach \x in {#4,...,#5} {\node[vdot] at (\x,1) {};}
  \foreach \x in {#5,...,\stairN} {\node[vdot] at (\x,2) {};}}
\begin{tikzpicture}[
    font=\scriptsize,
    dot/.style={circle, fill=black!50, inner sep=0pt, minimum size=2.2pt},
    vdot/.style={circle, fill=stairblue, inner sep=0pt, minimum size=3.4pt},
    spath/.style={stairblue, line width=1.2pt, line cap=round, line join=round},
    band/.style={fill=black!10, draw=black!10, line width=4.5pt, line join=round},
    small glyph/.style={x=0.2cm, y=0.2cm,
      dot/.style={circle, fill=black!50, inner sep=0pt, minimum size=1.7pt},
      vdot/.style={circle, fill=stairblue, inner sep=0pt, minimum size=2.6pt},
      spath/.style={stairblue, line width=0.9pt, line cap=round, line join=round},
      band/.style={fill=black!10, draw=black!10, line width=3pt, line join=round}},
    cell/.style={draw=black, line width=0.5pt},
  ]
  \begin{scope}[x=1.1cm, y=1.1cm]
    \node[anchor=south west, font=\small] at (-0.1,3.05) {(a) \(\vek d=(0,0,0)\): \(3\Koz(3)\)};
    \filldraw[cell, fill=stairblue!18] (0,0) -- (3,0) -- (3,1) -- cycle;
    \filldraw[cell, fill=stairblue!18] (3,1) -- (6,1) -- (6,2) -- cycle;
    \filldraw[cell, fill=stairblue!18] (6,2) -- (9,2) -- (9,3) -- cycle;
    \filldraw[cell, fill=white] (3,0) rectangle (6,1);
    \filldraw[cell, fill=white] (6,0) rectangle (9,1);
    \filldraw[cell, fill=white] (6,1) rectangle (9,2);
    \foreach \x in {0,3,6,9} \draw (\x,0) -- (\x,-0.08) node[below] {\(\x\)};
    \foreach \y in {0,1,2,3} \draw (9,\y) -- (9.08,\y) node[right] {\(\y\)};
    \node[below] at (4.5,-0.3) {\(x_2\)};
    \node[right] at (9.3,1.5) {\(x_3\)};
    \foreach \a/\b/\ta/\tb in {0/0/1/1, 3/1/1/3, 6/2/3/3} {
      \node[anchor=east] at (\a+2.33,\b+0.22) {\((1,\ta,\tb,3)\)};
      \begin{scope}[shift={(\a+2.47,\b+0.12)}, small glyph, shift={(-1,0)}]
        \fill[white] (0.55,-0.45) rectangle (3.45,2.45);
        \stairglyph{3}{0}{0}{\ta}{\tb}
      \end{scope}}
    \foreach \a/\b/\ta/\tb in {3/0/1/2, 6/0/2/2, 6/1/2/3} {
      \node at (\a+1.95,\b+0.5) {\((1,\ta,\tb,3)\)};
      \begin{scope}[shift={(\a+0.42,\b+0.32)}, small glyph, shift={(-1,0)}]
        \stairglyph{3}{0}{0}{\ta}{\tb}
      \end{scope}}
  \end{scope}
  \node[anchor=south west, font=\small] at (-0.1cm,-1.3cm) {(b) \(\vek d=(0,1,2)\): \(\Koz(3;0,1,2)\)};
  \begin{scope}[shift={(0cm,-2.85cm)}, x=0.45cm, y=0.6cm]
    \filldraw[band] (1,0) -- (3,0) -- (3,1) -- (4,1) -- (4,2) -- (5,2) -- (3,2) -- (3,1)
      -- (2,1) -- (2,0) -- cycle;
    \foreach \i/\d in {1/0, 2/1, 3/2}
      \foreach \j in {1,2,3} {
        \node[dot] at ({\j+\d},\i-1) {};
        \node[font=\tiny, above=0.2pt] at ({\j+\d},\i-1) {\(\j\)};}
    \foreach \x in {1,...,5} \node[font=\tiny, text=black!60] at (\x,-0.5) {\(\x\)};
    \foreach \i in {1,2,3} \node[font=\tiny, anchor=east] at (0.55,\i-1) {\(i=\i\)};
    \node[font=\tiny, text=black!60, anchor=east] at (0.55,-0.5) {\(k\)};
  \end{scope}
  \foreach \ta/\tb/\nv [count=\k] in {2/3/324, 2/4/108, 3/3/54, 3/4/36} {
    \begin{scope}[shift={({0.65cm+\k*2.0cm},-2.55cm)}, x=0.3cm, y=0.3cm]
      \stairglyph{3}{1}{2}{\ta}{\tb}
      \node at (3,-1.3) {\((1,\ta,\tb,5)\)};
      \node at (3,-2.5) {\(\NV=\nv\)};
    \end{scope}}
\end{tikzpicture}
\caption{Staircase indices as lattice paths (Example~\ref{ex-staircase-paths} and
Remark~\ref{rem-staircase-cayley}). In each small
diagram, row \(i\), counted from below, has dots in the columns \(d_i+1,\dots,d_i+n\), the dot in
column \(k\) standing for the vertex \(\tilde F_{k-d_i}\) of the \(i\)-th summand; the shaded band
lies between the two extreme paths. The path of \(\vek t\) runs along row \(i\) from column
\(t_{i-1}\) to column \(t_i\), and the dots it visits there span the \(i\)-th summand of
\(C_{\vek t}\). (a) \(\Koz(3;0,0,0)=3\Koz(3)\) in the coordinates \((x_2,x_3)\), with \(x_0=3\),
\(x_1=9\): its six cells, each with its index and its path. The shaded triangles are translates of
\(\Koz(3)\), with \(\NV=3\); the rectangles are sums of the edges
\(\conv\lbrace\tilde F_1,\tilde F_2\rbrace\) and \(\conv\lbrace\tilde F_2,\tilde F_3\rbrace\) of two
different summands, with \(\NV=6\). (b) \(\Koz(3;0,1,2)\): the band, with the own index of each
vertex above its dot and the column below, and the four paths with their \(\NV(C_{\vek t})\), which
add up to \(522=\NV(\Koz(3;0,1,2))\).}
\label{fig:staircase-paths}
\end{figure}

\begin{mexample}
(\emph{Staircase indices as lattice paths.}) Put the vertex \(\shift^{d_i}\tilde F_m\) of summand \(i\)
at the grid point \((m+d_i,i)\): column = global degree, row = summand. The vertices fill the \textbf{band}
\begin{equation}\label{eq:band}
  \mathcal B=\lbrace(k,i)\in\mathbb Z^2:1\le i\le r,\ d_i+1\le k\le d_i+n\rbrace .
\end{equation}
To a staircase index \(\vek t\) attach the east--north lattice path with corners
\[
  (t_0,1)\to(t_1,1)\to(t_1,2)\to(t_2,2)\to\dots\to(t_{r-1},r)\to(t_r,r).
\]
Its north step from row \(i\) to row \(i+1\) lies in \(\mathcal B\) if and only if
\(d_{i+1}+1\le t_i\le d_i+n\), which is \eqref{eq:staircase-index}. So staircase indices are the
east--north paths in \(\mathcal B\) from \((1,1)\) to \((N,r)\), and
\(\vek t^{\uparrow}\le\vek t\le\vek t^{\downarrow}\) for
\(\vek t^{\uparrow}=(1,d_2+1,\dots,d_r+1,N)\) and \(\vek t^{\downarrow}=(1,n,d_2+n,\dots,d_{r-1}+n,N)\),
which step north as early and as late as possible. In row \(i\) the path visits \((k,i)\) for
\(t_{i-1}\le k\le t_i\), that is \(\tilde F_{\alpha_i},\dots,\tilde F_{\beta_i}\) of
\eqref{eq:cell}. Give the east step \((k-1,i)\to(k,i)\) the weight \(c_{k-d_i}\), the leg
\(\tilde F_{k-d_i}-\tilde F_{k-d_i-1}=c_{k-d_i}\vek e_{k-d_i}\). Then the term of \(\vek t\) in
Theorem \ref{thm-staircase} is
\[
  \NV(C_{\vek t})=\binom{D}{t_1-t_0,\ \dots,\ t_r-t_{r-1}}\prod_{(k-1,i)\to(k,i)}c_{k-d_i} ,
\]
the product over the \(D=N-1\) east steps. For \(\Koz(3;0,1)\) the rows are \([1,3]\) and
\([2,4]\), \(t_1\in\lbrace2,3\rbrace\), and the two paths are the two cells of
Example \ref{ex-volume}. \(\Koz(3;0,1,2)\) has four (Figure \ref{fig:staircase-paths}(b)), and
\(\NV=324+108+54+36=522\).
\label{ex-staircase-paths}
\end{mexample}

\begin{remark}
(\emph{Equal shifts: the staircase triangulation.}) At \(\vek d=\vek 0\) every row of
Example \ref{ex-staircase-paths} is \([1,n]\), the band is the whole \(r\times n\) grid, and the paths
are all its \(\binom{n+r-2}{r-1}\) monotone paths from corner to corner. Let
\(\Delta^{n-1}=\conv\lbrace\vek e_1,\dots,\vek e_n\rbrace\) be the standard simplex, and draw the
vertex \((\vek e_i,\vek e_j)\) of \(\Delta^{r-1}\times\Delta^{n-1}\) as the dot in row \(i\) and
column \(j\). Then these paths are the maximal simplices of the \emph{staircase triangulation} of
\(\Delta^{r-1}\times\Delta^{n-1}\) \autocite[Example 2.2]{CeballosPadrolSarmiento2015}. The Cayley
trick sends a simplex \(s\) of a triangulation of \(\Delta^{r-1}\times\Delta^{n-1}\) to the
Minkowski cell \(\sum_i\conv\lbrace\vek e_j:(\vek e_i,\vek e_j)\in s\rbrace\) of \(r\Delta^{n-1}\),
and is a bijection from triangulations to fine mixed subdivisions \cite[Theorem 1.4 and
Section 1.4]{Santos2005}, \cite[Section 2]{CeballosPadrolSarmiento2015}. The linear map \(\vek e_j\mapsto\tilde F_j\)
takes \(\Delta^{n-1}\) onto \(\rasimp(\vek c)\), and the cell of the path of \(\vek t\) onto
\(C_{\vek t}\). So at \(\vek d=\vek 0\) the staircase tiling is the image of the fine mixed
subdivision of \(r\Delta^{n-1}\) that corresponds to the staircase triangulation
(Figure \ref{fig:staircase-paths}(a)). There it needs no log-concavity:
\(\operatorname{diag}(c_0,\dots,c_n)\) maps \(\rasimp(\vek 1)\) onto \(\rasimp(\vek c)\), and each
cell for \(\vek 1\), which is log-concave, onto the cell for \(\vek c\) with the same index. For
\(\vek d\ne\vek 0\) the diagonal map does not commute with the shifts, and the non-log-concave
leg vector after Theorem \ref{thm-staircase} shows that log-concavity cannot then be dropped.
\label{rem-staircase-cayley}
\end{remark}

\begin{proposition}
\leanmark{\leanfull}(\emph{The product regime.}) If every gap is at least \(n-1\), \(\Koz(n;\vek d)\) is unimodularly
equivalent to \(\Koz(n)^r\). Consequently \(\vol(\Koz(n;\vek d))=\vol(\Koz(n))^r\), the lattice-point
count is \(\lvert\mathcal L_n\rvert^r\), the Ehrhart polynomial (Section \ref{sec-ehrhart}) is
\(L_{\Koz(n)}(t)^r\), and
\(\lvert\FF(n;\vek d)\rvert=\lvert S_n\rvert^r\).
\label{prop-product-regime}
\end{proposition}

\begin{proof}
Consecutive windows overlap only in coordinates where the later summand is constant (\(x_0=1\),
\(x_1=n\)), so the coordinate projection onto the summands' blocks is an affine bijection onto
\(\Koz(n)^r\) that maps lattice points onto lattice points, a unimodular equivalence. Volumes and lattice points of a product, also of its dilates,
multiply. For the last statement, the sumset map \(S_n^r\to\FF(n;\vek d)\) is onto by definition,
and the block projection is a left inverse of it on \(\Koz(n)^r\), so it is a bijection; this is
also \autocite[Proposition 42]{SnellmanArticle2}.
\end{proof}

\begin{center}\small
\begin{tabular}{@{}rrrrrr@{}}
\toprule
\(n\) & \(d=0\) & \(d=1\) & \(d=2\) & \(d=3\) & \(d\ge n-1\)\\
\midrule
2 & 2 & 1 & & & 1\\
3 & 6 & 6 & 9/4 & & 9/4\\
4 & 32 & 52 & 56 & 16 & 16\\
5 & 1000/3 & 5125/6 & 64375/36 & 34375/18 & 15625/36\\
6 & 7200 & 29475 & 107550 & 256500 & 50625\\
\bottomrule
\end{tabular}

\smallskip
\begin{tabular}{@{}rrrrrr@{}}
\toprule
\(n\) & \(d=0\) & \(d=1\) & \(d=2\) & \(d=3\) & \(d\ge n-1\)\\
\midrule
3 & 12 & 36 & 54 & & 54\\
4 & 192 & 1248 & 6720 & 11520 & 11520\\
5 & 8000 & 102500 & 1287500 & 9625000 & 17500000\\
\bottomrule
\end{tabular}
\captionof{table}{Euclidean volume (top) and normalised volume (bottom) of \(\Koz(n;0,d)\), from
Corollary~\ref{cor-rank2-volume}; every entry was also computed independently with Normaliz. For
\(n=6\), \(d=4\) the Euclidean volume is \(270000\), larger than both neighbours.}
\label{tab:volumes}
\end{center}

The Euclidean volume is not monotone in \(d\), but nothing should be read into this: the dimension
\(D=n-1+d\) changes with \(d\), so the entries measure different kinds of volume. The normalised
volume, although it too depends on the dimension, is monotone in every computed case.

\begin{conjecture}
\(\NV(\Koz(n;0,d))\) is weakly increasing in \(d\). (Verified for \(n\le7\) and
\(0\le d\le n-1\); for \(d\ge n-1\) the value is constant, by Proposition \ref{prop-product-regime}.)
\label{conj-nv-monotone}
\end{conjecture}
\subsection{Hilbert functions, lattice points and volume}
\label{sec:orgb22c76e}

Volume compares the polytope with itself across \(n\) and \(\vek d\); to compare it with the set of
Hilbert functions, a dimension-free size is needed. Lattice-point counts are dimension-free, and
they are the right comparison for \(\lvert\FF\rvert\). The \textbf{fill rate} is
\begin{equation}\label{eq:fill}
  \operatorname{fill}(n;\vek d)=\frac{\lvert\FF(n;\vek d)\rvert}{\#\bigl(\Koz(n;\vek d)\cap\mathbb Z^{N+1}\bigr)},
  \qquad \operatorname{fill}(n)=\frac{\lvert S_n\rvert}{\lvert\mathcal L_n\rvert} .
\end{equation}

\begin{center}\small
\begin{tabular}{@{}rrrrr@{\qquad}rrrrr@{}}
\toprule
\(n\) & \(\lvert S_n\rvert\) & \(\lvert\mathcal L_n\rvert\) & \(\vol\) & fill &
\(n\) & \(d\) & \(\lvert\FF\rvert\) & lattice & fill\\
\midrule
2 & 2 & 2 & 1 & 1 & 4 & 0 & 69 & 81 & 0.852\\
3 & 5 & 5 & 3/2 & 1 & 4 & 1 & 162 & 214 & 0.757\\
4 & 16 & 20 & 4 & 0.8 & 4 & 2 & 246 & 386 & 0.637\\
5 & 70 & 162 & 125/6 & 0.432 & 4 & 3 & 256 & 400 & 0.640\\
6 & 457 & 1981 & 225 & 0.231 & 5 & 0 & 627 & 1108 & 0.566\\
7 & 4908 & 60993 & 84035/16 & 0.0805 & 5 & 1 & 1858 & 4484 & 0.414\\
8 & 95248 & 3576524 & 2458624/9 & 0.0266 & 5 & 2 & 3805 & 14214 & 0.268\\
 & & & &  & 5 & 3 & 4857 & 26119 & 0.186\\
 & & & &  & 5 & 4 & 4900 & 26244 & 0.187\\
\bottomrule
\end{tabular}
\captionof{table}{Left: rank one. Right: \(\Koz(n;0,d)\). Counts of \(\FF\) from Part II, lattice
points by the transfer matrix of Corollary~\ref{cor-transfer} (Section~\ref{sec-ehrhart}) and
independently by Normaliz (for
\(\lvert\mathcal L_8\rvert\), by a separate count of integer LYM chains).}
\label{tab:fill}
\end{center}

\begin{figure}[htbp]
\centering
\includegraphics[width=\linewidth]{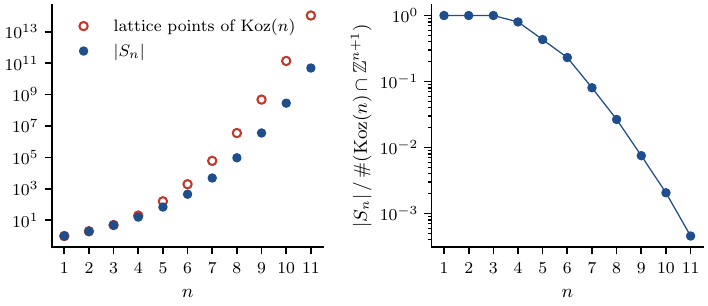}
\caption{The number of $f$-vectors of proper complexes on $[n]$ compared with the number of
lattice points of their convex hull $\Koz(n)$ (left, logarithmic scale), and the fraction of
lattice points that are $f$-vectors (right), for $n\le11$. The two agree for $n\le3$
(Corollary~\ref{cor-small-n}); from $n=4$ on the fraction decreases, to $4.6\cdot10^{-4}$ at
$n=11$.}
\label{fig:fraction}
\end{figure}

Three ways of letting the rank \(r\) grow, with \(n\) fixed, are within reach, and the fill rate
behaves differently in each. Cases 2 and 3 use results proved in Section \ref{sec-ehrhart} below.

\begin{enumerate}
\item \textbf{Product regime.} By Proposition \ref{prop-product-regime}, \(\operatorname{fill}(n;\vek d)=\operatorname{fill}(n)^r\) whenever
every gap is at least \(n-1\), with no limit needed; with gaps held in that regime the fill rate
tends to \(0\) geometrically in \(r\) for \(n\ge4\).
\item \textbf{Zero shifts.} Write \(0^r\) for the shift vector of \(r\) zeros. Then
\(\operatorname{fill}(n;0^r)\to1\) as \(r\to\infty\) (Corollary \ref{cor-zero-shift}): in the dilation
limit the ff-vectors fill the polytope to first order. The centroid and the covariance of the
ff-vectors converge too, after scaling by \(1/r\) and \(1/r^2\), to those of \(\Koz(n)\)
(Theorem \ref{thm-centroid-dilation}), and so does the whole uniform measure
(Remark \ref{rem-weak-convergence}).
\item \textbf{Constant positive gap.} Both counts grow like the Perron roots of transfer matrices, Part II's
for \(\FF\) and that of Corollary \ref{cor-transfer} for the lattice points.
\end{enumerate}

\begin{observation}
For a constant gap \(\delta\ge1\) the ratio of the two Perron roots is below \(1\) for \(n\ge4\), so
the fill rate tends to \(0\) geometrically: for \(n=4,5,6\) and \(\delta=1\) the ratios are
\(0.939\), \(0.890\), \(0.891\); for \(\delta=2\), \(0.796\), \(0.622\), \(0.609\). At rank two the
fill rate is not monotone in \(d\): it is smallest at \(d=n-2\), just before the product regime
(\(0.6373<0.64\) at \(n=4\), \(0.18596<0.18671\) at \(n=5\), \(0.053171<0.053219\) at \(n=6\)). Both
statements are verified by computation only, for \(n\le7\).
\label{obs-fill}
\end{observation}
\section{Kupreyeva's limit}
\label{sec-kupreyeva}
Kupreyeva \autocite[Section 6.2]{Kupreyeva2019} studied a limit of a different kind: the number of
vertices grows while the size of the faces stays bounded. For fixed \(n\) she counted two things:
\begin{itemize}
\item the \(f\)-vectors with \(f_1=p\) and no face of more than \(n\) elements, her \(F^n_p\);
\item the lattice points of their convex hull, projected to \((f_2,\dots,f_n)\), her \(T^n_p\).
\end{itemize}

She found closed forms for \(3\le n\le6\): a polynomial, and a quasi-polynomial in \(p\) of period
\(3\), \(12\), \(60\) and \(60\) for \(n=3,4,5,6\). She checked that the two have the same degree for
\(n\le30\), observed that the ratio of their leading coefficients is a number \(q_n\), and gave an
area argument at \(n=3\). Theorem \ref{thm-kupreyeva-limit} proves the observation and explains
\(q_n\) as a volume ratio, and Theorem \ref{thm-kupreyeva-rank-r} extends it to every rank.

In this section, as in \autocite[Section 9]{SnellmanArticle2}, \(p\ge n\) is the number of vertices and
\(n\) bounds the size of a face. Put

\begin{equation}\label{eq:kup-defs}
\begin{aligned}
  \vek c^{(p)}&=\Bigl(\binom p0,\binom p1,\dots,\binom pn\Bigr),\\
  \Koz^{(p)}(n)&=\rasimp\bigl(\vek c^{(p)}\bigr),\\
  \Koz^{(p)}(n;\vek d)&=\RA\bigl(\vek c^{(p)};\vek d\bigr).
\end{aligned}
\end{equation}
At rank \(r\) this truncates every summand at its own degree \(n\), which is the faithful analogue.
The canonical decomposition commutes with this truncation, in the sense that padding the parts by
zeros gives the canonical decomposition of the padded point (Lemma \ref{lem-kup-truncation}). The other
natural choice, truncation at global degree \(n\), gives the summands different leg vectors, and is
not covered by Theorem \ref{thm-chained}.

Let \(S_n^{(p)}\) be the set of \(f\)-vectors \((1,p,f_2,\dots,f_n)\) of the complexes on \([p]\)
that use every vertex and have no face of more than \(n\) elements, and put
\begin{equation}\label{eq:kup-sets}
  \mathcal L_n^{(p)}=\Koz^{(p)}(n)\cap\mathbb Z^{n+1},\qquad
  \FF^{(p)}(n;\vek d)=\sum_{i=1}^r\shift^{d_i}\bigl(S_n^{(p)}\bigr).
\end{equation}
The latter is Part II's \(\FF^{(p)}_{\mathrm{pr}}(n;\vek d)\) \autocite[Theorem 60]{SnellmanArticle2}. At \(p=n\) these are
\(\Koz(n)\), \(\Koz(n;\vek d)\), \(S_n\), \(\mathcal L_n\) and \(\FF(n;\vek d)\). In module terms,
\(\FF^{(p)}(n;\vek d)\) is the set of Hilbert functions of the quotients of
\(\bigoplus_i(E/E_{>n})g_i\), with \(E\) on \(p\) generators, by proper monomial submodules.

The leg vector \(\vek c^{(p)}\) is log-concave, so Theorem \ref{thm-chained} and its corollaries hold
for \(\Koz^{(p)}(n;\vek d)\). The \(f\)-vectors of complexes on \([p]\) satisfy \(\mu_{k+1}\le\mu_k\)
in the LYM coordinates \eqref{eq:mu}, here \(\mu_k=f_k/\binom pk\): a face of size \(k+1\) has \(k+1\) faces of size
\(k\), and a face of size \(k\) lies in at most \(p-k\) faces of size \(k+1\). The vertices of
\(\Koz^{(p)}(n)\) are the \(f\)-vectors of the skeleta of the full simplex. Hence
\(\Koz^{(p)}(n)=\conv S_n^{(p)}\): Kupreyeva's \(F^n_p\) is \(S_n^{(p)}\), and her simplex
\(T^n_p\) is \(\Koz^{(p)}(n)\), whose lattice points are \(\mathcal L_n^{(p)}\). We write
\begin{equation}\label{eq:fill-p}
  \operatorname{fill}^{(p)}(n;\vek d)=\frac{\lvert\FF^{(p)}(n;\vek d)\rvert}
  {\#\bigl(\Koz^{(p)}(n;\vek d)\cap\mathbb Z^{N+1}\bigr)} .
\end{equation}

In LYM coordinates \(\Koz^{(p)}(n)\) is, for every \(p\), the order simplex \(\Delta_n\) of
\eqref{eq:Rn} below (Proposition \ref{prop-rank-one-H}). The Kruskal--Katona bound \eqref{eq:kk-cond},
\(f_{k+1}\le f_k^{(k)}\), with the pseudo-power \eqref{eq:kk}, becomes
\(\mu_{k+1}\le\mu_k^{(k+1)/k}\) to leading order in \(p\) (Lemma \ref{lem-pseudo-power}). This cuts the
region \(R_n\) out of \(\Delta_n\):
\begin{equation}\label{eq:Rn}
\begin{aligned}
  \Delta_n&=\bigl\lbrace\mu\in\mathbb R^{n-1}:1\ge\mu_2\ge\dots\ge\mu_n\ge0\bigr\rbrace,\\
  R_n&=\bigl\lbrace\mu\in\mathbb R^{n-1}:0\le\mu_2\le1,\ 0\le\mu_{k+1}\le\mu_k^{(k+1)/k}
  \ (2\le k\le n-1)\bigr\rbrace .
\end{aligned}
\end{equation}
The region \(R_n\) is not convex: its upper boundary is the graph of a convex function.

\begin{theorem}
\leanmark{\leannone}(\emph{Kupreyeva's limit.}) Let \(n\ge2\). As \(p\to\infty\),
\[
  \bigl\lvert\mathcal L_n^{(p)}\bigr\rvert\sim\vol(\Delta_n)\prod_{k=2}^n\binom pk,\qquad
  \bigl\lvert S_n^{(p)}\bigr\rvert\sim\vol(R_n)\prod_{k=2}^n\binom pk,
\]
with \(\vol\Delta_n=1/(n-1)!\) and
\[
  \vol R_n=\prod_{k=2}^{n-1}\frac{2k}{(n-k+1)(n+k)}=\frac{2^{n-2}\,(n+1)!}{(2n-1)!}.
\]
Consequently
\begin{equation}\label{eq:qn}
  \lim_{p\to\infty}\frac{\lvert S_n^{(p)}\rvert}{\lvert\mathcal L_n^{(p)}\rvert}
  =\frac{\vol R_n}{\vol\Delta_n}=q_n=\frac{(n+1)!}{2\,(2n-1)!!},
\end{equation}
so \(q_2=1\), \(q_3=4/5\), \(q_4=4/7\), \(q_5=8/21\) and \(q_6=8/33\).
\label{thm-kupreyeva-limit}
\end{theorem}

The proof is in Appendix \ref{sec-app-kupreyeva}, and it was found by the language model (see the
disclosure section). It has three steps.
\begin{enumerate}
\item In LYM coordinates the lattice \(\mathbb Z^{n-1}\) becomes \(\prod_k\binom pk^{-1}\mathbb Z\),
whose mesh tends to \(0\) in every coordinate. So the lattice points of a fixed region with null
boundary, counted with weight \(\prod_k\binom pk^{-1}\), tend to its volume
(Lemma \ref{lem-anisotropic-count}).
\item The pseudo-power satisfies \(f^{(k)}/\binom p{k+1}=(f/\binom pk)^{(k+1)/k}+o(1)\), uniformly in
\(0\le f\le\binom pk\) (Lemma \ref{lem-pseudo-power}). This squeezes the \(f\)-vectors between the
lattice points of two regions whose volumes tend to \(\vol R_n\).
\item \(\vol R_n\) is an iterated integral (Lemma \ref{lem-rho}).
\end{enumerate}

The proof uses neither Theorem \ref{thm-chained} nor Part II. At \(n=3\) the value \(4/5\) is already
in Kupreyeva's closed forms, and she also obtained it by an interpolation-and-area argument
\autocite[Section 4.1]{Kupreyeva2019}; here \(R_3\) has area \(\int_0^1t^{3/2}\,dt=2/5\) against \(1/2\) for
the triangle \(\Delta_3\). A second proof, from Linusson's enumeration of \(f\)-vectors
\autocite{Linusson1999}, follows the first in the appendix, and Remark \ref{rem-kupreyeva-linusson} there
combines the two into the rate: \(\lvert S_n^{(p)}\rvert/\lvert\mathcal L_n^{(p)}\rvert=q_n+O(1/p)\).

The rank-\(r\) analogue holds for every shift vector. Group the shifts by value, into \(g\)
distinct values \(D_j\) with multiplicities \(m_j\), and put
\begin{equation}\label{eq:aj}
\begin{gathered}
  \lbrace d_1,\dots,d_r\rbrace=\lbrace D_1<\dots<D_g\rbrace,\qquad m_j=\#\lbrace i:d_i=D_j\rbrace,\\
  a_1=2,\qquad a_j=\max(2,\,n-D_j+D_{j-1}+1)\quad(j\ge2),\qquad \ell_j=n+1-a_j .
\end{gathered}
\end{equation}
Here \(a_j-1\) is the own degree up to which, for almost every lattice point, the summands with shift
\(D_j\) are full: below \(a_j\) they are constant (\(j=1\)) or held full by the summands of the
previous shift (\(j\ge2\)); \(\ell_j\) is the number of the remaining own degrees
\(a_j,\dots,n\), in which they are free. For \(2\le a\le b\) let
\begin{equation}\label{eq:rho}
  \rho(a,b)=\prod_{k=a}^{b-1}\frac{2k}{(b-k+1)(b+k)}
\end{equation}
be the volume of \(\lbrace0\le\mu_a\le1,\ 0\le\mu_{k+1}\le\mu_k^{(k+1)/k}\ (a\le k<b)\rbrace\), and
put \(\rho(b+1,b)=1\). For \(m\ge1\) let
\begin{equation}\label{eq:Phi}
  \Phi_n(a;m)=\sum_{a-1=s_0\le s_1\le\dots\le s_m=n}\ \prod_{l=1}^m\rho(s_{l-1}+1,s_l) .
\end{equation}

\begin{theorem}
\leanmark{\leannone}(\emph{Kupreyeva's limit at rank \(r\).}) For every \(n\ge2\) and every shift vector \(\vek d\),
\[
  \lim_{p\to\infty}\operatorname{fill}^{(p)}(n;\vek d)
  =\prod_{j=1}^g\frac{\ell_j!\,\Phi_n(a_j;m_j)}{m_j^{\ell_j}} .
\]
In particular:
\begin{enumerate}
\item If every gap \(\delta_i\) is positive, the limit is \(q_n\prod_{i=2}^rq(n,\delta_i)\), where
\begin{equation}\label{eq:qndelta}
  q(n,\delta)=\frac{2^{\delta-1}\,(n-1)!\,(2n-\delta)!}{(n-\delta)!\,(2n-1)!}
  \quad(1\le\delta\le n-1),\qquad q(n,\delta)=q_n\quad(\delta\ge n-1).
\end{equation}
So \(q(n,1)=1\), and in the product regime the limit is \(q_n^r\).
\item If \(\vek d=(0,\dots,0)\), the limit is
\(\vol(R_n+\dots+R_n)/\vol(r\Delta_n)=(n-1)!\,\Phi_n(2;r)/r^{n-1}\), the Minkowski sum having
\(r\) terms. For \(n=3\) it is \(1-1/(5r)\).
\end{enumerate}
\label{thm-kupreyeva-rank-r}
\end{theorem}

In words: at a positive gap \(\delta\), almost every lattice point has its later summand full in
its own degrees up to \(n-\delta\). Only the top \(\min(\delta,n-1)\) degrees are free, and the
Kruskal--Katona conditions on them form a rank-one problem on fewer coordinates, with limit
\(q(n,\delta)\). At equal shifts the summands share their coordinates, and the canonical
decomposition stacks them, so the limit is a Minkowski sum of copies of the non-convex region
\(R_n\). As \(r\to\infty\) that ratio tends to \(1\), the \(p\)-direction counterpart of
Corollary \ref{cor-zero-shift}. The proof is in Appendix \ref{sec-app-kupreyeva}. It uses
Corollary \ref{cor-membership}, with \(p\) in place of \(n\), and so rests on Part II's chaining
reformulation; the volume formula at equal shifts also uses Theorem \ref{thm-staircase}.

\begin{figure}[htbp]
\centering
\includegraphics[width=\linewidth]{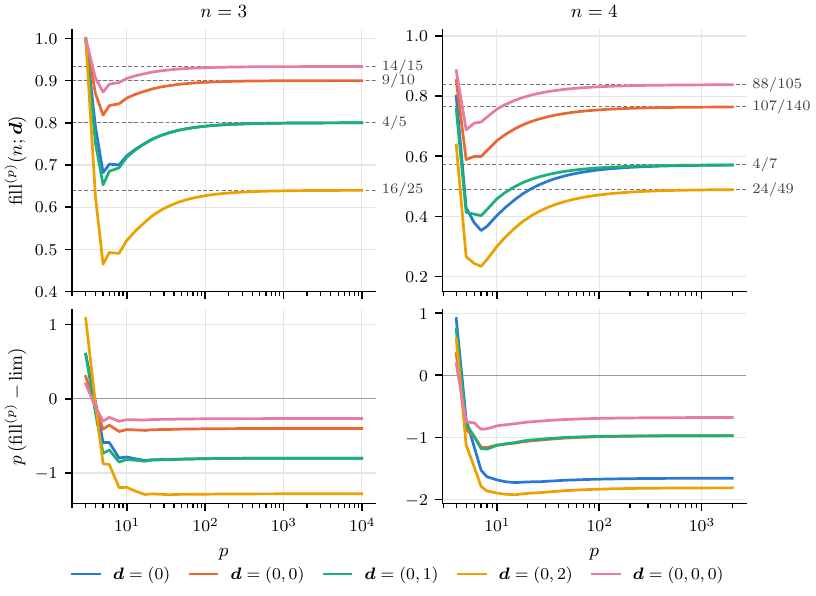}
\caption{Kupreyeva's limit (Example~\ref{ex-kupreyeva-values}): the fill rate
$\operatorname{fill}^{(p)}(n;\vek d)$ of \eqref{eq:fill-p} against the number $p$ of vertices,
on a logarithmic scale, for $n=3$ (left), $n=4$ (right) and five shift vectors, from exact counts.
Top: the dashed lines are the limits of Theorem~\ref{thm-kupreyeva-rank-r}, labelled on the right;
at $n=3$ the curves for $(0)$ and $(0,1)$, both with limit $q_3=4/5$, almost coincide. The left
ends, at $p=n$, are the fill rates $\operatorname{fill}(n;\vek d)$ of Section~\ref{sec-volume}.
Bottom: $p$ times the error, which settles to a negative constant: the approach is from below, at
rate $1/p$.}
\label{fig:kupreyeva-fill}
\end{figure}

\begin{mexample}
For \(n=3\), Theorem \ref{thm-kupreyeva-rank-r} gives the limits \(9/10\), \(4/5\) and \(16/25\) at
\(\vek d=(0,0)\), \((0,1)\) and \((0,2)\), and \(14/15\) and \(4/5\) at \((0,0,0)\) and
\((0,1,2)\). For \(n=4\) it gives \(107/140\), \(4/7\), \(24/49\) and \(16/49\) at \(\vek d=(0,0)\),
\((0,1)\), \((0,2)\) and \((0,3)\), and \(88/105\) at \((0,0,0)\). Among these only the product regime
gives \(q_n^r\). For \(n\ge3\) and positive gaps this is general: \(q(n,\delta+1)/q(n,\delta)=
2(n-\delta)/(2n-\delta)<1\), so \(q(n,\delta)>q_n\) for \(\delta<n-1\), and the limit exceeds \(q_n^r\)
unless every gap is at least \(n-1\). At \(\vek d=(0,1)\) the limit is \(q_n\) itself, because the shifted summand is
generically full except in its top degree, where Kruskal--Katona imposes nothing.

We compared these limits with the exact fill rates in eighteen cases: rank one for \(n=3,4,5\),
and fifteen shift vectors of rank \(2\) to \(4\) for \(n=3,4\), up to \(p=30000\) for \(n=3\),
\(p=2000\) for \(n=4\) and \(p=120\) for \(n=5\). In every case
\(p\bigl(\operatorname{fill}^{(p)}(n;\vek d)-\lim\bigr)\) settles to a constant, so the error
appears to be of order \(1/p\); Figure \ref{fig:kupreyeva-fill} shows ten of the cases. In all
eighteen the approach is from below, and slow. It need not be: at \(n=4\) and
\(\vek d=(0,1,1)\), not among them, it is from above. At rank one the constant is about \(-0.80\) for \(n=3\), \(-1.66\) for \(n=4\) and
\(-2.19\) for \(n=5\). For comparison, at \(p=n=5\) the fraction is \(0.43\) and at \(p=10\) it is
\(0.18\), against \(q_5\approx0.38\); the asymptotic regime starts later.
\label{ex-kupreyeva-values}
\end{mexample}

\begin{remark}
(\emph{The limits in \(p\) and in \(n\).}) Since \((2n-1)!!=(2n)!/(2^nn!)\), the value \eqref{eq:qn} is
\(q_n=2^{n-1}(n+1)/\binom{2n}n\), and Stirling's formula, through \(\binom{2n}n\sim4^n/\sqrt{\pi n}\),
gives \(q_n\sim(n+1)\sqrt{\pi n}/2^{n+1}\). The limits \(p\to\infty\) and \(n\to\infty\) do not
commute in any simple way: the fill rate at \(p=n\) falls much faster (Table \ref{tab:fill}: \(\operatorname{fill}(8)\approx0.027\) against
\(q_8=64/715\approx0.090\)).
\label{rem-kupreyeva-limit}
\end{remark}
\section{Ehrhart theory}
\label{sec-ehrhart}
For a lattice polytope \(P\) of dimension \(D\), the Ehrhart polynomial \(L_P\) and the
\(h^*\)-polynomial \(h^*_P\) are given by
\begin{equation}\label{eq:ehrhart}
  L_P(t)=\#\bigl(tP\cap\mathbb Z^{N+1}\bigr),\qquad
  \sum_{t\ge0}L_P(t)z^t=\frac{h^*_P(z)}{(1-z)^{D+1}} .
\end{equation}
\(L_P\) is a polynomial of degree \(D\), whose leading coefficient is the relative volume of \(P\):
its volume normalised so that a fundamental domain of the lattice parallel to the affine hull of
\(P\) has volume \(1\). For the Kozlov polytopes this is \(\vol P\), as in Section \ref{sec-volume}.

The coefficients of \(h^*_P\) are non-negative. Write
\(L^\circ_P(t)=\#\bigl(t\relint P\cap\mathbb Z^{N+1}\bigr)\) for the number of lattice points in the
relative interior of \(tP\). \emph{Ehrhart--Macdonald reciprocity} says that
\begin{equation}\label{eq:reciprocity}
  L_P(-t)=(-1)^DL^\circ_P(t)\quad(t\ge1),\qquad
  \sum_{t\ge1}L^\circ_P(t)z^t=\frac{z^{D+1}h^*_P(1/z)}{(1-z)^{D+1}} .
\end{equation}
The \textbf{codegree} of \(P\) is the least \(t\ge1\) for which \(tP\) has an interior lattice point,
that is, with \(L^\circ_P(t)>0\). Let \(s=\deg h^*_P\). The lowest term on the right of
\eqref{eq:reciprocity} is \(h^*_sz^{D+1-s}\), so the codegree is \(D+1-s\), and the number of
interior lattice points of that dilate is the top coefficient \(h^*_s\). For all of this see
\autocite[Chapters 3 and 4, in particular Theorems 4.1, 4.4 and 4.5]{BeckRobins2015}. We call \(P\)
\textbf{Gorenstein} if \(h^*_P\) is palindromic, \(h^*_i=h^*_{s-i}\); by a theorem of Stanley this is
equivalent to the Ehrhart ring of \(P\) being a Gorenstein ring \autocite{Stanley1996}.
\subsection{The Kozlov simplex}
\label{sec:orgb616b5a}

\begin{theorem}
\leanmark{\leannone}Write \(c_k=\binom nk\) (or take any integral leg vector). For residue vectors \(\vek z=(z_2,\dots,z_n)\) with
\(0\le z_k<c_k\), put \(\rho_1=0\), \(\rho_k=z_k/c_k\), and
\(\asc(\vek z)=\#\lbrace k\in[1,n-1]:\rho_k<\rho_{k+1}\rbrace\). Then
\(h^*_i(\Koz(n))=\#\lbrace\vek z:\asc(\vek z)=i\rbrace\).
\label{thm-hstar}
\end{theorem}

This is the specialisation to \(\Koz(n)\) of Savage and Schuster's formula for the \(h^*\)-polynomial
of a lecture hall polytope as an \(\vek s\)-Eulerian polynomial \autocite{SavageSchuster2012}, in the
form in which it is quoted in \autocite{HibiOlsenTsuchiya2018,SavageVisontai2015}; we include the
short proof, which is the standard fundamental-parallelepiped count.

\begin{proof}
Cone over the projected simplex in \(\mathbb R\times\mathbb R^{n-1}\), with generators
\(w_j=(1,c_2,\dots,c_j,0,\dots,0)\), \(j=1,\dots,n\). For a lattice simplex, \(h^*_i\) is the number
of lattice points at height \(i\) in \(\Pi=\lbrace\sum\lambda_jw_j:0\le\lambda_j<1\rbrace\)
\autocite[Chapter 3]{BeckRobins2015}. A point \(\sum\lambda_jw_j\) has coordinates \(x_k=c_k\mu_k\) with
\(\mu_k=\sum_{j\ge k}\lambda_j\) for \(k\ge2\), and height \(t=\lambda_1+\mu_2\); it is integral iff
\(t\in\mathbb Z\) and \(\mu_k\in\rho_k+\mathbb Z\), where \(\rho_k=(x_k\bmod c_k)/c_k\). Conversely the
residues determine the point: \(\lambda_k=\mu_k-\mu_{k+1}\in[0,1)\) with \(\mu_{n+1}=0\) forces
\(\mu_n=\rho_n\) and, descending, \(\mu_k=\rho_k+m_k\) with
\(m_k=m_{k+1}+[\rho_k<\rho_{k+1}]\); and \(\lambda_1=t-\mu_2\in[0,1)\) forces
\(t=\lceil\mu_2\rceil=m_2+[\rho_2>0]=\asc(\vek z)\). So \(x\mapsto\vek z\) is a bijection from
\(\Pi\cap\mathbb Z^n\) onto \(\prod_k\mathbb Z/c_k\), and the height is \(\asc(\vek z)\).
\end{proof}

\begin{remark}
(\emph{A misprint in the literature.}) Some restatements of Theorem \ref{thm-hstar}, for instance
\autocite[Lemma 1.1]{HibiOlsenTsuchiya2018}, count ascents of the residue vector after an auxiliary
entry \(e_0\) with \(s_0=1\), whose ratio \(e_0/s_0\) plays the part of our \(\rho_1\). They print the
convention \(e_0=1\) where \(e_0=0\), that is \(\rho_1=0\), is meant; with \(e_0=1\) one gets
\((14,10)\) instead of \((1,16,7)\) at \(n=4\), with \(\vek s=(6,4,1)\). (In the order \((1,4,6)\) of
Proposition \ref{prop-lecture-hall} the two conventions agree.)
\label{rem-hstar-misprint}
\end{remark}

\begin{mexample}
At \(n=3\) (dimension \(2\)), \((c_2,c_3)=(3,1)\): \(z_3=0\), and \(\asc=[0<\rho_2]\), so \(h^*=(1,2)\) and
\(L_{\Koz(3)}(t)=(t+1)(3t+2)/2=1,5,12,22,35,\dots\). At \(n=4\), \((c_2,c_3,c_4)=(6,4,1)\), \(24\) residue vectors:
\(h^*=(1,16,7)\) and \(L_{\Koz(4)}(t)=(t+1)^2(4t+1)\). The \(h^*\)-vectors for \(5\le n\le8\), of dimension \(n-1\), are
\[\small
\begin{array}{ll}
n=5: & (1,157,308,34)\\
n=6: & (1,1975,14184,10253,587)\\
n=7: & (1,60986,1054938,2043353,607906,14391)\\
n=8: & (1,3576516,143878949,624510000,523655643,80268620,939711).
\end{array}
\]
The sums are \(\NV(\Koz(n))=\prod_{k\ge2}\binom nk\), and \(h^*_1=\lvert\mathcal L_n\rvert-n\). The
sequence \(\lvert\mathcal L_n\rvert=1,2,5,20,162,1981,60993,3576524\) (\(n=1,\dots,8\)), the number of
lattice points of the Kozlov simplex, was not in the OEIS \autocite{OEIS} when checked in October 2026.
\label{ex-hstar}
\end{mexample}

\begin{lemma}
\leanmark{\leannone}Let \(n\ge3\) and \(c_k=\binom nk\). For \(2\le j\le n\) let \(y_j=1\), let
\(y_k=\lfloor c_ky_{k+1}/c_{k+1}\rfloor+1\) for \(k=j-1,\dots,2\), downward, and let \(y_k=0\) for
\(k>j\).
\begin{enumerate}
\item If \(j\le n-1\), the projection of \(\Koz(n)\) to the coordinates \(x_0,\dots,x_j\) is the
\((j-1)\)-simplex \(\lbrace1=\mu_1\ge\mu_2\ge\dots\ge\mu_j\ge0\rbrace\), and
\((1,n,y_2,\dots,y_n)\) is a lattice point of \(\Koz(n)\) whose projection lies in the relative
interior of that simplex.
\item If \(j=n\), then for every \(m\ge2\) the two lattice points \((m,mn,y_2,y_3,\dots,y_n)\) and
\((m,mn,y_2+1,y_3,\dots,y_n)\) lie in the relative interior of \(m\Koz(n)\).
\end{enumerate}
\label{lem-interior-points}
\end{lemma}

\begin{proof}
In both parts \(\mu_j=y_j/c_j=1/c_j>0\); \(\mu_k>\mu_{k+1}\) for \(2\le k<j\) by the choice of
\(y_k\); and \(\mu_k\le\mu_{k+1}+1/c_k\), so \(\mu_2\le1/c_j+\sum_{k=2}^{j-1}1/c_k\le1/c_j+(j-2)/n\),
since \(c_k\ge n\) for \(1\le k\le n-1\).
(1) The projection sends \(\tilde F_1,\dots,\tilde F_j\) to their truncations and
\(\tilde F_{j+1},\dots,\tilde F_n\) to that of \(\tilde F_j\); by Proposition \ref{prop-rank-one-H} the
hull is the chain set stated, with barycentric coordinates \(\mu_k-\mu_{k+1}\) (\(\mu_{j+1}=0\)), so its
relative interior is \(1>\mu_2>\dots>\mu_j>0\). Here \(c_j\ge n\), so \(\mu_2\le(j-1)/n<1\). Setting
\(\mu_{j+1}=\dots=\mu_n=0\) keeps the chain, so the point lies in \(\Koz(n)\).
(2) By Proposition \ref{prop-rank-one-H}, the relative interior of \(m\Koz(n)\) is the set of points
with \(x_0=m\), \(x_1=mn\) and \(m>\mu_2>\dots>\mu_n>0\). For the first point \(c_n=1\), so
\(\mu_2\le1+(n-2)/n\). For the second, \(\mu_2\) is larger by \(1/c_2=2/(n(n-1))\), which keeps
\(\mu_2>\mu_3\), and \(1+(n-2)/n+2/(n(n-1))<2\), an inequality equivalent to \(n>2\). So both points
satisfy the strict chain, with \(\mu_2<2\le m\).
\end{proof}

\begin{proposition}
\leanmark{\leannone}Let \(n\ge3\).
\begin{enumerate}
\item \(h^*_{\Koz(n)}\) is real-rooted, hence log-concave and unimodal \autocite{SavageVisontai2015}.
\item \(\Koz(n)\) has no interior lattice point and \(2\Koz(n)\) has one: the codegree is \(2\) and
\(\deg h^*=n-2\). Consequently \((t+1)\) divides \(L_{\Koz(n)}(t)\).
\item \(\Koz(n)\) is not Gorenstein.
\item The fundamental-parallelepiped group \(\mathbb Z^n/\langle w_j\rangle\), with the cone generators
\(w_j\) of the proof of Theorem \ref{thm-hstar}, is \(\bigoplus_{k=2}^n\mathbb Z/\binom nk\).
\end{enumerate}
\label{prop-hstar-props}
\end{proposition}

\begin{proof}
(2) On the relative interior the valid inequalities \(x_n\ge0\) and \(x_{n-1}\le c_{n-1}=n\) are
strict, since neither is tight on all of \(\Koz(n)\); for a lattice point that means \(x_n\ge1\) and
\(x_{n-1}\le n-1\), contradicting \(x_n\le x_{n-1}/n\). So \(\Koz(n)\) has no interior lattice point,
while \(2\Koz(n)\) has one by Lemma \ref{lem-interior-points}(2). Reciprocity \eqref{eq:reciprocity}
gives \(L(-1)=\pm L^\circ(1)=0\).

(3) By Lemma \ref{lem-interior-points}(2), \(2\Koz(n)\) has at least two interior lattice points. By
\eqref{eq:reciprocity} their number is \(h^*_{n-2}\), whereas a Gorenstein polytope has palindromic \(h^*\), so
\(h^*_{n-2}=h^*_0=1\). (This is the case \(r=1\) of Corollary \ref{cor-gorenstein}.)
(4) The basis change \(w_1\), \(w_j-w_{j-1}=c_je_j\) diagonalises the lattice.
\end{proof}

The polynomials \(L_{\Koz(n)}\) have positive coefficients for \(n\le6\); for example
\[
  L_{\Koz(5)}(t)=\tfrac16(t+1)(125t^3+250t^2+105t+6).
\]

\begin{observation}
The Ehrhart polynomial of \(\Koz(7)\) is
\begin{equation*}
\begin{split}
  L_{\Koz(7)}(t)={}&1-594t-\tfrac{5397}8t^2+\tfrac{173313}{16}t^3\\
  &+\tfrac{414687}{16}t^4+\tfrac{324135}{16}t^5+\tfrac{84035}{16}t^6
\end{split}
\end{equation*}
(computed by the transfer matrix of Corollary \ref{cor-transfer} and independently by Normaliz), so
\(\Koz(7)\) is not Ehrhart positive. Lecture hall simplices that are not Ehrhart positive were first exhibited by Liu and Solus
\autocite[Theorem 4.3]{LiuSolus2019}, and a second family, \(\vek s=(a,\dots,a,a+1)\), by Caicedo,
Juhnke and Poullot \autocite[Theorem 3.5]{CaicedoJuhnkePoullot2025}; \(\Koz(7)\), with
\(\vek s=(21,35,35,21,7,1)\), belongs to neither family and is a natural new instance. By
Proposition \ref{prop-linear-coeff} below, no \(\Koz(7;\vek d)\) is Ehrhart positive.
\label{obs-positivity}
\end{observation}
\subsection{Lattice points and Ehrhart polynomials by transfer matrix}
\label{sec:org87e83bc}

Let \(A^{\mathcal L}[\tau][\iota]=\#\lbrace x\in\mathcal L_n:\tp(x)=\tau,\ \indeg(x)=\iota\rbrace\), and
build from it \(T^{\mathcal L}(\delta)\) and \(\vek u^{\mathcal L}\) exactly as Part II builds
\(T'_n(\delta)\) and \(\vek u'_n\) from the counts on \(S_n\)
\autocite[Definition 28 and Theorem 41]{SnellmanArticle2}:
\begin{equation}\label{eq:transfer}
\begin{aligned}
  T^{\mathcal L}(\delta)[\tau][\tau']&=\sum_{\iota>\tau-\delta}A^{\mathcal L}[\tau'][\iota]
  =\#\lbrace x\in\mathcal L_n:\tp(x)=\tau',\ \indeg(x)+\delta>\tau\rbrace,\\
  u^{\mathcal L}[\tau]&=\sum_\iota A^{\mathcal L}[\tau][\iota]=\#\lbrace x\in\mathcal L_n:\tp(x)=\tau\rbrace .
\end{aligned}
\end{equation}
In the next corollary \(\vek u^{\mathcal L}\) is a row vector and \(\vek 1\) the all-ones column vector
of length \(n\), both indexed by \(\tau\in[1,n]\) like the rows of \(T^{\mathcal L}(\delta)\).

\begin{corollary}
\leanmark{\leanfull}\(\#\bigl(\Koz(n;\vek d)\cap\mathbb Z^{N+1}\bigr)=\vek u^{\mathcal L}\,T^{\mathcal L}(\delta_2)\cdots
T^{\mathcal L}(\delta_r)\,\vek 1\), and, with \(T_0=T^{\mathcal L}(0)\),
\[
  L_{\Koz(n;\vek d)}(t)=\vek u^{\mathcal L}\,T_0^{\,t-1}\,T^{\mathcal L}(\delta_2)\,T_0^{\,t-1}\,T^{\mathcal L}(\delta_3)\cdots
  T^{\mathcal L}(\delta_r)\,T_0^{\,t-1}\,\vek 1 .
\]
\label{cor-transfer}
\end{corollary}

\begin{proof}
By Corollary \ref{cor-additivity}, lattice points are in bijection with chained \(\mathcal L_n\)-tuples,
and counting chained tuples by a transfer-matrix product is Part II's rearrangement
\autocite[Section 7]{SnellmanArticle2}, which uses nothing about \(S_n\). The dilate is
\(t\Koz(n;\vek d)=\Koz(n;\vek d^{\times t})\), with \(\vek d^{\times t}\) as in the proof of
Corollary \ref{cor-idp}.
\end{proof}

\(T_0\) is upper unitriangular: a non-zero entry needs \(\tau'\ge\tau\), because
\(\indeg\le\tp+1\); the diagonal entry counts only the vertex \(\tilde F_\tau\); and the superdiagonal
entry is \(T_0[\tau][\tau+1]=c_{\tau+1}\), the choices of \(x_{\tau+1}\in[1,c_{\tau+1}]\). The same holds
for Part II's \(T'_n(0)\) on its principal block \(\tau\in[1,n]\) (Part II, Lemma 44), since a full \(\tau\)-skeleton plus any \(m\le\binom n{\tau+1}\) faces of the
next size is a complex. We write
\begin{equation}\label{eq:T0}
  T_0=T^{\mathcal L}(0)=I+\mathcal N ,
\end{equation}
so that \(\mathcal N\) is strictly upper triangular, with superdiagonal \((c_2,\dots,c_n)\).

\begin{corollary}
\leanmark{\leanpartial}\(L_{\Koz(n)}(t)=\vek u^{\mathcal L}T_0^{\,t-1}\vek 1\), a polynomial of degree \(n-1\) with leading
coefficient \(\prod_{k\ge2}\binom nk/(n-1)!\). The count
\begin{equation}\label{eq:ar}
  a(r)=\bigl\lvert\FF(n;0^r)\bigr\rvert ,
\end{equation}
is a polynomial in \(r\) with the same degree and the same leading coefficient
\autocite[Theorem 36]{SnellmanArticle2}. Hence
\(\operatorname{fill}(n;0^r)=a(r)/L_{\Koz(n)}(r)\to1\) as \(r\to\infty\).
\label{cor-zero-shift}
\end{corollary}

\begin{proof}
By \autocite[Theorem 41]{SnellmanArticle2}, \(a(r)=\vek u'_n\,T'_n(0)^{r-1}\vek 1\), and by
Corollary \ref{cor-transfer}, \(L_{\Koz(n)}(t)=\vek u^{\mathcal L}T_0^{\,t-1}\vek 1\).
By \eqref{eq:T0}, \(T_0^{t-1}=\sum_k\binom{t-1}k\mathcal N^k\), and \(\mathcal N^{n-1}\) has the single entry
\(\prod_\tau\mathcal N[\tau][\tau+1]=\prod_{k\ge2}c_k\) at \((1,n)\). Both \(\vek u^{\mathcal L}\) and
Part II's \(\vek u'_n\), restricted to that block, have first entry \(1\), the only point with \(\tp=1\) being \((1,n,0,\dots,0)\),
and both superdiagonals are \((c_2,\dots,c_n)\).
\end{proof}

The statement about \(a(r)\), degree and leading coefficient included, is
\autocite[Theorem 36]{SnellmanArticle2}; what is new is the lattice-point count and the comparison.
Since \(\Koz(n;0^r)=r\Koz(n)\), the difference \(L_{\Koz(n)}(r)-a(r)\) counts the lattice points of
the dilate that are not ff-vectors. By Corollary \ref{cor-zero-shift} it is a polynomial in \(r\) of degree at most
\(n-2\), and it vanishes for \(n\le3\) by Corollary \ref{cor-small-n}. For \(n=4\) and \(n=5\) it
is therefore decided by four of its values.

\begin{corollary}
\leanmark{\leannone}For every \(r\ge1\),
\[
  \begin{alignedat}{2}
  n=4:&\qquad& a(r)&=(r+1)(4r^2+3r+1),\\
      &&       L_{\Koz(4)}(r)&=(r+1)^2(4r+1),\\
      &&       L_{\Koz(4)}(r)-a(r)&=2r(r+1);\\[4pt]
  n=5:&&       a(r)&=\tfrac16(r+1)(125r^3+45r^2+34r+6),\\
      &&       L_{\Koz(5)}(r)&=\tfrac16(r+1)(125r^3+250r^2+105r+6),\\
      &&       L_{\Koz(5)}(r)-a(r)&=\tfrac16r(r+1)(205r+71).
  \end{alignedat}
\]
\label{cor-deficiency}
\end{corollary}

\begin{proof}[Proof (computational)]
In the proof of Corollary \ref{cor-zero-shift}, \(T_0\) and \(T'_n(0)\) are upper unitriangular, so
\(a(r)\) and \(L_{\Koz(n)}(r)\) are \(\sum_{k=0}^{n-1}\binom{r-1}k w_k\) for constants \(w_k\): polynomials
in \(r\) of degree at most \(n-1\), for every \(r\ge1\), with equal leading coefficients. A polynomial
of degree at most \(n-1\) is determined by \(n\) values. At \(r=1,\dots,5\),
\(a=16,69,184,385,696\) and \(L_{\Koz(4)}=20,81,208,425,756\) for \(n=4\), and
\(a=70,627,2592,7385,16926\) and \(L_{\Koz(5)}=162,1108,3964,10355,22406\) for \(n=5\); the stated
polynomials take these values. The values were computed independently of the transfer matrices,
\(a(r)\) as the \(r\)-fold sumset of \(S_n\) listed by the Kruskal--Katona theorem and
\(L_{\Koz(n)}(r)\) as the number of integer points \(x\) with \(x_0=r\), \(x_1=nr\) and
\(r\ge x_2/\binom n2\ge\dots\ge x_n/\binom nn\ge0\) (Proposition \ref{prop-rank-one-H}); see
Section \ref{sec-methods}.
\end{proof}

At \(n=5\) the cubic through the differences at \(r=1,\dots,4\) predicts the fifth, \(5480\), as the
degree bound says it must.
\subsection{Rank two and higher: examples}
\label{sec:orgac5c9e7}

\begin{mexample}
Three Ehrhart polynomials at rank two, with the dimension \(D\) of each polytope from \eqref{eq:dim}:
\[\small
\begin{array}{lcl}
 & D & L(t)\\[2pt]
\Koz(3;0,1) & 3 & 1+5t+10t^2+6t^3\\
\Koz(3;0,2) & 4 & 1+5t+\tfrac{37}4t^2+\tfrac{15}2t^3+\tfrac94t^4\\
\Koz(4;0,1) & 4 & 1+12t+54t^2+95t^3+52t^4
\end{array}
\]
In each case the degree of the polynomial is the dimension of the polytope. The polytope
\(\Koz(3;0,2)\) is in the product regime, so its Ehrhart polynomial is the square of
\(L_{\Koz(3)}(t)=\tfrac12(t+1)(3t+2)\) (Proposition \ref{prop-product-regime} and Example \ref{ex-hstar}). Its
\(h^*\)-vector \((1,20,29,4)\) is not the square \((1,4,4)\) of \(h^*_{\Koz(3)}=(1,2)\): the Ehrhart
polynomial of a product is the product of the Ehrhart polynomials, but the \(h^*\)-polynomial is not.
The linear coefficients are \(5=2\cdot\tfrac52\) for the two polytopes with \(n=3\) and \(12=2\cdot6\)
for \(\Koz(4;0,1)\), twice those of \(L_{\Koz(3)}\) and \(L_{\Koz(4)}\), as Proposition \ref{prop-linear-coeff}
below requires.

The next table lists \(h^*\)-vectors at rank two and three, each with the dimension \(D\) of its
polytope:
\[\small
\begin{array}{lcl@{\qquad}lcl}
 & D & h^* & & D & h^*\\[2pt]
\Koz(3;0,0) & 2 & (1,9,2) & \Koz(4;0,0) & 3 & (1,77,107,7)\\
\Koz(3;0,1) & 3 & (1,18,17) & \Koz(4;0,1) & 4 & (1,209,773,265)\\
\Koz(3;0,2) & 4 & (1,20,29,4) & \Koz(4;0,2) & 5 & (1,380,2940,3008,391)\\
\Koz(3;0,1,2) & 4 & (1,90,318,113) & \Koz(4;0,3) & 6 & (1,393,3782,5702,1593,49)\\
 & & & \Koz(4;0,1,2) & 5 & (1,2052,22428,29456,5103).
\end{array}
\]
\label{ex-ehrhart-rank2}
\end{mexample}
\subsection{Linear and quadratic coefficients}
\label{sec:org7b21e94}

The proofs of the next two propositions stay inside the transfer matrix. They rest on one identity
between the matrices of Corollary \ref{cor-transfer}, indexed by \(\tau\in[1,n]\) (Lemma
\ref{lem-row-shift}), and on a polynomial form of the powers of \(T_0\). As noted after Corollary
\ref{cor-transfer}, \(T_0\) is upper unitriangular, \(T_0=I+\mathcal N\) by \eqref{eq:T0}. Put
\begin{equation}\label{eq:U}
  U(t)=\sum_{k=0}^{n-1}\binom{t-1}k\mathcal N^k=U_0+tU_1+t^2U_2+\cdots ,
\end{equation}
a matrix whose entries are polynomials in \(t\), with \(U(t)=T_0^{\,t-1}\) for \(t\ge1\); here
\(U_j\) is the matrix of the coefficients of \(t^j\). Since \(\binom{-1}k=(-1)^k\) and \(\mathcal N\)
is nilpotent, \(U_0=\sum_k(-\mathcal N)^k=T_0^{-1}\). Write \(\vek e_1\) for the first unit row
vector.

\begin{lemma}
\leanmark{\leanfull}For \(\delta\ge0\), row \(\tau\) of \(T^{\mathcal L}(\delta)\) is row \(\max(\tau-\delta,1)\) of \(T_0\).
Consequently \(\vek e_1T^{\mathcal L}(\delta)=\vek u^{\mathcal L}=\vek e_1T_0\) and
\(T^{\mathcal L}(\delta)\,T_0^{-1}\vek 1=\vek 1\).
\label{lem-row-shift}
\end{lemma}

\begin{proof}
Row \(\tau\) of \(T^{\mathcal L}(\delta)\) counts the \(x\in\mathcal L_n\) with \(\tp(x)=\tau'\) and
\(\indeg(x)>\tau-\delta\). Every \(x\in\mathcal L_n\) has \(x_0=1\) and \(x_1=n\), so \(\indeg(x)\ge2\),
and the condition \(\indeg(x)>\sigma\) is vacuous for \(\sigma\le1\). Hence the condition with
\(\tau-\delta\) is the condition with \(\max(\tau-\delta,1)\), which is that of row
\(\max(\tau-\delta,1)\) of \(T_0\). For \(\tau=1\) this gives \(\vek e_1T^{\mathcal L}(\delta)=\vek e_1T_0\),
and row \(1\) of \(T_0\), with its vacuous condition, is \(\vek u^{\mathcal L}\). Finally, entry
\(\tau\) of \(T^{\mathcal L}(\delta)T_0^{-1}\vek 1\) is entry \(\max(\tau-\delta,1)\) of
\(T_0T_0^{-1}\vek 1=\vek 1\).
\end{proof}

\begin{proposition}
\leanmark{\leanfull}\([t^1]\,L_{\Koz(n;\vek d)}(t)=r\cdot[t^1]\,L_{\Koz(n)}(t)\) for every \(\vek d\).
\label{prop-linear-coeff}
\end{proposition}

\begin{proof}
Substituting the matrix \(U(t)\) of \eqref{eq:U} for every \(T_0^{\,t-1}\) in the
display of Corollary \ref{cor-transfer} gives a polynomial that agrees with \(L_{\Koz(n;\vek d)}\) for
\(t\ge1\), hence equals it:
\(L_{\Koz(n;\vek d)}(t)=\vek u^{\mathcal L}U(t)A_2U(t)\cdots A_rU(t)\vek 1\) with
\(A_i=T^{\mathcal L}(\delta_i)\), and \(U_0=T_0^{-1}\).

The \(t^1\)-coefficient is the sum over \(j\in[1,r]\) of the product in which the \(j\)-th factor
\(U(t)\) is replaced by \(U_1\) and every other by \(U_0\). All gaps \(\delta_i\) are non-negative,
so Lemma \ref{lem-row-shift} applies to every \(A_i\). On the left,
\(\vek u^{\mathcal L}U_0=\vek e_1T_0T_0^{-1}=\vek e_1\) and \(\vek e_1A_i=\vek u^{\mathcal L}\), so
\(\vek u^{\mathcal L}U_0A_2U_0\cdots U_0A_j=\vek u^{\mathcal L}\). On the right, \(A_iU_0\vek 1=\vek 1\),
so \(A_{j+1}U_0\cdots A_rU_0\vek 1=\vek 1\). The \(j\)-th term is therefore
\(\vek u^{\mathcal L}U_1\vek 1\), the case \(r=1\), which is \([t^1]\,L_{\Koz(n)}(t)\), and there are
\(r\) terms.
\end{proof}

The same reduction gives the constant term:
\(L_{\Koz(n;\vek d)}(0)=\vek u^{\mathcal L}U_0A_2U_0\cdots A_rU_0\vek 1=\vek u^{\mathcal L}U_0\vek 1=\vek e_1\vek 1=1\),
as Ehrhart's theorem requires. Proposition \ref{prop-linear-coeff} also follows from a theorem of
McMullen \autocite{McMullen1977}, found independently by Bernstein, which we quote in the form of
\autocite[Theorem 2.3]{HaaseJuhnkeSanyalTheobald2017}.

\begin{theorem}
(\emph{McMullen}.) Let \(P_1,\dots,P_r\subseteq\mathbb R^M\) be lattice polytopes, of any dimensions. There is
a polynomial \(E(\lambda_1,\dots,\lambda_r)\) such that
\[
  \#\bigl((\lambda_1P_1+\dots+\lambda_rP_r)\cap\mathbb Z^M\bigr)=E(\lambda_1,\dots,\lambda_r)
  \qquad\text{for all }(\lambda_1,\dots,\lambda_r)\in\mathbb Z_{\ge0}^r,
\]
where \(0P_i=\lbrace0\rbrace\). Its degree in \(\lambda_i\) is \(\dim P_i\).
\label{thm-mcmullen}
\end{theorem}

The polytopes need not have the same dimension, and some \(\lambda_i\) may vanish. Setting
\(\lambda_l=0\) for \(l\ne i\) gives \(E(0,\dots,0,\lambda_i,0,\dots,0)=L_{P_i}(\lambda_i)\), so the part of
\(E\) of degree one is \(\sum_i\ell_i\lambda_i\), with \(\ell_i\) the linear coefficient of \(L_{P_i}\),
and the constant term is \(1\). Now take \(M=N+1\) and \(P_i=\shift^{d_i}\Koz(n)\), which has the
Ehrhart polynomial of \(\Koz(n)\) because \(\shift^{d_i}\) maps \(\mathbb Z^{n+1}\) bijectively onto the
lattice points of its image. Then \(E(t,\dots,t)=L_{\Koz(n;\vek d)}(t)\), whose linear coefficient is
\(r\ell_1\): this is Proposition \ref{prop-linear-coeff}. The argument given above avoids
Theorem \ref{thm-mcmullen}, and is the one checked in Lean (Section \ref{sec-methods}).

\begin{proposition}
\leanmark{\leannone}Let \(a_1=[t^1]\,L_{\Koz(n)}\), \(a_2=[t^2]\,L_{\Koz(n)}\) and, for \(\delta\ge0\),
\begin{equation}\label{eq:Qn}
  Q_n(\delta)=\vek u^{\mathcal L}U_1\,T^{\mathcal L}(\delta)\,U_1\vek 1,
\end{equation}
with \(U_1\) as in \eqref{eq:U}. For every \(\vek d\),
\[
  [t^2]\,L_{\Koz(n;\vek d)}(t)=r\,a_2+\sum_{1\le j<k\le r}Q_n(d_k-d_j).
\]
Moreover \(Q_n(\delta)=[t^2]\,L_{\Koz(n;0,\delta)}-2a_2\), \(Q_n(0)=2a_2\), and \(Q_n(\delta)=a_1^2\) for
\(\delta\ge n-1\).
\label{prop-quadratic-coeff}
\end{proposition}

\begin{proof}
Use the expansion \eqref{eq:U} to second order in the product
\[
  \vek u^{\mathcal L}U(t)A_2U(t)\cdots A_rU(t)\vek 1
\]
of the proof of Proposition \ref{prop-linear-coeff}. Its \(t^2\)-coefficient is the sum of the \(r\)
products with one factor \(U_2\) and \(U_0\) elsewhere, and of the \(\binom r2\) products with \(U_1\)
in positions \(j<k\) and \(U_0\) elsewhere. The two reductions of that proof,
\(\vek u^{\mathcal L}U_0A_2U_0\cdots U_0A_j=\vek u^{\mathcal L}\) on the left and
\(A_{k+1}U_0\cdots A_rU_0\vek 1=\vek 1\) on the right, turn a product of the first kind into
\(\vek u^{\mathcal L}U_2\vek 1=a_2\) (the case \(r=1\)), and one of the second kind into
\(\vek u^{\mathcal L}U_1M_{jk}U_1\vek 1\) with \(M_{jk}=A_{j+1}U_0A_{j+2}U_0\cdots U_0A_k\).

By Lemma \ref{lem-row-shift}, \(A_i=S(\delta_i)T_0\), where \(S(\delta)\) is the \(0/1\)-matrix whose row
\(\tau\) has its single \(1\) in column \(\max(\tau-\delta,1)\). Hence \(A_iU_0=S(\delta_i)\), and
\(S(\delta)S(\delta')=S(\delta+\delta')\) because \(\max(\max(\tau-\delta,1)-\delta',1)=\max(\tau-\delta-\delta',1)\).
So \(M_{jk}=S(\delta_{j+1})\cdots S(\delta_k)\,T_0=S(d_k-d_j)\,T_0=T^{\mathcal L}(d_k-d_j)\), and the
term is \(Q_n(d_k-d_j)\). The case \(r=2\) gives the first identity. At \(\delta=0\),
\(L_{\Koz(n;0,0)}(t)=L_{\Koz(n)}(2t)\) has \(t^2\)-coefficient \(4a_2\). For \(\delta\ge n-1\) every row of
\(S(\delta)\) is \(\vek e_1\), so \(T^{\mathcal L}(\delta)=\vek 1\,\vek e_1T_0=\vek 1\,\vek u^{\mathcal L}\)
and \(Q_n(\delta)=(\vek u^{\mathcal L}U_1\vek 1)^2=a_1^2\), as Proposition \ref{prop-product-regime}
requires.
\end{proof}

As for Proposition \ref{prop-linear-coeff}, Theorem \ref{thm-mcmullen} gives the same shape: the part of degree
two of the polynomial \(E\) of that theorem, with \(P_i=\shift^{d_i}\Koz(n)\), is
\(\sum_i\alpha_i\lambda_i^2+\sum_{j<k}\beta_{jk}\lambda_j\lambda_k\), with \(\beta_{jk}\) determined by the
pair \(P_j,P_k\): \(E\) is valid when some \(\lambda_i\) vanish, and deleting the first
\(d_j\) coordinates, the inverse of \(\shift^{d_j}\), carries the pair to \(\Koz(n)\),
\(\shift^{d_k-d_j}\Koz(n)\).

\begin{remark}
So the \(t^2\)-coefficient at every rank is read off rank two, through the \(n\) numbers
\(Q_n(0),\dots,Q_n(n-1)\):
\[\small
\begin{array}{rll}
n & a_2 & Q_n(0),\dots,Q_n(n-1)\\
3 & 3/2 & 3,\ 7,\ 25/4\\
4 & 9 & 18,\ 36,\ 36,\ 36\\
5 & 355/6 & 355/3,\ 949/4,\ 10231/36,\ 11021/36,\ 1369/4\\
6 & 1033/4 & 1033/2,\ 465/4,\ -91,\ -703/4,\ -2477/4,\ 1/4
\end{array}
\]
For \(3\le n\le5\) every entry is positive, so \(L_{\Koz(n;\vek d)}\) has positive coefficients of \(t^0\),
\(t^1\) and \(t^2\) for every \(\vek d\) (Proposition \ref{prop-linear-coeff}). At \(n=6\) one pair at
distance \(4\) outweighs two rank-one terms: the \(t^2\)-coefficient of \(L_{\Koz(6;0,4)}\) is
\(2\cdot\tfrac{1033}4-\tfrac{2477}4=-\tfrac{411}4\), so \(\Koz(6;0,4)\) is not Ehrhart positive, and
\(\Koz(6;0,2,4)\), with \(3\cdot\tfrac{1033}4-2\cdot91-\tfrac{2477}4=-\tfrac{53}2\), is another.
\label{rem-quadratic-coeff}
\end{remark}
\subsection{Codegree and the Gorenstein property}
\label{sec:orgf3e9535}

\begin{lemma}
\leanmark{\leannone}Let \(A,B\subseteq\mathbb R^M\) be polytopes, \(V_A\) the linear space parallel to the affine hull of
\(A\), and \(\pi\colon\mathbb R^M\to\mathbb R^M/V_A\) the quotient map. If \(a\in\relint A\), \(b\in B\)
and \(\pi(b)\in\relint\pi(B)\), then \(a+b\in\relint(A+B)\). In particular
\(\relint A+\relint B\subseteq\relint(A+B)\), and \(\relint A+B\subseteq\relint(A+B)\) when the affine hull of
\(B\) is parallel to a subspace of \(V_A\).
\label{lem-relint-sum}
\end{lemma}

\begin{proof}
A point \(x\) of a polytope \(C\) lies in \(\relint C\) if and only if every linear functional \(\ell\) with
\(\ell(x)=\max_C\ell\) is constant on \(C\). If \(x\in\relint C\) and \(\ell(y)<\ell(x)\) for some \(y\in C\),
then \(x+\varepsilon(x-y)\in C\) for small \(\varepsilon>0\), where \(\ell\) is larger; conversely, a point of
the relative boundary lies in a proper face, which is the set of maximisers of a functional that is
not constant on \(C\) \autocite[Chapter 2]{Ziegler1995}. Now let \(\ell(a+b)=\max_{A+B}\ell\). Since
\(\max_{A+B}\ell=\max_A\ell+\max_B\ell\), we have \(\ell(a)=\max_A\ell\) and \(\ell(b)=\max_B\ell\). The first
makes \(\ell\) constant on \(A\), so \(\ell\) vanishes on \(V_A\) and equals \(\bar\ell\circ\pi\) for a linear
\(\bar\ell\). The second says that \(\bar\ell\) attains its maximum over \(\pi(B)\) at \(\pi(b)\), so
\(\bar\ell\) is constant on \(\pi(B)\) and \(\ell\) on \(B\). Hence \(\ell\) is constant on \(A+B\). If
\(b\in\relint B\), the same argument shows \(\pi(b)\in\relint\pi(B)\); if the affine hull of \(B\) is
parallel to a subspace of \(V_A\), then \(\pi(B)\) is a point.
\end{proof}

The easy half of the codegree is immediate from Lemma \ref{lem-relint-sum}.

\begin{proposition}
\leanmark{\leannone}Let \(n\ge3\). The codegree of \(\Koz(n;\vek d)\) is at most \(2\), and equals \(2\) when \(r=1\) or
\(d_r\ge d_{r-1}+2\).
\label{prop-codegree}
\end{proposition}

\begin{proof}
By Lemma \ref{lem-relint-sum}, \(\relint(2\Koz(n;\vek d))\supseteq\sum_i\relint(2\shift^{d_i}\Koz(n))\), and
each summand contains a lattice point by Lemma \ref{lem-interior-points}(2). If \(d_r\ge d_{r-1}+2\), the coordinates \(N-1,N\) belong
to the last summand only, and the argument of Proposition \ref{prop-hstar-props}(2) applies to them.
\end{proof}

\begin{theorem}
\leanmark{\leannone}Let \(n\ge3\), and call an index \(i\in[1,r]\) \textbf{exposed} if \(d_i\) is not repeated from below
(\(i=1\) or \(d_{i-1}<d_i\)) and the gap after it is large (\(i=r\) or \(d_{i+1}-d_i\ge n-1\)).
\begin{enumerate}
\item The codegree of \(\Koz(n;\vek d)\) is \(2\) if some index is exposed, and \(1\) otherwise.
\item If \(c\) is the codegree, \(c\Koz(n;\vek d)\) has at least two lattice points in its relative
interior.
\end{enumerate}
In particular, if every gap is at most \(n-2\), the codegree is \(1\) if and only if \(r\ge2\) and
\(d_{r-1}=d_r\).
\label{thm-codegree}
\end{theorem}

Equivalently, in terms of the maximal runs of \(\vek d\) with gaps at most \(n-2\) (the factors of the
lattice product described after Theorem \ref{thm-staircase}), the codegree is \(1\) if and only if every
run ends in a repeated value.

\begin{proof}
Write \(P=\Koz(n;\vek d)\) and \(P_i=\shift^{d_i}\Koz(n)\). We use, as in the proof of Proposition
\ref{prop-hstar-props}(2), that a valid inequality that is not tight on all of a polytope is strict on its
relative interior.

\emph{An exposed index.} Let \(i\) be exposed and \(k=d_i+n\). For \(j<i\), \(d_j<d_i\), so summand \(j\) has
own index \(k-d_j>n\) and contributes \(0\) to \(x_k\); for \(j>i\), \(d_j\ge d_i+n-1\), so its own index
is at most \(1\) and it contributes a constant. Hence \(x_k=\beta+(\vek p_i)_n\) with \(\beta\in\mathbb Z\)
constant on \(P\), and \((\vek p_i)_n\) takes the values \(0\) and \(1\) (at \(\tilde F_1\) and
\(\tilde F_n\)). So \(\beta<x_k<\beta+1\) on \(\relint P\), which has no lattice point, and the codegree
is \(2\) by Proposition \ref{prop-codegree}.

\emph{No exposed index.} For an integer \(T\), or \(T=\infty\), let \(\tau_T\) project \(\mathbb R^{N+1}\) to the
coordinates \(x_0,\dots,x_T\), and put \(P_{\le s}=P_1+\dots+P_s\) and \(P_X=\sum_{i\in X}P_i\). Call
\(i\le s\) exposed at \((s,T)\) if \(d_i\) is not repeated from below and either \(i<s\) and
\(d_{i+1}-d_i\ge n-1\), or \(i=s\) and \(T-d_s\ge n\); exposed at \((r,\infty)\) is exposed. We show by
induction on \(s\): \emph{if \(T\ge d_s+2\) (no condition if \(s=0\)) and no index is exposed at \((s,T)\), there is a lattice point
\(b\in P_{\le s}\) with \(\tau_T(b)\in\relint\tau_T(P_{\le s})\).} For \(s=0\) take \(b=0\). Let \(s\ge1\).

If \(T-d_s\le n-1\), let \(k\) be the largest index with \(d_k\le T-n\) (\(k=0\) if there is none) and
\(A=\lbrace k+1,\dots,s\rbrace\). For \(i\in A\), \(2\le T-d_i\le n-1\), so \(\tau_T(P_i)\) is \(\shift^{d_i}\)
of the projection of \(\Koz(n)\) to its first \(T-d_i+1\) coordinates, and Lemma
\ref{lem-interior-points}(1) gives a lattice point \(\vek y_i\in P_i\) with
\(\tau_T(\vek y_i)\in\relint\tau_T(P_i)\). By Lemma \ref{lem-relint-sum}, \(a=\sum_{i\in A}\vek y_i\) has
\(\tau_T(a)\in\relint\tau_T(P_A)\). The linear space parallel to \(\tau_T(P_A)\) is
\(\sum_{i\in A}\mathbb R^{[d_i+2,T]}=\mathbb R^{[d_{k+1}+2,T]}\), and the quotient by it is \(\tau_{T'}\) with
\(T'=d_{k+1}+1\). If \(k\ge1\), then \(d_k\le T-n<d_{k+1}\), so \(T'\ge d_k+2\), and \(k\) is exposed at
\((k,T')\) exactly when \(d_{k+1}-d_k\ge n-1\) and \(d_k\) is not repeated from below, that is, when
it is exposed at \((s,T)\); the indices below \(k\) are unaffected. By induction there is a lattice
point \(b\in P_{\le k}\) with \(\tau_{T'}(b)\in\relint\tau_{T'}(P_{\le k})\), and Lemma
\ref{lem-relint-sum} gives \(\tau_T(a+b)\in\relint\tau_T(P_{\le s})\).

If \(T-d_s\ge n\), then, \(s\) not being exposed, \(s\ge2\) and \(d_{s-1}=d_s\). Let \(m\ge2\) be the number
of \(i\le s\) with \(d_i=d_s\), \(k=s-m\) and \(A=\lbrace k+1,\dots,s\rbrace\). Here \(\tau_T\) is injective on
\(P_A=m\shift^{d_s}\Koz(n)\), which has a lattice point \(a\) in its relative interior by Lemma
\ref{lem-interior-points}(2). The quotient by \(\mathbb R^{[d_s+2,d_s+n]}\) keeps the coordinates up to
\(d_s+1\) and those above \(d_s+n\), which vanish on \(P_{\le k}\); on \(\tau_T(P_{\le k})\) it is
therefore \(\tau_{T'}\) with \(T'=d_s+1\), up to these zero coordinates. If \(k\ge1\), then
\(T'\ge d_k+2\), and \(k\) is exposed at \((k,T')\) exactly when it is exposed at \((s,T)\), since
\(d_{k+1}=d_s\). Induction and Lemma \ref{lem-relint-sum} finish as before.

At \((s,T)=(r,\infty)\) this gives a lattice point in \(\relint P\), so the codegree is \(1\).

\emph{Two points.} If the codegree is \(1\), the induction at \((r,\infty)\) starts in the second case, where
Lemma \ref{lem-interior-points}(2) offers two choices of \(a\); with the same \(b\) they give two distinct
interior lattice points. If it is \(2\), choose for each \(i\) a lattice point \(\vek q_i\) in
\(\relint(2P_i)\) (Lemma \ref{lem-interior-points}(2), shifted); \(\sum_i\vek q_i\in\relint(2P)\) by Lemma
\ref{lem-relint-sum}, and replacing \(\vek q_r\) by a second such point gives another.

If every gap is at most \(n-2\), only \(i=r\) can be exposed, and it is unless \(r\ge2\) and
\(d_{r-1}=d_r\).
\end{proof}

\begin{mexample}
\(\Koz(3;0,2,2)\) has \(d_2=d_3\), but its gap \(2=n-1\) makes index \(1\) exposed, and its codegree is
\(2\). The projection to \((x_2,x_3)\) and \((x_4,x_5)\) is a unimodular equivalence onto a translate
of \(\Koz(3)\times2\Koz(3)\), by the remark on runs after Theorem \ref{thm-staircase}: the first summand
contributes \(x_2-2\) and \(x_3-6\), the other two contribute \(x_4\) and \(x_5\). So the condition
\(d_{r-1}=d_r\) alone does not decide the codegree: \(\Koz(3;0,0,2)\) is the same product, with
\(d_2\ne d_3\), while \(\Koz(3;0,0,2,2)\), with the same last pair as \(\Koz(3;0,2,2)\), has no exposed
index and codegree \(1\). \(\Koz(3;0,2,2)\) and \(\Koz(3;0,0,2)\) both have \(L(t)=L_{\Koz(3)}(t)\,L_{\Koz(3)}(2t)=\tfrac12(t+1)(3t+2)(2t+1)(3t+1)\), dimension
\(4\), \(h^*=(1,55,130,30)\), and \(30=2\cdot15\) interior lattice points in the second dilate, as
\eqref{eq:reciprocity} requires of the top coefficient.
\label{ex-codegree-product}
\end{mexample}

\begin{corollary}
\leanmark{\leannone}For \(n\ge3\), no Kozlov polytope \(\Koz(n;\vek d)\) is Gorenstein: \(h^*_s\ge2\), \(s=\deg h^*\).
\label{cor-gorenstein}
\end{corollary}

\begin{proof}
By \eqref{eq:reciprocity} \(h^*_s\) is the number of lattice points in the relative interior of
\(c\Koz(n;\vek d)\), \(c\) the codegree, which is at least \(2\) by Theorem \ref{thm-codegree}(2); and
\(h^*_0=1\).
\end{proof}

\begin{observation}
In all \(42\) computed cases (\(n\le7\), ranks \(2\) to \(4\)) the \(h^*\)-polynomial of \(\Koz(n;\vek d)\)
is real-rooted, so its \(h^*\)-vector is log-concave and unimodal, as at rank one (Proposition
\ref{prop-hstar-props}). In the \(39\) of them with \(n\ge3\) it is not palindromic, as Corollary \ref{cor-gorenstein} requires;
the three with \(n=2\) have \(h^*=(1,1)\), which is. Ehrhart
positivity fails at \(n=6\) as well: \(L_{\Koz(6;0,4)}(t)\) has \(t^2\)-coefficient \(-411/4\), which is \(2a_2+Q_6(4)\)
(Remark \ref{rem-quadratic-coeff}), although
\(L_{\Koz(6)}\) and \(L_{\Koz(6)}^2=L_{\Koz(6;0,5)}\) are positive (by the transfer matrix, and for
\(\Koz(6;0,4)\) independently by Barvinok's algorithm).
\label{obs-hstar-unimodal}
\end{observation}
\section{The centroid of the Hilbert functions and the centroid of the polytope}
\label{sec-centroid}
Three uniform measures live on \(\Koz(n;\vek d)\): on the ff-vectors \(\FF(n;\vek d)\) (which is
\(S_n\) at rank one), on all lattice points, and the continuous one. For the question ``where do the
Hilbert functions sit inside their convex hull'' the first is the natural one, with each distinct
Hilbert function counted once. At rank one we work in the LYM coordinates
\(\mu_k=x_k/\binom nk\) of \eqref{eq:mu}, in which the barycentric coordinates are \(\lambda_j=\mu_j-\mu_{j+1}\), so that
``where a centroid sits relative to the vertices'' and ``its \(\mu\)-profile'' are the same information.
At rank \(r\) we use
\begin{equation}\label{eq:nu}
  \nu_k=\frac{x_k}{\dim F_k},\qquad \dim F_k=\sum_{i=1}^r\binom n{k-d_i},
\end{equation}
the fraction of the degree-\(k\) part of \(F\) that survives in the quotient \(F/M\); at rank one
\(\nu_k=\mu_k\). These are not barycentric coordinates: for \(r\ge2\) the polytope
\(\Koz(n;\vek d)\) is not a simplex, each \(\nu_k\) lies in \([0,1]\), and their sum is not \(1\). As
at rank one, they say how full the quotient is in each degree. We write
\begin{equation}\label{eq:centroid}
  c(X)=\frac1{\lvert X\rvert}\sum_{x\in X}x\quad(X\text{ finite}),\qquad
  c(X)=\frac1{\vol X}\int_Xx\,dx\quad(X\text{ a polytope})
\end{equation}
for the centroid of the uniform measure on \(X\), the integral taken with Lebesgue measure on the
affine hull of \(X\). It is not to be confused with the leg vector \(\vek c\).

Table \ref{tab:centroid-summary} sorts the results of this section and of Section \ref{sec-moments} by
polytope: \(\Koz(n)\) itself, its dilates \(r\Koz(n)=\Koz(n;0^r)\), and \(\Koz(n;\vek d)\) in general.
Equal shifts \(\vek d=(d,\dots,d)\) give \(\shift^d\bigl(r\Koz(n)\bigr)\), a copy of the dilate.

\begin{table}[tbp]
\centering\footnotesize
\renewcommand{\arraystretch}{1.15}
\begin{tabular}{@{}p{0.12\linewidth}p{0.42\linewidth}p{0.39\linewidth}@{}}
\toprule
polytope & centroid & spread\\
\midrule
\(\Koz(n)\)
& \raggedright All three located: the polytope's exactly (Prop.~\ref{prop-centroid-simplex}). The
  \(\beta\)-symmetry holds for the polytope (Thm.~\ref{thm-beta-continuous}), for the lattice points
  iff \(n\) is prime (Thm.~\ref{thm-beta-lattice}, Obs.~\ref{obs-beta-composite}), and for the
  \(f\)-vectors up to \(\varepsilon_n\) (Thm.~\ref{thm-alexander-core}), which sit one vertex lower.
  Table~\ref{tab:centroids-rank1}, Figure~\ref{fig:mu}.
& \raggedright The polytope exactly (Prop.~\ref{prop-second-moments}). \(f\)-vectors, lattice points and polytope
  compared by computation only (Obs.~\ref{obs-concentration},
  Tables~\ref{tab:concentration}--\ref{tab:profile-n24}).\tabularnewline[3pt]
\(r\Koz(n)\)
&\raggedright The ff-vectors' centroid, scaled by \(1/r\), tends to the polytope's with error \(O(1/r)\)
  (Thm.~\ref{thm-centroid-dilation}(1)), exactly so at \(n=3\) (Thm.~\ref{thm-n3-equal-shifts}). The
  lattice-point centroid of \(t\Koz(n)\) is known exactly for \(n\le5\) (Prop.~\ref{prop-dilates}).
  Figure~\ref{fig:mu-rank-r}.
& \raggedright The covariance, scaled by \(1/r^2\), tends to the polytope's with error \(O(1/r)\)
  (Thm.~\ref{thm-centroid-dilation}(2)), and the whole measure converges weakly
  (Rem.~\ref{rem-weak-convergence}). The polytope's covariance grows like \(r^2\), the independent
  sum's like \(r\) (Prop.~\ref{prop-pushforward}(3)).\tabularnewline[3pt]
\(\Koz(n;\vek d)\)
& \raggedright Five centroids, which pair off in the product regime (Prop.~\ref{prop-pushforward}(1),(2));
  outside it the identities fail in every case computed (Obs.~\ref{obs-pushforward-fails},
  Ex.~\ref{ex-centroid-rank2}).
& \raggedright The same for the covariances (Prop.~\ref{prop-pushforward}(2), Obs.~\ref{obs-pushforward-fails}).
  The ff-vectors are more spread than the polytope in all \(21\) cases computed
  (Obs.~\ref{obs-moments-rank-r}, Table~\ref{tab:rank-r-variance}). At a constant positive gap
  nothing is proved (Section~\ref{sec-open}).\tabularnewline
\bottomrule
\end{tabular}
\caption{What Sections~\ref{sec-centroid} and~\ref{sec-moments} show, polytope by polytope.}
\label{tab:centroid-summary}
\end{table}
\subsection{Rank one: exact symmetries}
\label{sec:orgb28d9fc}

\begin{proposition}
\leanmark{\leanfull}The centroid of \(\Koz(n)\) has \(\mu_k=(n-k+1)/n\); all its barycentric coordinates equal \(1/n\),
and its mean vertex index \(\sum_jj\lambda_j=\sum_k\mu_k\) is \((n+1)/2\).
\label{prop-centroid-simplex}
\end{proposition}

\begin{proof}
The centroid of a simplex is the average of its vertices, and \(\tilde F_j\) has \(\mu_k=1\) iff
\(k\le j\).
\end{proof}

\begin{lemma}
\leanmark{\leanfull}Let \(n\ge3\). The affine map \(\beta\) of \(\mathbb R^{n+1}\) fixing \(x_0,x_1,x_{n-1},x_n\) and, in LYM
coordinates, sending \(\mu_k\mapsto1+\mu_{n-1}-\mu_{n-k}\) for \(2\le k\le n-2\), is an involution of
\(\Koz(n)\). It permutes the vertices, \(\tilde F_j\leftrightarrow\tilde F_{n-1-j}\) for
\(1\le j\le n-2\), and fixes \(\tilde F_{n-1}\) and \(\tilde F_n\). On the codimension-two face
\(\Koz(n)^\circ=\Koz(n)\cap\lbrace\mu_{n-1}=0\rbrace=\conv(\tilde F_1,\dots,\tilde F_{n-2})\) it is
\(x_k\mapsto\binom nk-x_{n-k}\), the Alexander duality map \(D_n\) of \autocite[Definition 13]{SnellmanArticle2}.
\label{lem-beta}
\end{lemma}

\begin{proof}
Since \(\mu_{n-1}\) is fixed, applying the map twice returns \(\mu_k\). The vertex \(\tilde F_j\),
\(j\le n-2\), has \(\mu_{n-1}=0\) and \(\mu_k=[k\le j]\), and its image has
\(\mu_k=1-[n-k\le j]=[k\le n-1-j]\); \(\tilde F_{n-1}\) and \(\tilde F_n\) have \(\mu_k=1\) for
\(k\le n-1\), and are fixed. An affine map permuting the vertices of a simplex maps it onto itself.
\end{proof}

\begin{theorem}
\leanmark{\leanfull}For \(n\ge3\) the uniform measure on \(\Koz(n)\) is \(\beta\)-invariant. Hence its centroid satisfies
\(\mu_k+\mu_{n-k}=1+\mu_{n-1}\) for \(1\le k\le n-1\), and its covariance matrix is \(\beta\)-invariant:
\(\bar\beta\Sigma\bar\beta^{\mathsf T}=\Sigma\), where \(\Sigma\) is the covariance matrix in the coordinates
\((x_0,\dots,x_n)\) and \(\bar\beta\) the matrix of the linear part of \(\beta\).
\label{thm-beta-continuous}
\end{theorem}

\begin{proof}
An affine bijection of a polytope onto itself preserves the normalised relative Lebesgue measure.
The centroid of a \(\beta\)-invariant measure is fixed by \(\beta\), which in LYM coordinates is the
identity stated, and the image of a measure with covariance \(\Sigma\) under \(\beta:x\mapsto\bar\beta x+b\) has
covariance \(\bar\beta\Sigma\bar\beta^{\mathsf T}\).
\end{proof}

\begin{theorem}
\leanmark{\leanfull}Let \(n\ge3\), and let \(\beta\) be the involution of Lemma \ref{lem-beta}. Whether it preserves the
lattice points of \(\Koz(n)\) depends on the arithmetic of \(n\).
\begin{enumerate}
\item If \(n\) is prime, \(\beta\) maps \(\mathcal L_n\) bijectively onto itself; so the lattice-point
centroid satisfies \(\mu_k+\mu_{n-k}=1+\mu_{n-1}\) exactly, and the covariance matrix
\(\Sigma\) of the uniform measure on \(\mathcal L_n\) is \(\beta\)-invariant in the sense of
Theorem \ref{thm-beta-continuous}, \(\bar\beta\Sigma\bar\beta^{\mathsf T}=\Sigma\).
\item If \(n\ge4\) is composite, \(\beta\) does not preserve \(\mathcal L_n\), so the uniform measure on
\(\mathcal L_n\) is not \(\beta\)-invariant.
\end{enumerate}
\label{thm-beta-lattice}
\end{theorem}

\begin{proof}
(1) \(\beta(x)_k=\binom nk+\binom nk\,x_{n-1}/n-x_{n-k}\), and \(n\mid\binom nk\) for \(0<k<n\) when \(n\)
is prime. (2) Let \(\ell\) be a prime divisor of \(n\), so \(2\le\ell\le n/2\le n-2\). Since
\(\binom n\ell=(n/\ell)\binom{n-1}{\ell-1}\) and the last base-\(\ell\) digit of \(n-1\) is \(\ell-1\),
Lucas's theorem gives \(\binom{n-1}{\ell-1}\equiv1\pmod\ell\), so \(v_\ell\binom n\ell=v_\ell(n)-1\) and
\(\binom n\ell/n\notin\mathbb Z\). The \(f\)-vector with \(f_k=\binom nk\) for \(k\le n-2\),
\(f_{n-1}=1\), \(f_n=0\) lies in \(\mathcal L_n\), and \(\beta(x)_\ell=\binom n\ell/n\).
\end{proof}

For instance at \(n=4\), \(\ell=2\): the \(f\)-vector \((1,4,6,1,0)\) has
\(\beta(x)_2=6+6\cdot\tfrac14-6=\tfrac32\notin\mathbb Z\). On the face \(\Koz(n)^\circ\), where \(\beta\) is
Alexander duality, \(\beta\) maps lattice points to lattice points for every \(n\). Part (2) does not by
itself show that the centroid identity fails for composite \(n\), since a non-symmetric set can have a
symmetric mean.

\begin{observation}
For every composite \(n\le18\) the lattice-point centroid of \(\Koz(n)\) violates the identity
\(\mu_k+\mu_{n-k}=1+\mu_{n-1}\) for some \(k\), and for every prime \(n\le19\) it satisfies it, as
Theorem \ref{thm-beta-lattice}(1) requires.
\label{obs-beta-composite}
\end{observation}

\begin{theorem}
\leanmark{\leanfull}(\emph{The Alexander core.}) Let \(n\ge3\), \(S_n^\circ=\lbrace\vek s\in S_n:s_{n-1}=0\rbrace\) and
\(S_n^{\mathrm{top}}=S_n\setminus S_n^\circ\).
\begin{enumerate}
\item Alexander duality restricts to an involution of \(S_n^\circ\); hence the centroid of \(S_n^\circ\)
satisfies \(\mu_k+\mu_{n-k}=1\) for \(1\le k\le n-1\) and \(\mu_{n-1}=\mu_n=0\).
\item \(\lvert S_n^{\mathrm{top}}\rvert=\sum_{m=0}^{n-1}\lvert S_m\rvert\), with \(\lvert S_0\rvert=1\).
\item With \(\varepsilon_n=\lvert S_n^{\mathrm{top}}\rvert/\lvert S_n\rvert\), the centroid of \(S_n\)
satisfies \(0\le\mu_k+\mu_{n-k}-1\le\varepsilon_n\) for \(1\le k\le n-1\) and
\(0\le\mu_n\le\mu_{n-1}\le\varepsilon_n\), and its mean vertex index
\(\sum_jj\lambda_j=\sum_{k=1}^n\mu_k\) (Proposition \ref{prop-centroid-simplex}) lies in
\([(n-1)/2,\ (n-1)/2+(n+1)\varepsilon_n/2]\).
\end{enumerate}
\label{thm-alexander-core}
\end{theorem}

\begin{proof}
(1) Let \(\vek s\in S_n^\circ\); then \(s_n=0\) too. By \autocite[Lemma 14]{SnellmanArticle2},
\(D_n(\vek s)\) is an \(f\)-vector of an order ideal, with \(D_n(\vek s)_0=1-s_n=1\),
\(D_n(\vek s)_1=n-s_{n-1}=n\) and \(D_n(\vek s)_{n-1}=n-s_1=0\), so \(D_n(\vek s)\in S_n^\circ\); and \(D_n\)
is an involution, \(\mu_k\mapsto1-\mu_{n-k}\). (2) If a proper \(\Delta\) has a face of size \(n-1\),
relabel so that \([n-1]\in\Delta\); then \(\Delta=2^{[n-1]}\cup\lbrace\sigma\cup\lbrace n\rbrace:\sigma\in\Gamma\rbrace\)
with \(\Gamma=\operatorname{lk}_\Delta(n)\) a non-void complex on a subset of \([n-1]\), and
\(f_k(\Delta)=\binom{n-1}k+f_{k-1}(\Gamma)\). So \(f(\Delta)\mapsto f(\Gamma)\) is a bijection from
\(S_n^{\mathrm{top}}\) onto the \(f\)-vectors of non-void complexes on vertex sets contained in
\([n-1]\), which are partitioned by their number \(m\) of vertices into copies of \(S_m\). (3) The
centroid of \(S_n\) is \((1-\varepsilon_n)\,c(S_n^\circ)+\varepsilon_n\,c(S_n^{\mathrm{top}})\). On
\(S_n^{\mathrm{top}}\), \(\mu_k+\mu_{n-k}-1=(f_{k-1}(\Gamma)+f_{n-k-1}(\Gamma))/\binom nk\in[0,1]\);
\(\mu_{n-1}\) vanishes on the core and is at most \(1\); and
\(\sum_k\mu_k=\tfrac12\sum_{k=1}^{n-1}(\mu_k+\mu_{n-k})+\mu_n\).
\end{proof}

Since \(\lvert S_n\rvert\) grows much faster than \(\sum_{m<n}\lvert S_m\rvert\),
\(\varepsilon_n\to0\): its values for \(n=3,\dots,8\) are \(0.8\), \(0.5625\), \(0.357\), \(0.208\),
\(0.112\) and \(0.0573\), and it is \(0.0026\) at \(n=12\). So, asymptotically, the centroid of the
\(f\)-vectors lies in the codimension-two face \(\Koz(n)^\circ\) and is Alexander self-dual up to
\(O(\varepsilon_n)\). Its mean vertex index \(\sum_jj\lambda_j\) is the average of \(j\) when the
vertex \(\tilde F_j\), the \(f\)-vector of the \((j-1)\)-skeleton, is drawn with probability
\(\lambda_j\); by Theorem \ref{thm-alexander-core}(3) it lies between \((n-1)/2\) and
\((n-1)/2+(n+1)\varepsilon_n/2\), against
\((n+1)/2\) for the centroid of the simplex (Proposition \ref{prop-centroid-simplex}). So, measured in
vertices, the \(f\)-vectors sit on average \textbf{one vertex lower} than the simplex: nearer to the skeleta
of one dimension less.

\begin{remark}
(\emph{One map, three behaviours.}) In LYM coordinates the identity
\(\mu_k+\mu_{n-k}=1+\mu_{n-1}\) says that the profile \(k\mapsto\mu_k\) is point-symmetric about
\((n/2,(1+\mu_{n-1})/2)\). It holds exactly for the polytope, where \(\mu_{n-1}=2/n\) and the centre
is \((n/2,\tfrac12+\tfrac1n)\); exactly for the lattice points when \(n\) is prime; and up to
\(\varepsilon_n\) for the \(f\)-vectors, where \(\mu_{n-1}\to0\). The centre is at \(k=n/2\) in every
case; what differs is its height. For lattice points \(\mu_{n-1}\) stays large (\(0.126\) at \(n=7\)),
because the facet \(\mu_n=0\) carries almost every lattice point while the face \(\mu_{n-1}=0\) does
not.
\label{rem-one-map}
\end{remark}

\begin{mexample}
At \(n=3\), where \(\Koz(3)\) is a triangle, \(S_3=\mathcal L_3\) in \((x_2,x_3)\) is \(\lbrace(0,0),(1,0),(2,0),(3,0),(3,1)\rbrace\). Its
centroid is \((9/5,1/5)\), that is \(\mu=(3/5,1/5)\), with barycentric coordinates
\((2/5,2/5,1/5)\): by hand, in \((x_2,x_3)\), \((1,0)=\tfrac23\tilde F_1+\tfrac13\tilde F_2\) and
\((2,0)=\tfrac13\tilde F_1+\tfrac23\tilde F_2\), so, listing the points in order,
\(\lambda_1=(1+\tfrac23+\tfrac13+0+0)/5=\tfrac25\). Kupreyeva's row \(p=3\) of
\autocite[Table 4.2]{Kupreyeva2019}, \((1.80,0.20)\), is this centroid \((9/5,1/5)\) of \(S_3\), since the
\(f\)-vectors in her \(T^3_3\) are \(S_3\). The centroid of the triangle \(\Koz(3)\) is \((2,1/3)\), \(\mu=(2/3,1/3)\),
barycentric \((\tfrac13,\tfrac13,\tfrac13)\). Her row \(p=4\) of Table 5.1 is the centroid
\((x_2,x_3,x_4)=(17/4,19/16,1/16)\) of \(S_4\).
\label{ex-centroid-n3}
\end{mexample}

\begin{center}\small
\begin{tabular}{@{}rlrrrrrr@{}}
\toprule
\(n\) & measure & \(\mu_2\) & \(\mu_3\) & \(\mu_4\) & \(\mu_5\) & \(\mu_6\) & \(\mu_7\)\\
\midrule
4 & \(S_n\) & 0.7083 & 0.2969 & 0.0625 & & &\\
4 & lattice & 0.6833 & 0.3375 & 0.0500 & & &\\
4 & \(\Koz(n)\) & 0.7500 & 0.5000 & 0.2500 & & &\\
5 & \(S_n\) & 0.8014 & 0.4300 & 0.1286 & 0.0143 & &\\
5 & lattice & 0.7407 & 0.4815 & 0.2222 & 0.0062 & &\\
5 & \(\Koz(n)\) & 0.8000 & 0.6000 & 0.4000 & 0.2000 & &\\
6 & \(S_n\) & 0.8709 & 0.5699 & 0.2392 & 0.0514 & 0.0022 &\\
6 & lattice & 0.8024 & 0.5814 & 0.3600 & 0.1508 & 0.0005 &\\
6 & \(\Koz(n)\) & 0.8333 & 0.6667 & 0.5000 & 0.3333 & 0.1667 &\\
7 & \(S_n\) & 0.9179 & 0.6904 & 0.3803 & 0.1292 & 0.0202 & 0.0002\\
7 & lattice & 0.8345 & 0.6508 & 0.4750 & 0.2913 & 0.1258 & 0.0000\\
7 & \(\Koz(n)\) & 0.8571 & 0.7143 & 0.5714 & 0.4286 & 0.2857 & 0.1429\\
\bottomrule
\end{tabular}
\captionof{table}{LYM coordinates of the three centroids at rank one. Exact values, for example
\(c(S_5)=(1,5,\tfrac{561}{70},\tfrac{43}{10},\tfrac9{14},\tfrac1{70})\), are in the
computation log (Section~\ref{sec-methods}).}
\label{tab:centroids-rank1}
\end{center}

\begin{figure}[htbp]
\centering
\includegraphics[width=\linewidth]{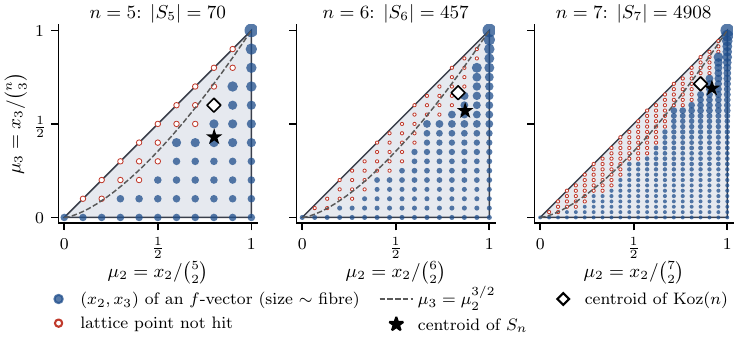}
\caption{The $f$-vectors of proper complexes on $[n]$, $n=5,6,7$, projected to the normalised
coordinates $(\mu_2,\mu_3)=(x_2/\binom n2,x_3/\binom n3)$, in which $\Koz(n)$ projects onto the
triangle $1\ge\mu_2\ge\mu_3\ge0$. The area of a filled dot grows with the number of $f$-vectors
over it; hollow dots are lattice points of the triangle hit by no $f$-vector ($21$, $51$, $119$
of $66$, $171$, $400$). The Kruskal--Katona boundary approaches $\mu_3=\mu_2^{3/2}$ (dashed);
the centroid of $S_n$ ($\star$) sits well inside it, below the centroid of the simplex
$\Koz(n)$ ($\diamond$).}
\label{fig:mu}
\end{figure}

The lattice points that are not \(f\)-vectors are those whose profile stays high in the upper
coordinates: a plateau \(\mu_k\approx\mu_{k+1}<1\) is exactly what the Kruskal--Katona bound, roughly
\(\mu_{k+1}\lesssim\mu_k^{(k+1)/k}\), forbids. The \(f\)-vectors avoid the top vertices
\(\tilde F_{n-1}\), \(\tilde F_n\), put less weight than the lattice points on \(\tilde F_1\), and pile up
on \(\tilde F_j\) with \(j\) near \(n/2\) (Figure \ref{fig:mu}).

\subsection{Rank two and higher: five centroids}
\label{sec:org3b0ad2c}

At rank \(r\) there are five natural centroids, of three uniform measures and two push-forwards:
the centroid of \(\FF(n;\vek d)\); the
\textbf{independent sum} \(\sum_i\shift^{d_i}c(S_n)\), which is the centroid of the push-forward of the
uniform measure on \(S_n^r\), that is, the expected value of \(\sum_i\shift^{d_i}\vek s_i\) for
independent uniform \(\vek s_i\in S_n\); the lattice-point
centroid; the centroid of \(\Koz(n;\vek d)\); and \(\sum_i\shift^{d_i}c(\Koz(n))\), the centroid of a
sum of independent uniform points of the summands.

\begin{proposition}
\begin{enumerate}
\item \leanmark{\leanfull}The two independent sums are the centroids of the push-forward measures, by linearity of
expectation, and the covariances of those measures are the sums of the shifted rank-one
covariances.
\item In the product regime the centroid of \(\FF\) equals the independent sum of the \(c(S_n)\), the
centroid of \(\Koz(n;\vek d)\) equals the independent sum of the \(c(\Koz(n))\), the same two
identities hold for the covariance matrices, and the lattice-point centroid is the sum of the shifted lattice-point centroids
of \(\Koz(n)\).
\item For equal shifts \(\vek d=0^r\), \(\Koz(n;0^r)=r\Koz(n)\) has centroid \(r\,c(\Koz(n))\), but
covariance \(r^2\) times that of \(\Koz(n)\), against \(r\) times for the independent sum.
\end{enumerate}
\label{prop-pushforward}
\end{proposition}

\begin{proof}
(1) is linearity and independence. (2) The sumset map \(S_n^r\to\FF\) is onto, and
\(\lvert\FF\rvert=\lvert S_n\rvert^r\) \autocite[Proposition 42]{SnellmanArticle2}, so it is a bijection, so the uniform measure on \(\FF\) is the push-forward of the
uniform measure on \(S_n^r\); the polytope is a lattice product (Proposition \ref{prop-product-regime})
and Fubini applies to it and to its lattice points. (3) Dilation.
\end{proof}

\begin{observation}
Outside the product regime the identities of Proposition \ref{prop-pushforward}(2) fail in each of the
\(16\) computed cases, those of Table \ref{tab:rank-r-variance} without an asterisk, with one exception:
at \(\vek d=0^r\) the centroid of the polytope is still the independent sum, by (3). For instance
\(c(\FF(n;0^r))\ne r\,c(S_n)\) already at \(n=3\), \(r=2\) (Example \ref{ex-centroid-rank2}).
\label{obs-pushforward-fails}
\end{observation}

\begin{theorem}
\leanmark{\leanfull}For every \(r\ge1\), \(\FF(3;0^r)=r\Koz(3)\cap\mathbb Z^4\). In \((x_2,x_3)\) these are the lattice
points of the triangle \(0\le3x_3\le x_2\le3r\); there are \((r+1)(3r+2)/2\) of them, and their
centroid is
\[
  (x_0,x_1,x_2,x_3)=\Bigl(r,\ 3r,\ \frac{3r(2r+1)}{3r+2},\ \frac{r^2}{3r+2}\Bigr),
  \quad \nu_2=\frac{2r+1}{3r+2}\to\frac23,\quad \nu_3=\frac r{3r+2}\to\frac13 .
\]
So the uniform measure on \(\FF(3;0^r)\) has the same limiting centroid as the polytope, while the
independent sum stays at \(\nu=(3/5,1/5)\) for every \(r\).
\label{thm-n3-equal-shifts}
\end{theorem}

\begin{proof}
The triangle \(\Koz(3)\) is the union of the three unimodular triangles
\(\conv\lbrace(i,0),(i+1,0),(3,1)\rbrace\), \(i=0,1,2\), so it is IDP, and \(r\Koz(3)\cap\mathbb Z^4\) is
\(r\)-fold sums of \(\mathcal L_3=S_3\), which is \(\FF(3;0^r)\); the reverse inclusion is
\(\FF\subseteq\conv\). (It is also Corollaries \ref{cor-additivity} and \ref{cor-small-n}.) Write
\(x_2=3q+s\), \(0\le q<r\), \(s\in\lbrace0,1,2\rbrace\), or \(x_2=3r\); there are
\(\lfloor x_2/3\rfloor+1\) choices of \(x_3\). So the number of points is
\(\sum_{q<r}3(q+1)+(r+1)=(r+1)(3r+2)/2\), \(\sum x_3=\sum_{q<r}3q(q+1)/2+r(r+1)/2=r^2(r+1)/2\), and
\(\sum x_2=\sum_{q<r}(q+1)(9q+3)+3r(r+1)=3r(r+1)(2r+1)/2\). Divide.
\end{proof}

\begin{mexample}
\(\FF(3;0,0)\) consists of the twelve lattice points of the triangle \(2\Koz(3)\), with multiplicities (the
number of ordered pairs of elements of \(S_3\) summing to the point) as exponents. The columns
headed \(\nu\) below are the coordinates \eqref{eq:nu}, \(x_k/\dim F_k\), not barycentric
coordinates, so they need not sum to \(1\):
\[
  (0,0)^1(1,0)^2(2,0)^3(3,0)^4(4,0)^3(5,0)^2(6,0)^1\quad(3,1)^2(4,1)^2(5,1)^2(6,1)^2\quad(6,2)^1 .
\]
In
\((x_2,x_3)\):
\[\small
\begin{array}{lcc}
 & (x_2,x_3) & (\nu_2,\nu_3)\\[2pt]
\FF(3;0,0)=\text{lattice points} & (15/4,\,1/2) & (0.625,\,0.250)\\
\text{independent sum } 2c(S_3) & (18/5,\,2/5) & (0.600,\,0.200)\\
\Koz(3;0,0)=\text{independent sum of }c(\Koz(3)) & (4,\,2/3) & (0.667,\,0.333)
\end{array}
\]
\(\FF(3;0,1)\) has \(22\) points (three of the \(25\) pairs collide) and equals the lattice-point
set of the \(3\)-polytope \(\Koz(3;0,1)\). In \((x_2,x_3,x_4)\):
\[\small
\begin{array}{lcc}
 & (x_2,x_3,x_4) & (\nu_2,\nu_3,\nu_4)\\[2pt]
\FF(3;0,1)=\text{lattice points} & (51/11,\,2,\,5/22) & (0.773,\,0.500,\,0.227)\\
\text{independent sum of }c(S_3) & (24/5,\,2,\,1/5) & (0.800,\,0.500,\,0.200)\\
\Koz(3;0,1) & (37/8,\,7/3,\,3/8) & (0.771,\,0.583,\,0.375)\\
\text{independent sum of }c(\Koz(3)) & (5,\,7/3,\,1/3) & (0.833,\,0.583,\,0.333)
\end{array}
\]
The continuous centroid of \(\Koz(3;0,1)\) can be found by hand from a triangulation into five
tetrahedra of total volume \(6\), and was confirmed by an independent integration. In the product regime,
\(\FF(3;0,2)\), with \(\Koz(3;0,2)\) of dimension \(4\), all five centroids pair off as Proposition \ref{prop-pushforward} says: \(\FF\), the
independent sum and the lattice points at \((1,3,\tfrac{14}5,\tfrac{16}5,\tfrac95,\tfrac15)\), the polytope
and its independent sum at \((1,3,3,\tfrac{10}3,2,\tfrac13)\).
\label{ex-centroid-rank2}
\end{mexample}

\begin{theorem}
\leanmark{\leanpartial}Fix \(n\ge2\) and, for an integer \(t\ge1\), let
\begin{equation}\label{eq:Lambda}
  \Lambda_t=t\Koz(n)\cap\mathbb Z^{n+1} .
\end{equation}
On \(r\Koz(n)\) the
coordinates \(x_0=r\) and \(x_1=nr\) are constant; for the others use \(\nu_k=x_k/(r\binom nk)\),
\(2\le k\le n\), the LYM coordinates divided by \(r\). The target values are those of
\(c(\Koz(n))\) at \(r=1\), \(\nu_k=\mu_k=(n-k+1)/n\) (Proposition \ref{prop-centroid-simplex}).
\begin{enumerate}
\item For every \(r\ge1\) and \(2\le k\le n\),
\[
     \bigl\lvert\nu_k\bigl(c(\FF(n;0^r))\bigr)-\nu_k\bigl(c(\Lambda_r)\bigr)\bigr\rvert\le1-\operatorname{fill}(n;0^r),
   \]
and both \(1-\operatorname{fill}(n;0^r)\) and \(\nu_k(c(\Lambda_r))-(n-k+1)/n\) are \(O(1/r)\). Hence
\(c(\FF(n;0^r))/r\to c(\Koz(n))\) as \(r\to\infty\), with error \(O(1/r)\) in every coordinate.
\item The covariance matrix of the uniform measure on \(\FF(n;0^r)\), divided by \(r^2\), tends to that
of \(\Koz(n)\), with error \(O(1/r)\) in every entry.
\item For every \(r\ge1\) and \(2\le k\le n\), \(\bigl\lvert\nu_k\bigl(c(\FF(n;0^r))\bigr)-(n-k+1)/n\bigr\rvert\) is
less than \(13/(16r)\) when \(n=4\), and less than \(199/(100r)\) when \(n=5\).
\end{enumerate}
\label{thm-centroid-dilation}
\end{theorem}

\begin{proof}
\emph{This proof was found by the model; see the disclosure section.} Write \(F_r=\FF(n;0^r)\),
\(f_r=\operatorname{fill}(n;0^r)\), \(L=L_{\Koz(n)}\), \(a(r)=\lvert F_r\rvert\) as in \eqref{eq:ar}, and
\(\bar x_k=x_k\bigl(c(\Koz(n))\bigr)=\binom nk(n-k+1)/n\) for \(2\le k\le n\).

\emph{Inclusion and range.} \(F_r=S_n+\dots+S_n\subseteq\Koz(n)+\dots+\Koz(n)=r\Koz(n)\), and \(F_r\)
consists of integer points, so \(F_r\subseteq\Lambda_r\) (this is also Corollary \ref{cor-additivity}).
By Proposition \ref{prop-rank-one-H} every point of \(r\Koz(n)\) has \(x_0=r\), \(x_1=nr\) and
\(r\ge x_2/\binom n2\ge\dots\ge x_n/\binom nn\ge0\), so \(0\le\nu_k\le1\) on \(\Lambda_r\).

\emph{Mixing.} Let \(A\subseteq B\) be finite with \(A\ne\emptyset\), let \(g\) take values in \([0,1]\) on
\(B\), write \(\mathbb E_Xg\) for the average over \(X\), and \(f=\lvert A\rvert/\lvert B\rvert\). If
\(C=B\setminus A\) is non-empty, then \(\mathbb E_Bg=f\,\mathbb E_Ag+(1-f)\,\mathbb E_Cg\), so
\(\mathbb E_Ag-\mathbb E_Bg=(1-f)(\mathbb E_Ag-\mathbb E_Cg)\) and
\(\lvert\mathbb E_Ag-\mathbb E_Bg\rvert\le1-f\); if \(C\) is empty, both sides vanish. For \(A=F_r\),
\(B=\Lambda_r\) and the linear function \(g=\nu_k\) this is the inequality in (1).

\emph{The fill rate.} In the proof of Corollary \ref{cor-zero-shift}, \(a(r)\) and \(L(r)\) are, for \(r\ge1\),
of the form \(\sum_{j=0}^{n-1}\binom{r-1}jw_j\) with constants \(w_j\) (different for the two), and the
top constants agree: both are \(\prod_{k=2}^n\binom nk\), the single non-zero entry of
\(\mathcal N^{n-1}\), times the first entry \(1\) of the row vector. So \(L-a\) is a polynomial in
\(r\) of degree at most \(n-2\), while \(L\) has degree \(n-1\) and leading coefficient
\(\prod_{k\ge2}\binom nk/(n-1)!=\vol\Koz(n)>0\). Hence \(1-f_r=(L(r)-a(r))/L(r)=O(1/r)\), and
\(r(1-f_r)\) tends to \(\theta_n=[r^{n-2}](L-a)/\vol\Koz(n)\).

\emph{The lattice-point centroid.} Forgetting \(x_0\) and \(x_1\) maps \(t\Koz(n)\) onto \(tK\), where
\(K=\lbrace1\ge y_2/\binom n2\ge\dots\ge y_n/\binom nn\ge0\rbrace\subseteq\mathbb R^{n-1}\) is a
full-dimensional lattice simplex; as in Section \ref{sec-volume} the map is a lattice-preserving
isometry on the affine hull, so it maps \(\Lambda_t\) onto \(tK\cap\mathbb Z^{n-1}\), \(c(\Koz(n))\)
onto \(c(K)\), and \(\vol\Koz(n)=\vol K\). For \(2\le k\le n\) let
\begin{equation}\label{eq:Qk}
  Q_k=\bigl\lbrace(y,z)\in\mathbb R^{n-1}\times\mathbb R:\ y\in K,\ 0\le z\le y_k\bigr\rbrace .
\end{equation}
It is the convex hull of the points \((v,0)\) and \((v,v_k)\), \(v\) a vertex of \(K\), since a point
\((y,sy_k)\) with \(y=\sum_v\lambda_vv\) and \(0\le s\le1\) equals
\(\sum_v\lambda_v\bigl((1-s)(v,0)+s(v,v_k)\bigr)\). So \(Q_k\) is a lattice polytope; it has dimension
\(n\), because \(y_k>0\) on the interior of \(K\), and volume \(\int_Ky_k\,dy=\vol K\cdot\bar x_k\) by
Fubini. As \(tQ_k=\lbrace(y,z):y\in tK,\ 0\le z\le y_k\rbrace\) and \(y_k\) is a non-negative integer
at a lattice point,
\begin{equation}\label{eq:Sigmak}
\begin{gathered}
  \#\bigl(tQ_k\cap\mathbb Z^n\bigr)=\sum_{x\in\Lambda_t}(x_k+1)=\Sigma_k(t)+L(t),\qquad
  \Sigma_k(t)=\sum_{x\in\Lambda_t}x_k ,\\
  B_k(t)=\Sigma_k(t)-\bar x_k\,t\,L(t) .
\end{gathered}
\end{equation}
By Ehrhart's theorem \autocite{BeckRobins2015}, applied to \(Q_k\) and to \(K\), \(\Sigma_k\) agrees for
\(t\ge1\) with a polynomial of degree \(n\) with leading coefficient \(\vol\Koz(n)\,\bar x_k\), so
\(B_k\) has degree at most \(n-1\). Therefore
\begin{equation}\label{eq:betak}
  \nu_k\bigl(c(\Lambda_t)\bigr)-\frac{n-k+1}n=\frac{B_k(t)}{\binom nk\,t\,L(t)}
  =\frac{\beta_k}t+O\Bigl(\frac1{t^2}\Bigr),\qquad \beta_k=\frac{[t^{n-1}]B_k}{\binom nk\vol\Koz(n)} ,
\end{equation}
which completes (1); moreover
\(\limsup_{r\to\infty}r\,\lvert\nu_k(c(F_r))-(n-k+1)/n\rvert\le\theta_n+\lvert\beta_k\rvert\). The
constants \(\beta_k\) computed below at \(n=4,5\) are the first-order corrections of
Proposition \ref{prop-dilates}.

\emph{Second moments.} For \(2\le k,l\le n\) the mixing step with \(g=\nu_k\nu_l\in[0,1]\) gives
\(\lvert\mathbb E_{F_r}[\nu_k\nu_l]-\mathbb E_{\Lambda_r}[\nu_k\nu_l]\rvert\le1-f_r\). The polytope
\(Q_{kl}=\lbrace(y,z,w):y\in K,\ 0\le z\le y_k,\ 0\le w\le y_l\rbrace\) is, by the same
interpolation in \(z\) and in \(w\) separately, the convex hull of the points
\((v,\varepsilon v_k,\eta v_l)\) with \(\varepsilon,\eta\in\lbrace0,1\rbrace\); it has dimension \(n+1\)
(also for \(k=l\)) and volume \(\int_Ky_ky_l\,dy\), and
\(\#(tQ_{kl}\cap\mathbb Z^{n+1})=\sum_{x\in\Lambda_t}(x_k+1)(x_l+1)\). So
\(\sum_{x\in\Lambda_t}x_kx_l\) agrees for \(t\ge1\) with a polynomial of degree \(n+1\) with leading
coefficient \(\int_Ky_ky_l\,dy\), and \(\mathbb E_{\Lambda_t}[\nu_k\nu_l]\) is within \(O(1/t)\) of the
corresponding moment \(\mathbb E[\mu_k\mu_l]\) of \(\Koz(n)\). With (1), and since all these numbers
lie in \([0,1]\), where \(\lvert ab-a'b'\rvert\le\lvert a-a'\rvert+\lvert b-b'\rvert\), the covariance
\(\mathbb E_{F_r}[\nu_k\nu_l]-\mathbb E_{F_r}[\nu_k]\,\mathbb E_{F_r}[\nu_l]\) is within \(O(1/r)\) of the
covariance of \(\mu_k\) and \(\mu_l\) on \(\Koz(n)\). This is (2) in the coordinates \(\nu\); in the
coordinates \(x\) both covariances are multiplied by \(\binom nk\binom nl\), and \(x_0\), \(x_1\) have
variance \(0\) on both sides.

\emph{The constants in (3).} By Corollary \ref{cor-deficiency}, \(1-f_r=2r/((r+1)(4r+1))<1/(2r)\) at \(n=4\),
and \(1-f_r=r(205r+71)/(125r^3+250r^2+105r+6)<41/(25r)\) at \(n=5\); the two bounds are \(\theta_4/r\)
and \(\theta_5/r\). By the degree bounds just proved, \(L\) and \(\Sigma_k\) are determined by their values
at \(t=1,\dots,n+1\), and computed from these (and checked at \(t=n+2\)) they give, at \(n=4\),
with \(L(t)=(t+1)^2(4t+1)\),
\[
  B_2=-\tfrac12t(t+1)(5t+3),\qquad B_3=-\tfrac12t(t+1)(8t+5),\qquad B_4=-\tfrac14t(t+1)(5t+3),
\]
so \(\beta=(-\tfrac5{48},-\tfrac14,-\tfrac5{16})\); and at \(n=5\), with
\(6L(t)=(t+1)(125t^3+250t^2+105t+6)\) and \(v(t)=125t^2+165t-2\),
\[
  B_2=-\tfrac16t(t+1)v,\qquad B_3=-\tfrac13t(t+1)v,\qquad B_4=-\tfrac14t(t+1)v,
\]
\[
  B_5=-\tfrac1{120}t(t+1)(875t^2+1095t-86),
\]
so \(\beta=(-\tfrac1{10},-\tfrac15,-\tfrac3{10},-\tfrac7{20})\). Every \(B_k\) is negative for \(t\ge1\),
and \(\lvert\beta_k\rvert\binom nkL+B_k\) equals, at \(n=4\),
\(\tfrac18(t+1)(13t+5)\), \(\tfrac12(t+1)(5t+2)\), \(\tfrac1{16}(t+1)(13t+5)\), and at \(n=5\),
\(\tfrac16(t+1)(17t+1)(5t+6)\), \(\tfrac13(\dots)\), \(\tfrac14(\dots)\) (the same product) and
\(\tfrac1{120}(t+1)(131t+7)(5t+6)\), all positive. So
\(\lvert\nu_k(c(\Lambda_t))-(n-k+1)/n\rvert=-B_k/(\binom nk\,t\,L)\le\lvert\beta_k\rvert/t\) for every
\(t\ge1\), with \(\lvert\beta_k\rvert\le5/16\) at \(n=4\) and \(\le7/20\) at \(n=5\). Adding:
\(1/2+5/16=13/16\) and \(41/25+7/20=199/100\).
\end{proof}

Theorem \ref{thm-n3-equal-shifts} is the case \(n=3\), exactly; for \(n\le3\) every lattice point of the
dilate is an ff-vector (Corollary \ref{cor-small-n}), so the content is at \(n\ge4\). There the missing
lattice points are of lower order (Corollary \ref{cor-zero-shift}), and that alone forces the centroid;
which lattice points are missing does not matter. The bounds are far from sharp. At \(n=4\), \(r=10\),
the centroid of the ff-vectors, \(\nu=(0.744,0.468,0.217)\), is within \(0.008\) of the lattice-point
centroid and within \(0.033\) of \((3/4,1/2,1/4)\), against \(1-\operatorname{fill}=0.044\) and
\(13/160\approx0.081\); at \(n=5\), \(r=9\), it is within \(0.018\) of the lattice-point centroid,
against \(1-\operatorname{fill}=0.154\). Part (2) says that at equal shifts the uniform measure on
the Hilbert functions spreads like the polytope \(r\Koz(n)\), quadratically in \(r\), and not like the
independent sum, whose covariance grows linearly (Proposition \ref{prop-pushforward}(3)). At a constant
positive gap the fill rate tends to \(0\) (Observation \ref{obs-fill}), and the argument says nothing.

In the barycentric reading of Theorem \ref{thm-alexander-core}(3), part (1) says that the offset found
there does not survive dilation. At rank one the \(f\)-vectors sit about one vertex lower than the
simplex as \(n\) grows, whereas at fixed \(n\), as \(r\to\infty\), the mean vertex index of
\(c(\FF(n;0^r))/r\) tends to \((n+1)/2\), that of \(c(\Koz(n))\), with error \(O(1/r)\). At \(n=4\) it is
\(2.07\) for \(r=1\) and \(2.43\) for \(r=10\), against \(2.5\).

\begin{figure}[htbp]
\centering
\includegraphics[width=\linewidth]{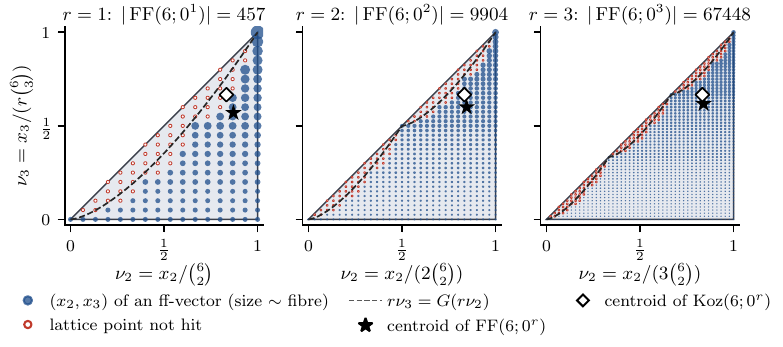}
\caption{The ff-vectors of $\FF(6;0^r)$, $r=1,2,3$, projected to the coordinates
$(\nu_2,\nu_3)=\bigl(x_2/(r\binom62),x_3/(r\binom63)\bigr)$ of \eqref{eq:nu}, in which
$\Koz(6;0^r)=r\Koz(6)$ projects onto the triangle $1\ge\nu_2\ge\nu_3\ge0$. The panel $r=1$ is the
middle panel of Figure~\ref{fig:mu}; the conventions are those of that figure, with dot areas scaled
by $1/r^2$ along with the grid. The dashed curve $r\nu_3=G(r\nu_2)$,
$G(v)=\lfloor v\rfloor+(v-\lfloor v\rfloor)^{3/2}$, is the upper edge of the Minkowski sum of $r$
copies of the region below $\mu_3=\mu_2^{3/2}$, scaled by $1/r$: the Kruskal--Katona boundary
becomes $r$ arches. The lattice points hit by no ff-vector, $51$ of $171$, $102$ of $641$ and $153$
of $1411$, shrink to thin lenses along the diagonal (Corollary~\ref{cor-zero-shift}), and the
centroid of the ff-vectors ($\star$) approaches that of the polytope ($\diamond$)
(Theorem~\ref{thm-centroid-dilation}).}
\label{fig:mu-rank-r}
\end{figure}

Figure \ref{fig:mu-rank-r} shows Corollary \ref{cor-zero-shift} and Theorem \ref{thm-centroid-dilation} at
\(n=6\), in the coordinates of Figure \ref{fig:mu}.

\begin{remark}
(\emph{Weak convergence.}) The proof of Theorem \ref{thm-centroid-dilation} gives more. Write, as there,
\(K\subseteq\mathbb R^{n-1}\) for the projection of \(\Koz(n)\) to the coordinates
\((x_2,\dots,x_n)\), and \(\Lambda_t\) as in \eqref{eq:Lambda}. For a monomial
\(x^\alpha=\prod_{k\ge2}x_k^{\alpha_k}\), add to \(K\) one coordinate \(0\le z_i\le y_{k(i)}\) for each of
its \(\lvert\alpha\rvert\) factors, \(k(i)\) being the index of the \(i\)-th factor. Over \(y\in K\) the
fibre is a box, so the resulting polytope is the
convex hull of the points \((v,\text{corner of the box over }v)\), \(v\) a vertex of \(K\), a lattice
polytope whose dilate by \(t\) has \(\sum_{x\in\Lambda_t}\prod_k(x_k+1)^{\alpha_k}\) lattice points. By
induction on \(\lvert\alpha\rvert\), \(\sum_{x\in\Lambda_t}x^\alpha\) agrees for \(t\ge1\) with a polynomial
of degree \(n-1+\lvert\alpha\rvert\) whose leading coefficient is \(\int_Ky^\alpha\,dy\). With the mixing
step of that proof for \(g=\nu^\alpha\in[0,1]\), every moment of the uniform measure on \(\FF(n;0^r)/r\) converges to
the corresponding moment of the uniform probability measure on \(\Koz(n)\). All these measures live on
the compact set \(\Koz(n)\), where polynomials are dense in the continuous functions
(Stone--Weierstrass), so the uniform measure on \(\FF(n;0^r)/r\) converges weakly to the uniform
measure on \(\Koz(n)\). Parts (1) and (2) of the theorem are its first two moments, with a rate.
\label{rem-weak-convergence}
\end{remark}
\subsection{Lattice points of dilates}
\label{sec:org84c90fa}

For a lattice polytope \(P\) of dimension \(D\), the Ehrhart polynomial \(L_P\) of
Section \ref{sec-ehrhart} has leading coefficient the relative volume of \(P\), and next coefficient
half the sum of the relative volumes of its facets \autocite[Theorem 5.6]{BeckRobins2015}. The same
holds for the sum of a linear form \(\ell\) with integer coefficients, nonnegative on \(P\). As in
Remark \ref{rem-weak-convergence}, the lattice polytope \(Q=\{(x,z):x\in P,\ 0\le z\le\ell(x)\}\) has
\(\sum_{x\in tP}(\ell(x)+1)\) lattice points in \(tQ\). Its facets are \(P\); the graph of \(\ell\)
over \(P\), of the same relative volume; and, over each facet \(F\) of \(P\), the region under
\(\ell\), of relative volume \(\int_F\ell\). So
\[
  \sum_{x\in tP\cap\mathbb Z^{N+1}}\ell(x)=L_Q(t)-L_P(t)
  =t^{D+1}\int_P\ell+\frac{t^D}2\sum_F\int_F\ell+O(t^{D-1}),
\]
all integrals with the lattice-normalised measures. For polynomial weights on a simple lattice
polytope, such as \(\Koz(n)\), the corresponding statement follows from the Euler--Maclaurin formula of Brion and Vergne
\autocite[Theorem 2.15]{BrionVergne1997}. For \(\Koz(n)\), with \(\ell=x_k\), this predicts the lattice
centroid of \(t\Koz(n)\) to be \(t\,c(\Koz(n))+b+O(1/t)\) for a constant vector \(b\).

\begin{proposition}
\leanmark{\leannone}(\emph{The lattice-point centroid of \(t\Koz(n)\), \(n\le5\).}) Let \(\Lambda_t\) be as in
\eqref{eq:Lambda} and, as in Theorem \ref{thm-centroid-dilation}, \(\nu_k=x_k/(t\binom nk)\), the LYM coordinates of
\(t\Koz(n)\) divided by \(t\). For every integer \(t\ge1\) the coordinates of the lattice-point centroid
\(c(\Lambda_t)\) are as follows.
\begin{enumerate}
\item For \(n=3\): \(\nu_2=\frac23-\frac1{3(3t+2)}\) and \(\nu_3=\frac13-\frac2{3(3t+2)}\).
\item For \(n=4\): \(\nu_2=\frac34-\frac{5t+3}{12(t+1)(4t+1)}\), \(\nu_3=\frac12-\frac{8t+5}{8(t+1)(4t+1)}\) and
\(\nu_4=\frac14-\frac{5t+3}{4(t+1)(4t+1)}\).
\item For \(n=5\), with \(u(t)=125t^3+250t^2+105t+6\) and \(v(t)=125t^2+165t-2\):
\(\nu_2=\frac45-\frac{v}{10u}\), \(\nu_3=\frac35-\frac v{5u}\), \(\nu_4=\frac25-\frac{3v}{10u}\) and
\(\nu_5=\frac15-\frac{875t^2+1095t-86}{20u}\).
\end{enumerate}
Every correction is negative, and to first order the centroid is
\[\small
\begin{array}{ll}
n=3: & \bigl(\tfrac23-\tfrac1{9t},\ \tfrac13-\tfrac2{9t}\bigr)\\[2pt]
n=4: & \bigl(\tfrac34-\tfrac5{48t},\ \tfrac12-\tfrac1{4t},\ \tfrac14-\tfrac5{16t}\bigr)\\[2pt]
n=5: & \bigl(\tfrac45-\tfrac1{10t},\ \tfrac35-\tfrac1{5t},\ \tfrac25-\tfrac3{10t},\ \tfrac15-\tfrac7{20t}\bigr)
\end{array}
\]
up to \(O(1/t^2)\).
\label{prop-dilates}
\end{proposition}

\begin{proof}
As in the proof of Theorem \ref{thm-centroid-dilation}, the Ehrhart polynomial
\(L_{\Koz(n)}(t)=\lvert\Lambda_t\rvert\) and the sum \(\Sigma_k(t)\) of \eqref{eq:Sigmak} agree for
\(t\ge1\) with polynomials of degree \(n-1\) and \(n\), so they are determined by their values at
\(t=1,\dots,n+1\). Computing these values by listing \(\Lambda_t\) (and checking at \(t=n+2\) and
\(t=n+3\)) gives \(L_{\Koz(n)}\) and the \(B_k\) of \eqref{eq:Sigmak}: for \(n=3\),
\(L_{\Koz(3)}(t)=\tfrac12(t+1)(3t+2)\), \(B_2=-\tfrac12t(t+1)\) and \(B_3=-\tfrac13t(t+1)\); for
\(n=4,5\), the polynomials displayed at the end of the proof of Theorem \ref{thm-centroid-dilation}.
Then, by \eqref{eq:betak},
\(\nu_k\bigl(c(\Lambda_t)\bigr)-(n-k+1)/n=B_k(t)/\bigl(\binom nk\,t\,L_{\Koz(n)}(t)\bigr)\), and
cancelling \(t(t+1)\) gives the formulas. Their numerators \(5t+3\), \(8t+5\), \(v(t)\) and
\(875t^2+1095t-86\) are
positive for \(t\ge1\), so every correction is negative; expanding in \(1/t\) gives the table.
\end{proof}

The \(1/t\) corrections agree with the boundary term computed independently from the facets.
Heuristically they are negative because, for small \(t\), the facet \(\mu_n=0\) carries nearly all
lattice points.
\section{Discrete moments versus moments}
\label{sec-moments}
\begin{proposition}
\begin{enumerate}
\item \leanmark{\leannone}For the uniform measure on a \(D\)-simplex with vertices \(v_0,\dots,v_D\) and \(s=\sum_av_a\),
\(\mathbb E[x_ix_j]=\bigl(\sum_av_{a,i}v_{a,j}+s_is_j\bigr)/\bigl((D+1)(D+2)\bigr)\).
\item For \(\Koz(n)\), in LYM coordinates and for \(2\le k\le l\le n\),
\(\operatorname{Cov}(\mu_k,\mu_l)=(k-1)(n-l+1)/\bigl(n^2(n+1)\bigr)\). More precisely,
\(\mu_2\ge\dots\ge\mu_n\) are the order statistics of \(n-1\) independent uniform variables on
\([0,1]\), \(\mu_k\) being the \((n-k+1)\)-th smallest.
\end{enumerate}
\label{prop-second-moments}
\end{proposition}

\begin{proof}
(1) A uniform point is \(\sum_a\lambda_av_a\) with \(\lambda\) Dirichlet\((1,\dots,1)\) on
\(m=D+1\) parameters, for which \(\mathbb E[\lambda_a\lambda_b]=(1+\delta_{ab})/(m(m+1))\). (2) In LYM
coordinates the vertices are \((1^j,0^{n-j})\), \(m=n\), \(\sum_jv_{j,k}v_{j,l}=n-l+1\) and
\(s_k=n-k+1\); with \(a=n-k+1\) and \(b=n-l+1\),
\(\mathbb E[\mu_k\mu_l]=(b+ab)/(n(n+1))\), and subtracting \(ab/n^2\) gives \(b(n-a)/(n^2(n+1))\). The
uniform measure on the order simplex is the law of the sorted sample.
\end{proof}

\begin{mexample}
At \(n=3\), in \((x_2,x_3)\), the covariance of the five \(f\)-vectors is
\(\left(\begin{smallmatrix}34/25&6/25\\6/25&4/25\end{smallmatrix}\right)\), and that of the triangle
\(\left(\begin{smallmatrix}1/2&1/12\\1/12&1/18\end{smallmatrix}\right)\): the discrete set is
\(2.7\) times more spread in \(x_2\), because its five points sit on the boundary. For
\(\Koz(3;0,1)\), in \((x_2,x_3,x_4)\), the covariance of the \(22\) ff-vectors and that of the polytope
are
\[\small
\begin{pmatrix}160/121&3/11&9/242\\3/11&18/11&3/11\\9/242&3/11&85/484\end{pmatrix}
\qquad\text{and}\qquad
\begin{pmatrix}47/64&7/48&1/64\\7/48&13/18&5/48\\1/64&5/48&13/192\end{pmatrix}.
\]
\label{ex-moments-n3}
\end{mexample}

For small \(n\), then, the discrete measure is more spread than the continuous one; for large \(n\)
the opposite holds, and markedly. Tables \ref{tab:concentration} to \ref{tab:profile-n24} show this
at rank one for \(n\le26\), and Table \ref{tab:rank-r-variance} at rank two and three.

\begin{observation}
(\emph{The \(f\)-vectors concentrate.}) For even \(n\le26\), the standard deviation of the middle
coordinate \(\mu_{n/2}\) under the uniform measure on \(S_n\) is about \(1.78/n\) (Table
\ref{tab:concentration}), whereas both the lattice-point and the continuous measures have standard
deviation of order \(n^{-1/2}\); by Proposition \ref{prop-second-moments} the continuous one is
\(\sqrt{(k-1)(n-k+1)/(n^2(n+1))}\approx1/(2\sqrt n)\). The standard deviations of the middle coordinate
\(\mu_{\lceil n/2\rceil}\) for \(S_n\) and for \(\Koz(n)\), as functions of \(n\), cross between \(n=14\) and
\(n=15\). So, as \(n\) grows, a uniformly random proper \(f\)-vector has an increasingly
deterministic LYM profile, a sigmoid symmetric about \(k=n/2\) (Table \ref{tab:profile-n24}: at
\(n=24\), \(\mu\approx0.915\) at \(k=n/4\) and \(\approx0.085\) at \(3n/4\)). Let \(J\) be the barycentric
index: the index of a vertex drawn with probability \(\lambda_J\) at a random point. For the simplex
\(J\) is uniform on \([n]\), with variance \((n^2-1)/12\). Under the uniform measure on \(S_n\),
\(\operatorname{Var}(J)/n^2\) increases slowly, from \(0.026\) to \(0.032\) over \(12\le n\le26\),
against \(1/12\) for the simplex (Table \ref{tab:barycentric-index}). All verified by exact dynamic
programming, \(n\le26\).
\label{obs-concentration}
\end{observation}

\begin{table}[htbp]\centering\small
\begin{tabular}{@{}rrrrrrr@{}}
\toprule
 & \multicolumn{2}{c}{\(S_n\)} & \multicolumn{2}{c}{lattice points} & \multicolumn{2}{c}{\(\Koz(n)\)}\\
\cmidrule(lr){2-3}\cmidrule(lr){4-5}\cmidrule(l){6-7}
\(n\) & sd & \(n\cdot\mathrm{sd}\) & sd & \(\sqrt n\cdot\mathrm{sd}\) & sd & \(\sqrt n\cdot\mathrm{sd}\)\\
\midrule
4 & 0.3146 & 1.258 & 0.2977 & 0.595 & 0.1936 & 0.387\\
6 & 0.2896 & 1.738 & 0.2325 & 0.570 & 0.1782 & 0.436\\
8 & 0.2344 & 1.875 & 0.1937 & 0.548 & 0.1614 & 0.456\\
10 & 0.1866 & 1.866 & 0.1679 & 0.531 & 0.1477 & 0.467\\
12 & 0.1528 & 1.834 & 0.1511 & 0.523 & 0.1367 & 0.474\\
14 & 0.1292 & 1.809 & 0.1390 & 0.520 & 0.1278 & 0.478\\
16 & 0.1121 & 1.794 & 0.1290 & 0.516 & 0.1203 & 0.481\\
18 & 0.0991 & 1.784 & 0.1211 & 0.514 & 0.1140 & 0.484\\
20 & 0.0889 & 1.778 & 0.1145 & 0.512 & 0.1086 & 0.486\\
22 & 0.0807 & 1.776 & 0.1089 & 0.511 & 0.1038 & 0.487\\
24 & 0.0740 & 1.776 & 0.1040 & 0.510 & 0.0997 & 0.488\\
26 & 0.0684 & 1.779 & 0.0998 & 0.509 & 0.0959 & 0.489\\
\bottomrule
\end{tabular}
\caption{Standard deviation of the middle LYM coordinate \(\mu_{n/2}\) under the uniform
measures on \(S_n\), on the lattice points of \(\Koz(n)\), and on \(\Koz(n)\) (the last by
Proposition~\ref{prop-second-moments}), with the scalings of Observation~\ref{obs-concentration}:
\(n\cdot\mathrm{sd}\) settles for \(S_n\), \(\sqrt n\cdot\mathrm{sd}\) for the other two.}
\label{tab:concentration}
\end{table}

\begin{table}[htbp]\centering\small
\setlength{\tabcolsep}{3.5pt}
\begin{tabular}{@{}lrrrrrrrrrrrr@{}}
\toprule
\(n\) & 4 & 6 & 8 & 10 & 12 & 14 & 16 & 18 & 20 & 22 & 24 & 26\\
\midrule
\(S_n\) & 0.048 & 0.032 & 0.027 & 0.026 & 0.026 & 0.027 & 0.028 & 0.029 & 0.030 & 0.031 & 0.031 & 0.032\\
lattice & 0.050 & 0.050 & 0.057 & 0.061 & 0.064 & 0.067 & 0.069 & 0.070 & 0.072 & 0.073 & 0.074 & 0.074\\
\(\Koz(n)\) & 0.078 & 0.081 & 0.082 & 0.083 & 0.083 & 0.083 & 0.083 & 0.083 & 0.083 & 0.083 & 0.083 & 0.083\\
\bottomrule
\end{tabular}
\caption{\(\operatorname{Var}(J)/n^2\) for the barycentric index \(J\) of
Observation~\ref{obs-concentration}, under the three uniform measures; for \(\Koz(n)\) it is
\((n^2-1)/(12n^2)\to1/12\).}
\label{tab:barycentric-index}
\end{table}

\begin{table}[htbp]\centering\small
\begin{tabular}{@{}rrrrrrr@{}}
\toprule
 & \multicolumn{2}{c}{\(S_{24}\)} & \multicolumn{2}{c}{lattice points} & \multicolumn{2}{c}{\(\Koz(24)\)}\\
\cmidrule(lr){2-3}\cmidrule(lr){4-5}\cmidrule(l){6-7}
\(k\) & mean & sd & mean & sd & mean & sd\\
\midrule
2 & 1.000 & 0.001 & 0.957 & 0.043 & 0.958 & 0.040\\
4 & 0.982 & 0.017 & 0.869 & 0.070 & 0.875 & 0.066\\
6 & 0.915 & 0.045 & 0.780 & 0.086 & 0.792 & 0.081\\
8 & 0.797 & 0.069 & 0.691 & 0.096 & 0.708 & 0.091\\
10 & 0.652 & 0.076 & 0.602 & 0.102 & 0.625 & 0.097\\
12 & 0.500 & 0.074 & 0.513 & 0.104 & 0.542 & 0.100\\
14 & 0.348 & 0.076 & 0.425 & 0.103 & 0.458 & 0.100\\
16 & 0.203 & 0.069 & 0.336 & 0.099 & 0.375 & 0.097\\
18 & 0.085 & 0.045 & 0.247 & 0.091 & 0.292 & 0.091\\
20 & 0.018 & 0.017 & 0.158 & 0.079 & 0.208 & 0.081\\
22 & 0.000 & 0.001 & 0.069 & 0.058 & 0.125 & 0.066\\
24 & 0.000 & 0.000 & 0.000 & 0.000 & 0.042 & 0.040\\
\bottomrule
\end{tabular}
\caption{The LYM profile at \(n=24\): mean and standard deviation of \(\mu_k\), even \(k\),
under the three uniform measures. The means are the LYM coordinates of the three centroids, as in
Table~\ref{tab:centroids-rank1}; for \(\Koz(24)\) they are \((25-k)/24\)
(Proposition~\ref{prop-centroid-simplex}).}
\label{tab:profile-n24}
\end{table}

At rank two and three, in the computed range, the small-\(n\) behaviour persists. Write the \textbf{total
\(\nu\)-variance} of a measure on \(\Koz(n;\vek d)\) for the trace of its covariance matrix in the
coordinates \(\nu_k\) of \eqref{eq:nu}.

\begin{observation}
(\emph{Rank two and three.}) In each of the \(21\) cases of Table \ref{tab:rank-r-variance}, the uniform
measure on \(\FF(n;\vek d)\) has larger total \(\nu\)-variance than the uniform measure on
\(\Koz(n;\vek d)\): the small-\(n\) boundary effect again.
\label{obs-moments-rank-r}
\end{observation}

\begin{table}[htbp]\centering\small
\setlength{\tabcolsep}{4.5pt}
\begin{tabular}{@{}rlrrr@{\qquad}rlrrr@{}}
\toprule
\multicolumn{5}{c}{rank two} & \multicolumn{5}{c}{rank three}\\
\cmidrule(r){1-5}\cmidrule(l){6-10}
\(n\) & \(\vek d\) & \(\FF\) & indep. & \(\Koz\) & \(n\) & \(\vek d\) & \(\FF\) & indep. & \(\Koz\)\\
\midrule
3 & \((0,0)\) & 0.2066 & 0.1556 & 0.1111 & 3 & \((0,0,0)\) & 0.1739 & 0.1037 & 0.1111\\
3 & \((0,1)\) & 0.3146 & 0.2928 & 0.1332 & 3 & \((0,0,1)\) & 0.3282 & 0.2608 & 0.1506\\
3 & \((0,2)^*\) & 0.4061 & 0.4061 & 0.1458 & 3 & \((0,1,1)\) & 0.2101 & 0.1556 & 0.1209\\
4 & \((0,0)\) & 0.1967 & 0.1382 & 0.1250 & 3 & \((0,1,2)\) & 0.3420 & 0.3138 & 0.1524\\
4 & \((0,1)\) & 0.2835 & 0.2273 & 0.1596 & 3 & \((0,1,3)\) & 0.4243 & 0.4197 & 0.1646\\
4 & \((0,2)\) & 0.3555 & 0.3538 & 0.1806 & 3 & \((0,2,4)^*\) & 0.5011 & 0.5011 & 0.1806\\
4 & \((0,3)^*\) & 0.4538 & 0.4538 & 0.1960 & 4 & \((0,0,0)\) & 0.1722 & 0.0921 & 0.1250\\
5 & \((0,0)\) & 0.1732 & 0.1153 & 0.1333 & 4 & \((0,1,2)\) & 0.3205 & 0.2400 & 0.1900\\
5 & \((0,1)\) & 0.2397 & 0.1706 & 0.1781 & 4 & \((0,3,6)^*\) & 0.6312 & 0.6312 & 0.2670\\
5 & \((0,2)\) & 0.2990 & 0.2812 & 0.2079 & & & & & \\
5 & \((0,3)\) & 0.3778 & 0.3801 & 0.2225 & & & & & \\
5 & \((0,4)^*\) & 0.4319 & 0.4319 & 0.2285 & & & & & \\
\bottomrule
\end{tabular}
\caption{Total \(\nu\)-variance of the uniform measure on \(\FF(n;\vek d)\), of the
independent sum (the push-forward of the uniform measure on \(S_n^r\),
Proposition~\ref{prop-pushforward}(1)), and of the uniform measure on \(\Koz(n;\vek d)\). An asterisk
marks the product regime, where the first two agree (Proposition~\ref{prop-pushforward}(2)). At
\(\vek d=0^r\) the last column is the rank-one value \(\sum_{k=2}^n(k-1)(n-k+1)/(n^2(n+1))\) of
Proposition~\ref{prop-second-moments}, since \(\nu\) on \(r\Koz(n)\) is \(\mu\) on \(\Koz(n)\).}
\label{tab:rank-r-variance}
\end{table}

At rank one the effect reverses near \(n=14\) (Observation \ref{obs-concentration}); at higher rank
nothing is known for large \(n\).
\section{Methods and verification}
\label{sec-methods}
Every theorem above has a written proof. The two original groups of results (lattice points and
volume, centroids) were in addition re-derived line by line in an adversarial review by a second,
separately prompted instance of the language model, and checked by exact computation with code
written for that review, without access to the code used for the paper. The later results
(Theorem \ref{thm-centroid-dilation}, Kupreyeva's limit, and Theorem \ref{thm-codegree} with its
companions) were each reviewed in the same way, line by line, and checked by computation as listed
below. None of these reviews was by a human colleague; the author's own check is described in the
disclosure section. The code, the run logs and the Lean library are in the public repository
\url{https://gitlab.liu.se/jansn19/exterior-convex}; below, file names refer to its directory
\texttt{verification-code/oneoff/}, and the run logs are under \texttt{verification-code/runs/}. A frozen snapshot of the code, the run logs and the Lean library is archived on Zenodo
\autocite{SnellmanZenodoSeries} under the DOI
\href{https://doi.org/10.5281/zenodo.24288404}{10.5281/zenodo.24288404}, which resolves to its latest version.

\begingroup\sloppy

\begin{itemize}
\item \textbf{Lattice points, volume and the staircase tiling.} Lattice-point counts by Normaliz
\autocite{BrunsIchim2010}, by direct sumsets and by chained tuples; normalised volumes by Normaliz
against Theorem \ref{thm-staircase} (\(97\) cases); the lattice-path reading of staircase indices
(Example \ref{ex-staircase-paths}) for \(n\le7\), \(r\le4\) and all gaps at most \(n-1\), and its
identification at \(\vek d=\vek 0\) with the staircase triangulation (Remark
\ref{rem-staircase-cayley}) for \(n\le5\), \(r\le4\) (\texttt{staircase\_lattice\_paths.py}).
\item \textbf{Ehrhart series and \(h^{\ast}\).} Ehrhart series by Normaliz against Corollary \ref{cor-transfer}
wherever \(\NV\le3\cdot10^7\), and for four rank-two polytopes, \(\Koz(6;0,4)\) among them, also by
Barvinok's algorithm, in the \texttt{barvinok} library \autocite{Verdoolaege2007}; \(h^*\) by the ascent
formula, by the \(\vek s\)-Eulerian sum (\(n\le6\)) and by Normaliz (\(n\le7\)), the \(n=8\) vector by
the ascent formula alone. Script: \texttt{kozlov\_volume\_ehrhart.py}, review
\texttt{review\_volume\_ehrhart\_check.py}. Real-rootedness of the \(h^*\)-polynomials in
Observation \ref{obs-hstar-unimodal} by exact real-root counting with Sturm sequences, a count equal
to the degree meaning that all roots are real (\texttt{hstar\_realrooted\_check.py}).
\item \textbf{Codegree, quadratic coefficient, deficiency.}
Theorem \ref{thm-codegree}, Corollary \ref{cor-gorenstein} and Proposition \ref{prop-quadratic-coeff} against a
separately written transfer-matrix computation for every \(\vek d\) with every gap at most \(n\),
\(3\le n\le7\) and rank up to \(5\), \(4\), \(3\) (\(1854\) cases), and against the stored \(h^*\)-vectors; the
codegree and the interior lattice points also from convex hulls computed in SageMath for \(n=3,4\),
and the interior points of Lemma \ref{lem-interior-points} by membership tests
(\texttt{hstar\_codegree\_t2\_check.py}, \texttt{hstar\_codegree\_sage\_crosscheck.sage},
\texttt{interior\_points\_direct\_check.py}). The values in the proof of Corollary \ref{cor-deficiency}: sumsets of
\(S_n\) listed by the Kruskal--Katona theorem, and integer LYM chains, \(n=4,5\), \(r\le5\)
(\texttt{zero\_shift\_deficiency.py}); at \(n=4\) also by enumerating every complex on \([4]\).
\item \textbf{Centroids and moments.} Exact dynamic programming over the chain sets (\(n\le26\)), sumsets,
triangulation integrals in SageMath, Normaliz integrals, Kupreyeva's tables. Scripts:
\texttt{centroids\_moments.py}, \texttt{centroids\_moments.sage}, \texttt{profile\_large\_n.py}, review
\texttt{review\_centroids\_check.py}. The tables of Section \ref{sec-moments} are typeset from a rerun of
\texttt{profile\_large\_n.py} (\texttt{moments\_tables.py}); Table \ref{tab:rank-r-variance},
Observation \ref{obs-pushforward-fails} and the rank-one columns for \(n=4,6,8\) were recomputed
without that code, from \(S_n\) built out of all complexes on \([n]\), plain sumsets and a pulling
triangulation of \(\Koz(n;\vek d)\) (\texttt{moments\_tables\_check.sage}). For Theorem \ref{thm-centroid-dilation}: the inequalities of its proof,
the polynomial degrees and leading coefficients, and the constants, by exact enumeration for
\(n\le6\) (\texttt{centroid\_dilation\_check.py}); the polynomials \(L\), \(\Sigma_k\), \(B_k\) of part (3), by
interpolation at \(t=1,\dots,n+1\) and a check at \(t=n+2\) (\texttt{centroid\_dilation\_Rk\_polys.py}). For
Proposition \ref{prop-dilates}: the same interpolation for \(n=3,4,5\), checked at \(t=n+2,n+3\)
(\texttt{lattice\_centroid\_dilates\_exact.py}).
\item \textbf{Kupreyeva's limit.} Exact counts of \(S_n^{(p)}\) by Linusson's recursion and, independently, by
sums of pseudo-powers; lattice points by floor sums; both checked against brute-force sumsets for
\(p\le7\). Fill rates for \(n=3\) (\(p\le30000\)), \(n=4\) (\(p\le2000\)) and, at rank one,
\(n=5\) (\(p\le120\)); the limits in exact rational arithmetic; the volumes of the Minkowski sums
\(R_n+\dots+R_n\) for \(n\le5\), \(r\le3\) also by numerical quadrature. Scripts:
\texttt{kupreyeva\_limit\_rank\_r.py}, \texttt{kupreyeva\_qn\_volume\_ratio.py}.
\item \textbf{Figures.} \texttt{ff\_pictures.py}; Figure \ref{fig:mu-rank-r} by \texttt{fig\_mu\_rank\_r.py}.
\end{itemize}
\endgroup
\subsection{Formal verification}
\label{sec-lean}
Much of this paper is also proved in the proof assistant Lean 4, on top of its mathematical
library Mathlib \autocite{deMouraUllrich2021,Mathlib2020}. The Lean files are in the directory
\texttt{lean/} of the public repository \url{https://gitlab.liu.se/jansn19/exterior-convex} and of its
snapshot on Zenodo \autocite{SnellmanZenodoSeries}, with the code and the run logs described above. A
mark in the margin beside a result, or beside one of its parts, records how far the formal proof
goes:
\begin{itemize}
\item \leanfull{} proved in Lean, with no \texttt{sorry}, from the three axioms that Mathlib uses throughout
(propositional extensionality, choice and the soundness of quotients) and nothing else;
\item \leancond{} proved in the same way, assuming the surjectivity half of Part II's chaining
reformulation \autocite[Theorems 23 and 41, with Lemma 40]{SnellmanArticle2}, which enters as an
explicit hypothesis of the formal statement, not as an axiom;
\item \leanpartial{} part of the statement proved as for \leanfull{}: a special case or the
combinatorial core;
\item \leannone{} not formalised.
\end{itemize}
The file \texttt{lean/README.md}, in the repository and in the snapshot,
pairs each result with its formal statement, says what a partial mark covers, and records where a
formal definition differs in form from the printed one (the normalised volume, the uniform measure
on a dilate).
\section{Open problems}
\label{sec-open}
\begin{enumerate}
\item \textbf{Ehrhart theory at rank \(r\).} Unimodality and real-rootedness of \(h^*_{\Koz(n;\vek d)}\)
(Observation \ref{obs-hstar-unimodal}); Conjecture \ref{conj-nv-monotone}; and the signs of the Ehrhart
coefficients beyond \(t^2\), which Propositions \ref{prop-linear-coeff} and \ref{prop-quadratic-coeff}
reduce to ranks one and two.
\item \textbf{Fill rates.} Prove that \(\lvert S_n\rvert/\lvert\mathcal L_n\rvert\to0\), find its rate, and prove
that the Perron-root ratio of Observation \ref{obs-fill} is below \(1\) for \(n\ge4\) and every gap
\(\delta\ge1\). At equal shifts, find the polynomial \(L_{\Koz(n)}(r)-a(r)\) for general \(n\):
it vanishes for \(n\le3\), and at \(n=4,5\) (Corollary \ref{cor-deficiency}) it is divisible by
\(r(r+1)\). Is it always?
\item \textbf{Kupreyeva's limit beyond the first order.} Theorems \ref{thm-kupreyeva-limit} and
\ref{thm-kupreyeva-rank-r} give the limit of \(\operatorname{fill}^{(p)}(n;\vek d)\) as
\(p\to\infty\). In every computed case (Example \ref{ex-kupreyeva-values}),
\(p\bigl(\operatorname{fill}^{(p)}(n;\vek d)-\lim\bigr)\) appears to converge, to about
\(-4/5\) at \(n=3\), \(r=1\) and \(-1.66\) at \(n=4\), \(r=1\). Find the second-order term,
which should come from the lattice points near the Kruskal--Katona boundary; at rank one the
error is \(O(1/p)\) (Remark \ref{rem-kupreyeva-linusson}), and at rank \(r\) its order is not proved.
And relate \(q_n\) to the fill rate at \(p=n\), where open problem 2 asks for the
rate.
\item \textbf{Centroids.} The analogue of Theorem \ref{thm-centroid-dilation} at a constant positive gap, where
the fill rate tends to \(0\) and its argument fails; the limiting LYM profile of the centroid of
\(S_n\) as \(n\to\infty\), and the constant in Observation \ref{obs-concentration}; a rank-\(r\) analogue
of the involution \(\beta\), which would have to come from a duality of modules, since
\(\Koz(n;\vek d)\) is not a simplex.
\end{enumerate}

\appendix
\section{Proof of the canonical chained decomposition}
\label{sec-app-chained}
This appendix proves Theorem \ref{thm-chained}. The notation is that of Section \ref{sec-lattice}: an
integral leg vector \(\vek c\) with \((c_2,\dots,c_n)\) log-concave where needed, a shift vector
\(\vek d\) with gaps \(\delta_i\) \eqref{eq:gaps}, the LYM coordinates \(\mu_k=p_k/c_k\) \eqref{eq:mu} of a point \(\vek p\) of
\(\rasimp(\vek c)\), its statistics \(\fl\) and \(\tp\) of \eqref{eq:fl} and \eqref{eq:tp}, and the
water-filling quantities \(Q_i(k)\) of \eqref{eq:water}. A remark on what is new closes it.

\begin{proof}[Proof of Theorem~\ref{thm-chained}]
\emph{Uniqueness and the formula.} Let the tuple be chained and fix a global degree \(k\). If summand
\(i\) is non-zero in degree \(k\), then \(k-d_i\le\tp(\vek p_i)\), and chaining gives
\(k\le\fl(\vek p_l)+d_l\) for every \(l>i\), so every later summand is full in degree \(k\) (or
has capacity \(c_{k-d_l}=0\) there). Hence the summands, read in the order of \(i\), are zero, then
at most one partial, then full, and the amount carried by summands \(1,\dots,i\) in degree \(k\) is
\(Q_i(k)\). The possible values of \(x_k\) are the partial capacity sums
\(C_{\ge a}=\sum_{l\ge a}c_{k-d_l}\) and the open intervals between consecutive ones, and these are
pairwise disjoint because every \(c_s\) present is positive. So \(x_k\) determines each summand's
coordinate in degree \(k\), and every own coordinate of every summand occurs in some global degree.

\emph{Existence.} The set of representations of \(x\), tuples \((\mu^{(1)},\dots,\mu^{(r)})\) of LYM
chains \(1=\mu_1\ge\mu_2\ge\dots\ge\mu_n\ge0\) (\(\mu_k=(\vek p)_k/c_k\) as in \eqref{eq:mu}, cf. Proposition
\ref{prop-rank-one-H}) with \(\sum_i\shift^{d_i}(\vek c\cdot\mu^{(i)})=x\), the product taken
entrywise, is a non-empty compact polytope. On it the linear functional
\begin{equation}\label{eq:Psi}
  \Psi=\sum_{i=1}^r i\,\bigl\lVert\vek p_i\bigr\rVert_1,\qquad
  \bigl\lVert\vek p_i\bigr\rVert_1=\sum_{j}c_j\mu^{(i)}_j ,
\end{equation}
attains its maximum. Take a maximiser and suppose it is not chained, at some \(i\): with
\(u=\fl(\vek p_{i+1})\), \(T=\tp(\vek p_i)\) and \(\delta=d_{i+1}-d_i\), we have \(T>u+\delta\).
Then \(u<T-\delta\le n\), so \(\mu^{(i+1)}_u=1>\mu^{(i+1)}_{u+1}\), and \(\mu^{(i)}_T>0=\mu^{(i)}_{T+1}\)
(with \(\mu_{n+1}=0\)). Put \(j_0=u+\delta+1\le T\). For small \(\varepsilon>0\) (two strict
inequalities, at the two ends below, fix how small) replace
\begin{equation}\label{eq:exchange}
  \mu^{(i)}_j\mapsto\mu^{(i)}_j-\varepsilon,\qquad
  \mu^{(i+1)}_{j-\delta}\mapsto\mu^{(i+1)}_{j-\delta}+\varepsilon\,c_j/c_{j-\delta}\qquad(j_0\le j\le T).
\end{equation}
The global coordinate \(j+d_i=(j-\delta)+d_{i+1}\) changes by \(-\varepsilon c_j+\varepsilon c_j=0\),
so \(x\) is unchanged. Summand \(i\) remains a chain: it is lowered by the same amount on the whole
interval \([j_0,T]\), so only the two ends matter; at the left end lowering only helps, at the right
end \(\mu^{(i)}_T-\varepsilon\ge0=\mu^{(i)}_{T+1}\), and \(\mu^{(i)}_1=1\) is untouched since
\(j_0\ge2\). Summand \(i+1\) is raised on its own indices \(m\in[u+1,T-\delta]\subseteq[2,n]\) by
\(\eta_m=\varepsilon c_{m+\delta}/c_m\). Inside, for \(u+1\le m<T-\delta\), the chain needs
\(\eta_m\ge\eta_{m+1}\), that is \(c_{m+\delta}/c_m\ge c_{m+1+\delta}/c_{m+1}\), a product of
\(\delta\) comparisons \(c_{q+1}/c_q\ge c_{q+2}/c_{q+1}\) for \(m\le q\le m+\delta-1\), so
\(3\le q+1\le m+\delta\le T-1\le n-1\), which is the log-concavity of \((c_2,\dots,c_n)\). At the
right end raising only helps, and at the left end \(1=\mu_u\ge\mu_{u+1}+\eta_{u+1}\) for small
\(\varepsilon\), because \(\mu_{u+1}<1\); by the monotonicity just shown, every raised entry is
then at most \(1\). The exchange moves the mass
\(\varepsilon\sum_{j=j_0}^{T}c_j>0\) from summand \(i\) to summand \(i+1\):
\(\lVert\vek p_i\rVert_1\) falls by this amount and \(\lVert\vek p_{i+1}\rVert_1\) rises by
\(\sum_m\eta_mc_m\), the same amount. So \(\Psi\) rises by it, contradicting maximality.

\emph{Integrality.} Each coordinate of each part is \(c_s\), \(0\), or \(x_k-C_{\ge a+1}\), an integer.
\end{proof}

\begin{remark}
(\emph{What is new.}) At \(\vek d=(0,\dots,0)\) the lattice-point version of Theorem \ref{thm-chained},
existence and uniqueness, is Brändén and Solus's decomposition theorem for lecture hall polytopes
\autocite[Theorem 2.1]{BrandenSolus2020} (Theorem 1 in the arXiv version): their relation
\(\lambda\trianglelefteq\gamma\), ``\(\lambda\le\gamma\), and \(\gamma_i=s_i\) wherever
\(\lambda_i\ne0\)'', is chaining with \(\delta=0\), and their existence proof is the decomposition
\(\lambda=(\lambda\wedge s)+((\lambda-s)\vee0)\). For \(f\)-vectors and all \(\vek d\), uniqueness and the
water-filling formula are \autocite[Lemma 24]{SnellmanArticle2}, and existence of a chained tuple of \emph{\(f\)-vectors}
for every ff-vector is the content of Part II's chaining reformulation, proved there from Amata and
Crupi's lex submodules. What is new here is positive gaps for arbitrary points, and that is where
log-concavity enters: for \(\delta=0\) the inside-block condition is vacuous. The hypothesis is
needed. For \(\vek c=(1,2,1,3,1,3)\) and \(\vek d=(0,1)\) the polytope has \(101\) lattice points,
the sumset of lattice points only \(70\), and existence fails.
\label{rem-chained-credit}
\end{remark}
\section{Proof of Kupreyeva's limit}
\label{sec-app-kupreyeva}
The proofs in this appendix were found by the model; see the disclosure section. Throughout,
\(n\ge2\), \(p\ge2n+2\), \(\mu_k=x_k/\binom pk\) are the LYM coordinates \eqref{eq:mu}, and
\begin{equation}\label{eq:gk}
  g_k(t)=t^{(k+1)/k},\qquad t\ge0 ,
\end{equation}
so that the upper boundary of \(R_n\) in \eqref{eq:Rn} is \(\mu_{k+1}=g_k(\mu_k)\).

\begin{lemma}
\leanmark{\leannone}Let \(\Omega\subseteq\mathbb R^D\) be bounded, with boundary of Lebesgue measure zero, and for
\(\vek h\in(0,\infty)^D\) let \(\Gamma_{\vek h}=h_1\mathbb Z\times\dots\times h_D\mathbb Z\). Then
\(h_1\cdots h_D\cdot\#(\Gamma_{\vek h}\cap\Omega)\to\vol\Omega\) as \(\max_jh_j\to0\).
\label{lem-anisotropic-count}
\end{lemma}

\begin{proof}
Put \(\eta=\max_jh_j\) and \(Q(y)=y+\prod_j[0,h_j)\) for \(y\in\Gamma_{\vek h}\). These boxes
partition \(\mathbb R^D\). Each has volume \(h_1\cdots h_D\) and contains \(y\), and any two of its
points differ by less than \(\eta\) in every coordinate. Let \(\overline\Omega_\eta\) be the set of
points at sup-distance at most \(\eta\) from the closure \(\overline\Omega\), and \(\Omega^\circ_{-\eta}\)
the set of points whose closed sup-ball of radius \(\eta\) lies in the interior \(\Omega^\circ\). If
\(y\in\Gamma_{\vek h}\cap\Omega\), then \(Q(y)\subseteq\overline\Omega_\eta\). If
\(z\in\Omega^\circ_{-\eta}\cap Q(y)\), then \(Q(y)\subseteq\Omega^\circ\), so \(y\in\Omega\). Hence
\[
  \vol\Omega^\circ_{-\eta}\le h_1\cdots h_D\cdot\#(\Gamma_{\vek h}\cap\Omega)\le\vol\overline\Omega_\eta .
\]
As \(\eta\to0\), \(\overline\Omega_\eta\) decreases to the compact set \(\overline\Omega\) and
\(\Omega^\circ_{-\eta}\) increases to the open set \(\Omega^\circ\). So the bounds tend to
\(\vol\overline\Omega\) and \(\vol\Omega^\circ\), both equal to \(\vol\Omega\) because \(\partial\Omega\) is null.
\end{proof}

\begin{lemma}
\leanmark{\leannone}Let \(k\ge1\) and \(p\ge2k+2\), and put \(\eta_s=s^2/(p-s+1)\). For every integer
\(0\le f\le\binom pk\), with \(\mu=f/\binom pk\),
\begin{equation}\label{eq:epsk}
  \biggl\lvert\frac{f^{(k)}}{\binom p{k+1}}-\mu^{(k+1)/k}\biggr\rvert\le\varepsilon_k(p)
  :=\eta_{k+1}+\frac{(k-1)\,(k+1)!\,2^{k+1}}{p}+(k+1)\Bigl(\eta_k+\frac kp\Bigr)^{1/k},
\end{equation}
and \(\varepsilon_k(p)\to0\) as \(p\to\infty\).
\label{lem-pseudo-power}
\end{lemma}

\begin{proof}
If \(f=0\), both sides vanish. Otherwise let
\(f=\binom{a_k}k+\binom{a_{k-1}}{k-1}+\dots+\binom{a_t}t\) be the Macaulay representation
\eqref{eq:macaulay}, and put \(a=a_k\). Then \(\binom ak\le f<\binom{a+1}k\). Since \(f\le\binom pk\),
we have \(a\le p\), with \(a=p\) only if \(f=\binom pk\). Put \(\alpha=a/p\) and
\(\beta_s(b)=\binom bs/\binom ps\).

First, \(\lvert\beta_s(b)-(b/p)^s\rvert\le\eta_s\) for integers \(0\le b\le p\) and \(1\le s\le p\).
\begin{itemize}
\item If \(b\ge s-1\), then \(\beta_s(b)=\prod_{i<s}(b-i)/(p-i)\) and \((b/p)^s\) are products of \(s\)
numbers in \([0,1]\). Their difference is at most
\(\sum_{i<s}\lvert b/p-(b-i)/(p-i)\rvert=\sum_{i<s}i(p-b)/(p(p-i))\le s^2/(p-s+1)\).
\item If \(b\le s-2\), then \(\beta_s(b)=0\) and \((b/p)^s\le s/p\le\eta_s\).
\end{itemize}

Next, \(\mu\) lies between \(\beta_k(a)\) and \(\beta_k(a+1)\), or equals \(1=\alpha^k\) if \(a=p\).
Since \(((a+1)/p)^k\le\alpha^k+k/p\), this gives \(\lvert\mu-\alpha^k\rvert\le\eta_k+k/p\). Using
\(\lvert u^{1/k}-v^{1/k}\rvert\le\lvert u-v\rvert^{1/k}\) for \(u,v\ge0\), and the fact that
\(t\mapsto t^{k+1}\) is \((k+1)\)-Lipschitz on \([0,1]\),
\(\lvert\mu^{(k+1)/k}-\alpha^{k+1}\rvert\le(k+1)(\eta_k+k/p)^{1/k}\).

Finally, \(f^{(k)}=\binom a{k+1}+\sum_{j<k}\binom{a_j}{j+1}\). The at most \(k-1\) terms of the sum
are each at most \(p^k\), and \(\binom p{k+1}\ge(p/2)^{k+1}/(k+1)!\). So
\(\bigl\lvert f^{(k)}/\binom p{k+1}-\alpha^{k+1}\bigr\rvert
\le\eta_{k+1}+(k-1)(k+1)!\,2^{k+1}/p\).
\end{proof}

\begin{lemma}
\leanmark{\leanpartial}For \(2\le a\le b\) the region
\(R_{a,b}=\lbrace\mu\in\mathbb R^{\lbrace a,\dots,b\rbrace}:0\le\mu_a\le1,\ 0\le\mu_{k+1}\le g_k(\mu_k)
\ (a\le k<b)\rbrace\) lies in \(\lbrace1\ge\mu_a\ge\dots\ge\mu_b\ge0\rbrace\), and
\[
  \vol R_{a,b}=\rho(a,b)=\prod_{k=a}^{b-1}\frac{2k}{(b-k+1)(b+k)}
  =\frac{2^{b-a}\,(b-1)!\,(a+b-1)!}{(a-1)!\,(b-a+1)!\,(2b-1)!} .
\]
In particular \(R_n=R_{2,n}\) has \(\vol R_n=2^{n-2}(n+1)!/(2n-1)!\), and
\((n-1)!\vol R_n=(n+1)!/\bigl(2(2n-1)!!\bigr)\).
\label{lem-rho}
\end{lemma}

\begin{proof}
The inclusion holds because \(g_k(t)\le t\) on \([0,1]\). Put \(E_k=\sum_{j=k+1}^bj\). Let
\(F_b\equiv1\) and \(F_k(t)=\int_0^{g_k(t)}F_{k+1}(u)\,du\) for \(a\le k<b\). By downward induction,
\(F_k(t)=C_kt^{E_k/k}\) with \(C_b=1\) and \(C_k=C_{k+1}(k+1)/E_k\), because
\[
  \int_0^{t^{(k+1)/k}}C_{k+1}u^{E_{k+1}/(k+1)}\,du=C_{k+1}\frac{k+1}{E_{k+1}+k+1}\,t^{(E_{k+1}+k+1)/k}
\]
and \(E_{k+1}+k+1=E_k\). By Tonelli, \(\vol R_{a,b}=\int_0^1F_a=C_a\,a/(E_a+a)\), which is
\(\prod_{k=a}^bk/E_{k-1}\) with \(E_{a-1}=E_a+a\). Its last factor is \(b/E_{b-1}=1\), and
\(E_{k-1}=(b-k+1)(b+k)/2\). The closed forms follow from \(\prod_{k=a}^{b-1}(b-k+1)=(b-a+1)!\),
\(\prod_{k=a}^{b-1}(b+k)=(2b-1)!/(a+b-1)!\) and \((2n-1)!=2^{n-1}(n-1)!\,(2n-1)!!\).
\end{proof}

For \(\vek x\in\mathcal L_n^{(p)}\) call \((\tp(\vek x),\indeg(\vek x))\) its \textbf{class}. Since
\(x_1=p>0\) and \(\mu_1=1\), the classes are the pairs \((\tau,\iota)\) with \(1\le\tau\le n\) and
\(2\le\iota\le\tau+1\): the \(n(n+1)/2\) lattice points of a triangle. In the terms of the
discussion after Definition \ref{def-chained-points}, a point of class \((\tau,\iota)\) is full in its
own degrees below \(\iota\), partial in \(\iota,\dots,\tau\), and zero above \(\tau\). So the class is
the block of degrees in which the point is partial, and \([\iota-1,\tau]\) is its barycentric
support. For \(\iota=\tau+1\) the block is empty and the point is the vertex \(\tilde F_\tau\); the
class then records only the degree \(\iota\) at which the point drops from full to zero. Within a
class the partial coordinates \((\mu_\iota,\dots,\mu_\tau)\) range over a grid in a simplex, which is
what Lemma \ref{lem-anisotropic-count} counts (proof of Lemma \ref{lem-kup-classes}). The numbers of points of \(\mathcal L_n^{(p)}\) and of \(S_n^{(p)}\) in a
class, and their common scale, are
\begin{equation}\label{eq:classcounts}
\begin{gathered}
  M^X_p(\tau,\iota)=\#\bigl\lbrace\vek x\in X_n^{(p)}:\bigl(\tp(\vek x),\indeg(\vek x)\bigr)=(\tau,\iota)\bigr\rbrace
  \qquad(X=\mathcal L,S),\\
  P_{\iota,\tau}(p)=\prod_{k=\iota}^\tau\binom pk ;
\end{gathered}
\end{equation}
and we put \(\rho(\tau+1,\tau)=1\).

\begin{lemma}
\leanmark{\leannone}As \(p\to\infty\), \(M^{\mathcal L}_p(\tau,\iota)/P_{\iota,\tau}(p)\to1/(\tau-\iota+1)!\) and
\(M^S_p(\tau,\iota)/P_{\iota,\tau}(p)\to\rho(\iota,\tau)\).
\label{lem-kup-classes}
\end{lemma}

\begin{proof}
For \(\iota=\tau+1\) the class is the vertex \(\tilde F_\tau\), an \(f\)-vector. Let
\(\iota\le\tau\). A point of the class has \(\mu_k=1\) for \(k<\iota\) and \(\mu_k=0\) for
\(k>\tau\), and is determined by its block \(\mu^\flat=(\mu_\iota,\dots,\mu_\tau)\). The block ranges
over \(\prod_{k=\iota}^\tau\binom pk^{-1}\mathbb Z\) intersected with
\(\Delta^\circ=\lbrace1>\mu_\iota\ge\dots\ge\mu_\tau>0\rbrace\), a set of volume
\(1/(\tau-\iota+1)!\) with null boundary. The first claim is therefore
Lemma \ref{lem-anisotropic-count}.

By the Kruskal--Katona theorem \eqref{eq:kk-cond}, and since \(\binom pk^{(k)}=\binom p{k+1}\), a
point of the class is an \(f\)-vector if and only if \(x_{k+1}\le x_k^{(k)}\) for \(\iota\le k<\tau\). For \(s>0\) let
\(R^{\pm s}=\lbrace\mu^\flat\in\Delta^\circ:\mu_{k+1}\le g_k(\mu_k)\pm s\ (\iota\le k<\tau)\rbrace\).
Once \(\varepsilon_k(p)\le s\) for these \(k\), Lemma \ref{lem-pseudo-power} gives two inclusions:
\begin{itemize}
\item every \(f\)-vector of the class has \(\mu^\flat\in R^{+s}\);
\item every lattice point of the class with \(\mu^\flat\in R^{-s}\) is an \(f\)-vector.
\end{itemize}

The boundaries of \(R^{\pm s}\) lie in \(\partial\Delta^\circ\) and in the graphs
\(\mu_{k+1}=g_k(\mu_k)\pm s\), which are null. So Lemma \ref{lem-anisotropic-count} gives
\[
  \vol R^{-s}\le\liminf_{p\to\infty}\frac{M^S_p(\tau,\iota)}{P_{\iota,\tau}(p)}
  \le\limsup_{p\to\infty}\frac{M^S_p(\tau,\iota)}{P_{\iota,\tau}(p)}\le\vol R^{+s}.
\]
As \(s\downarrow0\), \(R^{+s}\) decreases to \(\Delta^\circ\cap R_{\iota,\tau}\). \(R^{-s}\)
increases to a set that differs from it by a subset of the graphs \(\mu_{k+1}=g_k(\mu_k)\). Both
volumes therefore tend to \(\rho(\iota,\tau)\) (Lemma \ref{lem-rho}).
\end{proof}

\begin{proof}[Proof of Theorem~\ref{thm-kupreyeva-limit}]
\(\lvert\mathcal L_n^{(p)}\rvert\) and \(\lvert S_n^{(p)}\rvert\) are the sums of
\(M^{\mathcal L}_p\) and \(M^S_p\) over the finitely many classes. For the class \((n,2)\),
Lemma \ref{lem-kup-classes} gives \(P_{2,n}(p)\bigl(1/(n-1)!+o(1)\bigr)\) and
\(P_{2,n}(p)\bigl(\rho(2,n)+o(1)\bigr)\). Every other class has
\(\sum_{k=\iota}^\tau k<\sum_{k=2}^nk\), and at most \(\prod_{k=\iota}^\tau(\binom pk+1)\) points,
which is \(o(P_{2,n}(p))\). Lemma \ref{lem-rho} gives \(\rho(2,n)=\vol R_n\) and the value of
\(q_n\).
\end{proof}

\begin{proof}[Second proof of Theorem~\ref{thm-kupreyeva-limit}]
This proof takes over from the first only the count of lattice points,
\(\lvert\mathcal L_n^{(p)}\rvert\sim\prod_{k=2}^n\binom pk/(n-1)!\), and counts the \(f\)-vectors with
Linusson's formula instead of by volume. In his
notation \(\lvert S_n^{(p)}\rvert=F^p(n)-F^{p-1}(n)=\sum_{k=1}^nE^p(n,k)\), and \(E^p(n,k)\) is a
polynomial in \(p\) of degree \(M-\binom{k+1}2\), where \(M=\binom{n+1}2\). So
the term \(k=1\) dominates, and its leading coefficient is
\begin{equation}\label{eq:linusson}
  \gamma_n=\frac1{(M-1)!}\prod_{i=1}^{n-2}\binom{M-\binom{i+1}2-1}{i}
\end{equation}
\autocite[Theorem 3.3, and Proposition 2.8 with its proof]{Linusson1999}. Put
\(\alpha_i=M-\binom{i+1}2\). Then \(M-\binom{i+1}2-1-i=\alpha_{i+1}\), so the product telescopes to
\[
  \gamma_n=\frac{1}{(n-1)!\,\prod_{i=1}^{n-1}\alpha_i\,\prod_{i=1}^{n-2}i!} .
\]
Since \(\prod_{i=1}^{n-1}\alpha_i=(n-1)!\,(2n)!/\bigl(2^{n-1}(n+1)!\bigr)\),
\[
  \gamma_n\,(n-1)!\prod_{k=2}^nk!=\frac{2^{n-1}\,n!\,(n+1)!}{(2n)!}=q_n .
\]
Since \(M-1=\sum_{k=2}^nk\) and \(\prod_{k=2}^n\binom pk\sim p^{M-1}/\prod_{k=2}^nk!\), this says
\(\lvert S_n^{(p)}\rvert\sim\gamma_n\prod_{k=2}^nk!\,\prod_{k=2}^n\binom pk\), with constant
\(\gamma_n\prod_{k=2}^nk!=q_n/(n-1)!=2^{n-2}(n+1)!/(2n-1)!\), the closed form in the theorem; and
\(q_n\) is the ratio of the leading coefficients of \(\lvert S_n^{(p)}\rvert\) and
\(\lvert\mathcal L_n^{(p)}\rvert\), which is \eqref{eq:qn}. That the constant is the volume of \(R_n\)
is Lemma \ref{lem-rho}, which this proof does not need.
\end{proof}

So Kupreyeva's observation was already implied by Linusson's leading coefficient, together with the
first-order count of lattice points. What the first proof, by volume, adds is the reason for the
value, and the extension to rank \(r\).

\begin{remark}
(\emph{The rate.}) Together, Linusson's polynomiality and the volume argument give the rate. \(\lvert S_n^{(p)}\rvert\)
is a polynomial in \(p\) \autocite{Linusson1999},
and \(\lvert\mathcal L_n^{(p)}\rvert\) is a quasi-polynomial in \(p\) for large \(p\), because the
vertices of \(\Koz^{(p)}(n)\) are polynomials in \(p\) \autocite[Theorem 1.4]{ChenLiSam2012}. By
Theorem \ref{thm-kupreyeva-limit} both have degree \(\sum_{k=2}^nk\) and constant leading coefficients in
ratio \(q_n\), so \(\lvert S_n^{(p)}\rvert-q_n\lvert\mathcal L_n^{(p)}\rvert\) has lower degree, and
\(\lvert S_n^{(p)}\rvert/\lvert\mathcal L_n^{(p)}\rvert=q_n+O(1/p)\).
\label{rem-kupreyeva-linusson}
\end{remark}

A \textbf{type} is a tuple \(((\tau_1,\iota_1),\dots,(\tau_r,\iota_r))\) of classes. It is \textbf{chained} if
\(\tau_{i-1}<\iota_i+\delta_i\) for \(2\le i\le r\), which is the condition \eqref{eq:chained} of
Definition \ref{def-chained-points} read on classes.

In global degrees, summand \(i\) is partial in \(\iota_i+d_i,\dots,\tau_i+d_i\), and a type is
chained exactly when these blocks follow one another along the line of global degrees:
\[
  \iota_1+d_1\;\le\;\tau_1+d_1+1\;\le\;\iota_2+d_2\;\le\;\tau_2+d_2+1\;\le\;\dots\;\le\;\iota_r+d_r\;\le\;\tau_r+d_r+1 ,
\]
the inequalities \(\iota_i\le\tau_i+1\) holding within a class and the others being the chaining
condition. So a chained type is a weakly increasing sequence of \(2r\) integers, the \(i\)-th pair in
\([d_i+2,d_i+n+1]\): disjoint blocks of partial degrees in the order of the summands, which is the
profile ``zero, then at most one partial, then full'' of the discussion after
Definition \ref{def-chained-points}, read degree by degree.
A staircase index \eqref{eq:staircase-index} is the case in which the blocks are contiguous and
cover \([2,N]\), with \(\iota_i+d_i=t_{i-1}+1\) and \(\tau_i+d_i=t_i\). When every gap is at most
\(n-1\), the types of largest weight in Lemma \ref{lem-kup-dominant} below are of this kind.

\begin{lemma}
\leanmark{\leannone}Let \(p\ge n\).
\begin{enumerate}
\item Theorem \ref{thm-chained} and Corollary \ref{cor-additivity} hold for \(\Koz^{(p)}(n;\vek d)\).
\item A lattice point of \(\Koz^{(p)}(n;\vek d)\) lies in \(\FF^{(p)}(n;\vek d)\) if and only if every
part of its canonical decomposition lies in \(S_n^{(p)}\). Padding the parts with zeros gives the
canonical decomposition of the padded point in \(\Koz(p;\vek d)\).
\item Both \(\#\bigl(\Koz^{(p)}(n;\vek d)\cap\mathbb Z^{N+1}\bigr)\) and
\(\lvert\FF^{(p)}(n;\vek d)\rvert\) are sums over the chained types, of
\(\prod_iM^{\mathcal L}_p(\tau_i,\iota_i)\) and \(\prod_iM^S_p(\tau_i,\iota_i)\) respectively.
\end{enumerate}
\label{lem-kup-truncation}
\end{lemma}

\begin{proof}
(1) The sequence \((\binom p2,\dots,\binom pn)\) is log-concave.

(2) ``If'' is trivial, and for \(n=p\) the statement is Corollary \ref{cor-membership}, so let \(n<p\).
Let \(\operatorname{pad}\) append zeros, from \(\mathbb R^{n+1}\) to \(\mathbb R^{p+1}\) and from
\(\mathbb R^{N+1}\) to \(\mathbb R^{p+d_r+1}\). It commutes with shifts and maps \(S_n^{(p)}\) into
\(S_p\). It maps a point of \(\Koz^{(p)}(n)\) to a point of \(\Koz(p)\) with the same \(\fl\) and
\(\tp\): the LYM chain is extended by zeros, and \(\mu_{n+1}=0<1\).

Let \(\vek x=\sum_i\shift^{d_i}\vek s_i\) with \(\vek s_i\in S_n^{(p)}\), and let
\((\vek p_1,\dots,\vek p_r)\) be its canonical decomposition. Then
\(\operatorname{pad}\vek x\in\FF(p;\vek d)\). The tuple \((\operatorname{pad}\vek p_i)_i\) is a
chained tuple of points of \(\Koz(p)\) with sum \(\operatorname{pad}\vek x\), so by
Theorem \ref{thm-chained}(2) it is the canonical decomposition of \(\operatorname{pad}\vek x\). By
Corollary \ref{cor-membership}, with \(p\) in place of \(n\), each \(\operatorname{pad}\vek p_i\) lies
in \(S_p\). It vanishes above degree \(n\), so \(\vek p_i\in S_n^{(p)}\).

(3) By (1), (2) and the uniqueness in Theorem \ref{thm-chained}, the two sets are in bijection with
the chained tuples of points of \(\mathcal L_n^{(p)}\) and of \(S_n^{(p)}\) respectively. Whether a
tuple is chained depends only on the classes of its parts. The second formula is also Part II's
\autocite[Theorem 60]{SnellmanArticle2}, obtained there from lex submodules.
\end{proof}

\begin{lemma}
\leanmark{\leannone}For a type \(t\) put \(\omega(t)=\sum_i\sum_{k=\iota_i}^{\tau_i}k\). Among the chained types, \(\omega(t)\) is
largest exactly for the following ones. For each \(j\), the \(m_j\) summands with shift \(D_j\) \eqref{eq:aj}
have, in their order, \((\fl,\tp)=(s_{l-1},s_l)\) for \(l=1,\dots,m_j\), where
\(a_j-1=s_0\le s_1\le\dots\le s_{m_j}=n\). The multiset of the degrees \(k\) with
\(\iota_i\le k\le\tau_i\) is the same for all of them.
\label{lem-kup-dominant}
\end{lemma}

\begin{proof}
Let \(B_i=[\iota_i+d_i,\tau_i+d_i]\) be the global degrees in which summand \(i\) is partial, as
above. Chaining \eqref{eq:chained} gives \(\tau_{i-1}+d_{i-1}\le\iota_i-1+d_i\le\tau_i+d_i\), so the \(B_i\) are
disjoint and increasing. For \(g\in B_i\), \(2\le g-d_i\le n\). So \(g\) lies in
\(U=\bigcup_i[d_i+2,d_i+n]\), and \(d_i\ge d_*(g):=\min\lbrace d_l:d_l\ge g-n\rbrace\). Hence
\[
  \omega(t)=\sum_i\sum_{g\in B_i}(g-d_i)\le\sum_{g\in U}\bigl(g-d_*(g)\bigr),
\]
with equality if and only if every \(g\in U\) lies in a block \(B_i\) with \(d_i=d_*(g)\). The
\(g\in U\) with \(d_*(g)=D_j\) form the interval \([D_j+a_j,D_j+n]\). So equality means that the
own-degree blocks \((\fl,\tp]\) of the summands with shift \(D_j\) partition \((a_j-1,n]\).

Write \((f_l,t_l)\) for \((\fl,\tp)\) of the \(l\)-th summand with shift \(D_j\). Chaining at gap
\(0\) and \(\fl\le\tp\) give \(f_1\le t_1\le f_2\le\dots\le t_{m_j}\). The partition then forces
\(f_1\le a_j-1\) and \(t_{m_j}=n\). Conversely \(f_1\ge a_j-1\): \(f_1\ge1\) always, and for
\(j\ge2\), chaining with the last summand of shift \(D_{j-1}\), whose \(\tp\) is \(n\), gives
\(n+D_{j-1}\le f_1+D_j\). So \(f_1=a_j-1\), and the partition forces \(f_l=t_{l-1}\).

Every such type is chained, the inequality between groups being \(n<a_j+D_j-D_{j-1}\). The
multiset is \(\lbrace g-d_*(g):g\in U\rbrace\).
\end{proof}

\begin{proof}[Proof of Theorem~\ref{thm-kupreyeva-rank-r}]
Let \(\omega^*\) be the largest \(\omega(t)\) over chained types, and \(T^*\) the set of types attaining it.
For \(p\ge2n\) the counts \eqref{eq:classcounts} satisfy
\(M^X_p(\tau,\iota)\le\prod_{k=\iota}^\tau(\binom pk+1)\le2^nP_{\iota,\tau}(p)\), and
\(P_{\iota,\tau}(p)/p^{\,\iota+\dots+\tau}\to w(\iota,\tau)=\prod_{k=\iota}^\tau1/k!\). By
Lemmas \ref{lem-kup-truncation}(3) and \ref{lem-kup-classes}, \(p^{-\omega^*}\) times the number of lattice
points tends to
\[
  \sum_{t\in T^*}\prod_iw(\iota_i,\tau_i)/(\tau_i-\iota_i+1)!
\]
and \(p^{-\omega^*}\lvert\FF^{(p)}(n;\vek d)\rvert\) tends to
\[
  \sum_{t\in T^*}\prod_iw(\iota_i,\tau_i)\,\rho(\iota_i,\tau_i);
\]
the types outside \(T^*\) contribute \(O(p^{\omega^*-1})\). By Lemma \ref{lem-kup-dominant},
\(\prod_iw(\iota_i,\tau_i)\) is the same for all \(t\in T^*\), and \(T^*\) is a product over \(j\)
of the sets of sequences \(s_0\le\dots\le s_{m_j}\). So the limit of the fill rate is
\[
  \prod_j\frac{\sum_{\vek s}\prod_l\rho(s_{l-1}+1,s_l)}{\sum_{\vek s}\prod_l1/(s_l-s_{l-1})!} .
\]
The denominator is \(m_j^{\ell_j}/\ell_j!\) by the multinomial theorem.

(1) All \(m_j=1\), and \(\ell_j!\,\rho(n-\ell_j+1,n)=q(n,\ell_j)\) with \(\ell_j=\min(\delta_j,n-1)\)
and \(\ell_1=n-1\), by Lemma \ref{lem-rho}.

(2) Here \(g=1\), \(a_1=2\) and \(\vol(r\Delta_n)=r^{n-1}/(n-1)!\). It remains to show
\(\vol(R_n+\dots+R_n)=\Phi_n(2;r)\). Let \(G_k(v)=\lfloor v\rfloor+g_k(v-\lfloor v\rfloor)\) for
\(0\le v<r\), \(G_k(r)=r\), and
\(Z=\lbrace\nu:0\le\nu_2\le r,\ 0\le\nu_{k+1}\le G_k(\nu_k)\ (2\le k<n)\rbrace\).
\begin{itemize}
\item If \(\nu=\sum_j\mu^{(j)}\) with \(\mu^{(j)}\in R_n\), then
\(\nu_{k+1}\le\sum_jg_k(\mu^{(j)}_k)\le G_k(\nu_k)\). The last step holds because
\(\sum_jg_k(t_j)\) is convex on \(\lbrace t\in[0,1]^r:\sum_jt_j=\nu_k\rbrace\), hence maximal at a
vertex, where all \(t_j\) but at most one are \(0\) or \(1\).
\item Conversely, for \(\nu\in Z\) the parts
\(\mu^{(i)}_k=\min\bigl(1,\max(0,\nu_k-r+i)\bigr)\) of its canonical decomposition sum to \(\nu\)
and lie in \(R_n\): when \(r-i<\nu_k<r-i+1\), the bound \(\nu_{k+1}\le G_k(\nu_k)\) says exactly
\(\mu^{(i)}_{k+1}\le g_k(\mu^{(i)}_k)\); when \(\nu_k\ge r-i+1\), \(\mu^{(i)}_k=1\); and when
\(\nu_k\le r-i\), \(\mu^{(i)}_{k+1}=0\), because \(\nu_{k+1}\le G_k(\nu_k)\le\nu_k\).
\item So \(Z=R_n+\dots+R_n\). Now \(Z\subseteq r\Delta_n\), and on the interior of the cell
\(C_{\vek s}\) of Theorem \ref{thm-staircase} (all-ones legs, \(\vek d=0^r\)), \(Z\) is a translate
of \(\prod_iR^\circ_{s_{i-1}+1,s_i}\). Inside a block the condition reads
\(\mu_{k+1}\le g_k(\mu_k)\), and across blocks it holds automatically.
\item Since the cells tile \(r\Delta_n\), \(\vol Z=\Phi_n(2;r)\). For \(n=3\),
\(\Phi_3(2;r)=2r/5+\binom r2\).
\end{itemize}
\end{proof}
\section*{Acknowledgements}
\label{sec:org6182b05}
This paper owes its lattice-point question to Aliaksandra Kupreyeva's Bachelor's thesis
\autocite{Kupreyeva2019}, written under the author's supervision at Linköping University, where the fraction of
lattice points of the Kozlov simplex that are \(f\)-vectors, its volume and its centroid were first
studied. Luca Amata and Marilena Crupi's Macaulay2 package \autocite{AmataCrupi2021} was used
throughout.

Claude (Anthropic) assisted substantially with this paper; what that amounted to is set out in the
disclosure section below.

The computations used SageMath \autocite{SageMath}, Normaliz \autocite{BrunsIchim2010}, Macaulay2
\autocite{Macaulay2} and the \texttt{barvinok} library \autocite{Verdoolaege2007}; the formal verification used
Lean 4 \autocite{deMouraUllrich2021} and Mathlib \autocite{Mathlib2020}.
\section*{Disclosure of AI assistance}
\label{sec:org260b5cc}
This manuscript was produced with substantial assistance from a large language model, Claude
(Anthropic), and the extent of it is set out here rather than left to be inferred. All of it was
done at the author's request and under his direction. Where a result is said below to have been
reviewed, the review was an adversarial one by a second, separately prompted instance of the model
(Section \ref{sec-methods}); none was by a human colleague.

\begin{itemize}
\item \textbf{Sections \ref{sec-lattice}, \ref{sec-volume}, \ref{sec-ehrhart}.} The lecture hall identification, the
canonical chained decomposition at positive gaps (Theorem \ref{thm-chained}) and its corollaries, the
staircase tiling (Theorem \ref{thm-staircase}), the transfer-matrix Ehrhart formulas and the comparison
with \(\lvert\FF\rvert\), were found and proved by the model, and all tables were computed by it;
the literature on lecture hall polytopes was located by it. So were the lattice-path reading of
staircase indices and its identification at equal shifts with the staircase triangulation
(Example \ref{ex-staircase-paths}, Remark \ref{rem-staircase-cayley}, Figure \ref{fig:staircase-paths}),
and the sources for the latter.
\item \textbf{Sections \ref{sec-centroid} and \ref{sec-moments}.} The involution \(\beta\), Theorems
\ref{thm-beta-lattice}, \ref{thm-alexander-core} and \ref{thm-n3-equal-shifts}, and all centroid and moment
computations, were found and proved, or computed, by the model, and reviewed. So was the proof of
Theorem \ref{thm-centroid-dilation}, in November 2026; it had been stated as a conjecture in earlier
drafts. So were Remark \ref{rem-weak-convergence} and the exact centroids of
Proposition \ref{prop-dilates}, the latter stated as an observation in earlier drafts.
\item \textbf{Corollaries \ref{cor-zero-shift} and \ref{cor-deficiency}.} The statement about \(a(r)\) is Part II's; the
lattice-point count, the comparison and the closed forms at \(n=4,5\) are new here, and were
produced by the model.
\item \textbf{Kupreyeva's limit (Section \ref{sec-kupreyeva}, Appendix \ref{sec-app-kupreyeva}).} The
identification of Kupreyeva's \(q_n\) as a volume ratio, the proof of her observation
(Theorem \ref{thm-kupreyeva-limit}), the second proof from Linusson's leading coefficient and the rate
it gives (Appendix \ref{sec-app-kupreyeva}, Remark \ref{rem-kupreyeva-linusson}), and the limit at every rank
(Theorem \ref{thm-kupreyeva-rank-r}) with its proof, were found by the model in November 2026, and
reviewed.
\item \textbf{Theorem \ref{thm-codegree}, Corollary \ref{cor-gorenstein} and Proposition \ref{prop-quadratic-coeff}.}
Proved by the model in November 2026, from sketches in its own working notes; in the same check it
found that the codegree criterion conjectured in an earlier draft fails when a gap of at least
\(n-1\) follows an unrepeated value (Example \ref{ex-codegree-product}). Reviewed. The direct proof of
Lemma \ref{lem-interior-points}(2) was written by the model later, when the lemma was moved ahead of
Proposition \ref{prop-hstar-props}, and has been checked by computation only.
\item \textbf{Figures, code, formalisation.} All figures, all SageMath, Python and Normaliz code, and the Lean
formalisation described in Section \ref{sec-methods} were written by the model.
\item \textbf{Drafting.} This draft was written by the model from the working notes, with the author directing
the scope and the organisation of the series.
\end{itemize}

Every numerical claim in the paper was checked by exact computation along a route independent of the
formula being tested, in most cases by two independently written programs, and the reference list has
been checked against primary sources, except for the paper of Savage and Schuster, whose theorem is
cited in the form quoted by two other sources. The author has verified the results to the best of his
ability, made the final edits, and assumes full responsibility for any errors, mistakes or gaps that
remain.
\section*{Funding}
\label{sec:org49bd22b}
This research received no external funding. It was carried out as part of the author's regular duties
at Linköping University, within the fraction of that employment allocated to research.
\section*{Conflicts of interest}
\label{sec:org6ed2b0c}
The author declares no financial or non-financial conflicts of interest. The nature and extent of
large-language-model assistance in producing this manuscript is disclosed in full in the section
above.
\newcommand{\etalchar}[1]{$^{#1}$}

\end{document}